\documentclass[11pt,reqno,a4paper]{amsart}
\usepackage{color}
\usepackage{hyperref}
\usepackage{scalerel}
\usepackage{amsmath,amssymb}
\usepackage{mathtools}
\usepackage[latin1]{inputenc}
\usepackage[active]{srcltx}
\usepackage{graphicx}
\usepackage{exscale,relsize}
\usepackage{textgreek}
\usepackage{epsfig,graphics}
\usepackage{psfrag}
\usepackage{caption}
\usepackage{subcaption}
\usepackage[toc,page]{appendix}
\usepackage{bigints}

\def\N{\mathbb{N}}
\def\R{\mathbb{R}}
\def\m1{{I\!\!M}}

\newcommand{\parenthesis}[1]{\left(#1\right)} % round bracket
\newcommand{\braces}[1]{\left\{#1\right\}} % curly bracket
\def\ee{\`e}

\def\aa{\`a}

\renewcommand{\to}{\rightarrow}
\newcommand{\pa}{\partial}

\newcommand{\ino}{\int_{\Omega}}

\newcommand{\dx}{\, dx}

\newcommand{\opep}{\mbox{\rm o}_{\epsilon}(1)}

\newcommand{\rife}[1]{(\ref{#1})}
\newcommand{\ov}[1]{\overline{#1}}
\newcommand{\un}[1]{\underline{#1}}

\newcommand{\sscp}{\scriptscriptstyle}
\newcommand{\dsp}{\displaystyle}

\renewcommand{\dfrac}{\displaystyle\frac}
\newcommand{\finedim}{\hspace{\fill}$\square$}
\newcommand{\intbar}{\mathop{\int\makebox(-15.5,0){\rule[6pt]{.7em}{0.3pt}}\kern-6pt}\nolimits}

\renewcommand{\dx}{\, dx}

\newcommand{\ii}{\infty}

\newcommand{\eps}{\varepsilon}

\newcommand{\al}{\alpha}

\newcommand{\sg}{\sigma}

\newcommand{\om}{\Omega}
\newcommand{\lm}{\lambda}

\newcommand{\rl}{\mbox{\Large \textrho}_{\!\sscp \lm}}

\newcommand{\thl}{\theta_{\sscp \lm}}

\renewcommand{\rho}{\mbox{\Large \textrho}}
\newcommand{\rh}{\mbox{\Large \textrho}}

\newcommand{\pl}{\psi_{\sscp \lm}}
\newcommand{\xil}{\ul}
\newcommand{\ul}{u_{\sscp \lm}}

\newcommand{\ssl}{\sscp \lm}

\newcommand{\all}{\al_{\ssl}}

\newcommand{\el}{E_{\ssl}}

\newcommand{\lmn}{\lm_n}
\newcommand{\an}{\al_n}
\newcommand{\pn}{\psi_n}
\newcommand{\vn}{v_n}

\DeclareMathOperator{\dist}{dist}

\DeclareMathOperator{\Radius}{Radius}
\newcommand{\dd}{\mathop{}\!\mathrm{d}}

\newcommand{\fbi}{{\textbf{(}\mathbf F\textbf{)}_{\mathbf I}}}
\newcommand{\prl}{{\textbf{(}\mathbf P\textbf{)}_{\mathbf \lm}}}

\newtheorem{theorem}{Theorem}[section]
\newtheorem{proposition}[theorem]{Proposition}
\newtheorem{lemma}[theorem]{Lemma}
\newtheorem{corollary}[theorem]{Corollary}
\newtheorem{remark}[theorem]{Remark}
\newtheorem{definition}[theorem]{Definition}
\newcommand{\brm}{\begin{remark}\rm}
\newcommand{\erm}{\end{remark}}
\newcommand{\bdf}{\begin{definition}\rm}
\newcommand{\edf}{\end{definition}}
\newcommand{\bte}{\begin{theorem}}
\newcommand{\ete}{\end{theorem}}
\newcommand{\bpr}{\begin{proposition}}
\newcommand{\epr}{\end{proposition}}
\newcommand{\ble}{\begin{lemma}}
\newcommand{\ele}{\end{lemma}}
\newcommand{\bco}{\begin{corollary}}
\newcommand{\eco}{\end{corollary}}
\newcommand{\beq}{\begin{equation}}
\newcommand{\eeq}{\end{equation}}
\newcommand{\bdm}{\begin{displaymath}}
\newcommand{\edm}{\end{displaymath}}

\newcommand{\graf}[1]{\left\{\begin{array}{ll}#1\end{array}\right.}

\def\sideremark#1{\ifvmode\leavevmode\fi\vadjust{\vbox to0pt{\vss
 \hbox to 0pt{\hskip\hsize\hskip1em \vbox{\hsize2.1cm\tiny\raggedright\pretolerance10000 \noindent #1\hfill}\hss}\vbox to15pt{\vfil}\vss}}}

\begin{document}
\numberwithin{equation}{section}
\parindent=0pt
\hfuzz=2pt
\frenchspacing

\title[Sharp uniqueness and spikes condensation results]
{Sharp estimates, uniqueness and spikes condensation for superlinear free boundary problems arising in plasma physics}

\author[]{Daniele Bartolucci$^{(1,\dag)}$, Aleks Jevnikar$^{(2)}$, Ruijun Wu$^{(3)}$}


\thanks{$^{(1)}$Daniele Bartolucci, Department of Mathematics, University
of Rome {\it "Tor Vergata"},  Via della ricerca scientifica n.1, 00133 Roma,
Italy. e-mail:bartoluc@mat.uniroma2.it}

\thanks{$^{(\dag)}$ Research partially supported by the MIUR Excellence Department Project MatMod@TOV
awarded to the Department of Mathematics, University of Rome Tor Vergata and by PRIN project 2022, ERC PE1\_11,
"{\em Variational and Analytical aspects of Geometric PDEs}". D. Bartolucci is member of the INDAM Research Group
"Gruppo Nazionale per l'Analisi Matematica, la Probabilità e le loro Applicazioni".}

\thanks{$^{(2)}$Aleks Jevnikar, Department of Mathematics, Computer Science and Physics, University of Udine, Via delle Scienze 206, 33100 Udine, Italy.
e-mail:aleks.jevnikar@uniud.it}

\thanks{$^{(3)}$Ruijun Wu, School of mathematics and statistics, Beijing Institute of Technology, Zhongguancun South Street No. 5, 100081 Beijing, P.R.China.
e-mail:ruijun.wu@bit.edu.cn}

\begin{abstract}
We are concerned with Grad-Shafranov type equations, describing in dimension $N=2$ the equilibrium configurations of a plasma in a Tokamak.
We obtain a sharp superlinear generalization of the result of Temam (1977) about the linear case,
implying the first general uniqueness result ever for superlinear free boundary problems arising in plasma physics.
Previous general uniqueness results of Beresticky-Brezis (1980) were concerned with globally Lipschitz nonlinearities.
In dimension $N\geq 3$ the uniqueness result is new but not sharp, motivating the local analysis of
a spikes condensation-quantization phenomenon for {superlinear} and {subcritical} singularly perturbed Grad-Shafranov type free
boundary problems, implying among other things a converse of the results about spikes condensation in Flucher-Wei (1998) and Wei~(2001).
Interestingly enough, in terms of the ``physical'' global variables, we come up with a
concentration-quantization-compactness result sharing the typical features of critical problems (Yamabe $N\geq 3$, Liouville $N=2$)
but in a subcritical setting, the singular behavior being induced by a sort of
infinite mass limit, in the same spirit of Brezis-Merle (1991).\end{abstract}
\maketitle

%\tableofcontents

%\section{\bf ....}

\section{\bf Introduction}
\setcounter{section}{1}

\setcounter{equation}{0}
Let $\Omega\subset \R^N$,~$N\ge 2$, be a bounded domain of class~$C^{2,\beta}$ for some~$\beta\in (0,1)$, and set
$$
p_{_N}=\graf{+\ii,\;N=2\\ \frac{N}{N-2},\; N\geq 3.}
$$
For $p\in [1,p_{_N})$ and~$I>0$ we are concerned with the free boundary problem (\cite{BeBr,Te2})
$$
\graf{-\Delta {\rm v} =[{\rm v}]_+^p\quad \mbox{in}\;\;\om\\ \\
\quad\;\; {\rm v}=\gamma \qquad \mbox{on}\;\;\pa\om\\ \\
\bigints\limits_{\om} {\dsp [{\rm v}]_+^p}=I
}\qquad \fbi
$$
where $[t]_+$ denotes the positive part and the unknown is $(\gamma,{\rm v})\in \R \times
C^{2,\beta}(\ov{\om})$.
Let~$p>1$ and set $I=\lambda^q$, where~$q=\frac{p}{p-1}$ is the H\"older conjugate to~$p$, and consider the new variables~$(\alpha,\psi)\in \R\times C^{2,\beta}(\overline{\Omega})$ defined via
\begin{align}\label{eq:change variable}
   \begin{cases}
    \gamma= \lambda^{\frac{1}{p-1}}\alpha, \vspace{2mm}\\
    {\rm v} = \lambda^{\frac{1}{p-1}}(\alpha+\lambda\psi).
   \end{cases}
  \end{align}
  Then, for $p>1$, $\fbi$ is equivalent to the following system with parameter~$\lambda>0$,
$$
\graf{-\Delta \psi =[\al+{\lm}\psi]_+^p\quad \mbox{in}\;\;\om\\ \\
\quad\;\; \psi=0 \qquad\qquad \;\;\; \mbox{on}\;\;\pa\om\\ \\
\bigints\limits_{\om} {\dsp \left[\al+{\lm}\psi\right]_+^p}=1,
}\qquad \prl
$$
 where the unknown is $(\alpha,\psi)$ while, for $p=1$, $\prl$ is already equivalent to the corresponding problem classically discussed in literature (\cite{Te,Te2}). We will write~$(\alpha_\lambda,\psi_\lambda)$ to denote any solution of $\prl$ for fixed $\lm>0$.

\

  Due to its relevance to Tokamak's plasma physics \cite{Kad,Stac}, a lot of work has been done to understand
existence, uniqueness, multiplicity of solutions and existence/non existence/structure of the free boundary $\pa\{x\in\om\,:\,v>0\}$ of $\fbi$, see
\cite{AMN},\cite{BaMa}-\cite{BSp},\cite{BeBr,CF,CPY,dam},\cite{FL}-\cite{Gal},\cite{KNS}-\cite{Kor3},\cite{LP,Liu,Mar},\cite{Pudam}-\cite{Sij1},
\cite{suz}-\cite{wol}.\\
In particular due to the fundamental results of Beresticky-Brezis (\cite{BeBr}) it is well known that for any~$\lambda>0$
there exists at least one solution of~$\prl$.
  On the contrary, with the only exception of disks (\cite{BSp}) and balls (which we pursue in Appendix B below), uniqueness results are quite rare, as far as $p>1$.
  Partial information have been recently obtained about uniqueness and monotonicity of bifurcation diagrams in~\cite{BHJY,BJ1,BJ3}.
  Generalizing a sharp estimate obtained by Temam~\cite{Te2} for $p=1$, we obtain the very first, at least to our knowledge,
  general and neat uniqueness result for~$\fbi$ in dimension two for any~$p<+\infty$, extending that by Bandle--Sperb (\cite{BSp}) for the disk.
  Indeed, our approach also yields a general uniqueness result in higher dimension which unfortunately
  is not neat, motivating the analysis of a general ``infinite mass''-type singular limit.
Interestingly enough, in terms of the ``physical'' global variables, we come up with a
concentration-quantization-compactness result sharing the typical features of critical problems (Yamabe $N\geq 3$, Liouville $N=2$)
but in a subcritical setting, the singular behavior being induced by a sort of
infinite mass limit, in the same spirit of Brezis-Merle (\cite{BM}). Among other things, general a priori estimates for solutions of $\prl$ naturally follow 
in the same spirit of Li (\cite{yy})

\

From now on we assume that the domain has unit volume, $|\Omega|=1$.
  This can always be achieved without loss of generality in dimension two due to the scaling invariance of the problem.
  In higher dimension it would require a minor modification of the integral constraints, we disregard this issue for the sake of simplicity.
  For~$p\ge 1$ let us define
  \begin{align}\nonumber
   \Lambda(\Omega,p)
   =\inf\limits_{w\in H^1_0(\Omega), w\neq 0}
\frac{\int_\Omega |\nabla w|^2}{\left(\ino |w|^{p}\right)^{\frac{2}{p}}}\,,
  \end{align}
  which is related to the best constant in the Sobolev embeddings~$\|w\|_{p}\leq \mathcal{C}_S(\Omega, p)\|\nabla w\|_2$ via
  \begin{align}\nonumber
  \mathcal{C}_S(\Omega,p)=\Lambda(\om,p)^{-1/2}, \qquad \mbox{ for } p\in[1,2p_{N}).
  \end{align}

\

{\bf Definition} {\it We say that a solution $(\all,\pl)$ of {\rm $\prl$} is {\bf non-negative}
{\rm[}resp. {\bf positive}{\rm]} if
$\all\geq 0$ {\rm[}resp. $\all>0${\rm]}}.

\

It is well known since \cite{BeBr} that for any $\lm>0$ there exists at least one solution of
$\prl$. Remark that if $(\all,\pl)$ is a solution of $\prl$ then by the maximum principle
$\pl>0$ in $\om$ and in particular, if $\all\geq0$, then $|\om|=1$ readily implies $\all\in [0,1]$.

\

It has been recently proved in \cite{BJ1} that the following quantity is well-defined:
$$
\lm^*_+(\om,p)=
\sup\left\{\mu>0\,:\, \all>0 \mbox{ for any non-negative solution } (\all,\pl)
\mbox{ of } \prl \mbox{ with }\lm<\mu\right\}.
$$

\

{\bf Definition} {\it Let $(\all,\pl)$ be a solution of {\rm $\prl$}.
The {\bf energy} of $(\all,\pl)$ is
$$
\el=\frac12\ino |\nabla \pl|^2=\frac12 \ino  \mbox{\rm$\rl$}\pl,\quad \mbox{ where } \quad
\mbox{\rm$\rl$}\equiv [\all+\lm\pl]_+^p.
$$
}

The following result has been recently proved in~\cite{BJ1}.

\

{\bf Theorem A.}{(\cite{BJ1})} {\it Let $N\geq 2$ and $p\in [1,p_{_N})$, then $\lm^*_+(\om,p)\geq \frac1p {\Lambda(\om,2p)}$
where the equality holds if and only if $p=1$. Moreover, for any $\lm<\frac1p {\Lambda(\om,2p)}$ there exists a unique positive solution to {\rm $\prl$},
defining a
real analytic curve from $[0,\frac1p {\Lambda(\om,2p)})$ to $[0,1] \times C^{2,\beta}(\ov{\om})$, denoted by
$$
\mathcal{G}(\om)=\left\{(\all,\pl), \lm\in[0,\frac1p {\Lambda(\om,2p)})\right\},
$$
such that $\al_0=1$, $2 E_0=2 E_0(\om)$ is the torsional rigidity of $\om$ and
$$
\frac{d \all}{d\lm}<0,\quad \frac{d \el}{d\lm}>0,\quad \forall (\all,\pl)\in\mathcal{G}(\om).
$$}

It is worth to remark that $\gamma_{_I}$ cannot be a monotone function of $I$ in $[0,(\lm^*_+(\om,p))^{q}]$, 
while in fact $\all$ does, at least for $N=2$ in a disk and for 
$\lm \in [0,\frac1p {\Lambda(\om,2p)})$ from Theorem A, see \cite{BJ1} and Remark \ref{rem9.3} in Appendix B below.\\
With the unique exception of the case of the disk in \cite{BSp}, which we generalize here to balls for $N\geq 3$ (see Appendix B),
as far as $p>1$ we are not aware of any unconditional uniqueness result for solutions of $\prl$. Theorem A suggests a possible
solution to this problem, which is to understand whether or not the unique \un{positive} solutions of $\prl$
are in fact the unique solutions \un{at all} of $\prl$. In other words we ask whether or not any solution of
$\prl$ is non-negative as far as $\lm<\frac1p {\Lambda(\om,2p)})$. For $N=2$ we provide a first answer to this question
with a sharp estimate which improves Theorem A. Let $\mathbb{D}_N$ be the ball of unit area in $\R^N$ and for
$\om\subset \R^N$ define the \un{positivity threshold},
$$
\lm^*(\om,p)\coloneqq
\sup\left\{\mu>0\,:\, \all>0 \mbox{ for any solution } (\all,\pl) \mbox{ of } \prl \mbox{ with }\lm<\mu\right\}.
$$

For the ease of the presentation let us denote
\begin{align}\label{eq:lambda0}
 \lm_0(\om,p)\equiv \left(\frac{8\pi}{p+1}\right)^{\frac{p-1}{2p}}
\Lambda^{\frac{p+1}{2p}}(\om,p+1).
\end{align}

\bte\label{thm1.1} Let $N=2$ and $p\in [1,+\ii)$, then
$$
\lm^*(\om,p)\geq \lm_0(\om,p),
$$
and the equality holds if and only if either $p=1$ or $p>1$, $\om=\mathbb{D}_2$ and $\lm\pl$ is a minimizer of~$\Lambda(\mathbb{D}_2,p+1)$.
\ete
This is a sharp nonlinear generalization of the well known result by Temam (\cite{Te2}), providing at once a sharp estimate for
the value of the positivity threshold, whose existence was first proved for variational solutions in \cite{BaMa}, see also \cite{BSp}.
As an immediate consequence of Theorem \ref{thm1.1}, we deduce from Theorem A the following,
\bco\label{cor:1.12}
Let $N=2$,~ $p\in [1,+\ii)$ and~$\lm\leq \lm_0(\om,p)$. Then, any solution~$(\all, \pl)$ of {\rm $\prl$} satisfies $\all\geq 0$, where
$\all=0$ if and only if either $p=1$ or $p>1$,
$\om=\mathbb{D}_2$ and $\lm=\lm_0(\mathbb{D}_2,p)$.

As a consequence for any
$$
 0\leq \lm<  \min\left\{\frac{1}{p}\Lambda(\om,2p),\lm_0(\om,p)\,\right\}
$$
there exists one and only one solution to {\rm $\prl$}, which is a positive solution.
In particular the set of solutions in this range coincides with $\mathcal{G}(\om)$.
\eco
At least to our knowledge, with the unique exception of the above mentioned cases of disks and balls,
this is the \emph{first general uniqueness result ever} for solutions of $\prl$ with $p>1$. Previous pioneering
general uniqueness results were concerned with globally Lipschitz nonlinearities (\cite{BeBr}). We tend to believe that
$\frac{1}{p}\Lambda(\om,2p)<\lm_0(\om,p)$ for $p>1$, which we can just prove for $p$ large, see  Proposition \ref{lme:31} in
Section \ref{sec3}.
Note that, as a consequence of Theorem A and Theorem \ref{thm1.1}, for a disk this is actually true:
\begin{align}\nonumber
 \frac{1}{p}\Lambda(\mathbb{D}_2,2p)\leq \lm_0(\mathbb{D}_2,p), \quad \mbox{for all } 1\leq p <+\infty,
\end{align}
where the equality holds iff~$p=1$.

\

To prove Theorem \ref{thm1.1} we need a sharp energy estimate of independent interest which
improves a result in \cite{BJ1}. Here and in the rest of this paper, as far as $\all<0$,
we denote by
$$\om_+=\{x\in\om\,:\,\all+\lm\pl>0\},$$
known as the ``plasma region''.
Remark that, by a standard trick based on the Sard
Lemma, we have that if $\all<0$ then

\beq\label{enrg1.0}
2\el =\int\limits_{\om_+} |\nabla \pl|^2-\frac{\all}{\lm},
\eeq
see Lemma \ref{enrg0} below for a proof of this fact.

\ble\label{altype2.intro} Let $N=2$ and $(\all, \pl)$ be a solution of {\rm $\prl$} with $\all<0$.
Then
$$
\int\limits_{\om_+} |\nabla \pl|^2\leq \frac{p+1}{8\pi},
$$
where the equality holds if and only if $\om_+$ is a disk and $\pl$ is radial in $\om_+$.
\ele

For $N\geq 3$ we don't have a sharp estimate as in Lemma \ref{altype2.intro}.
However we can prove the following
\ble\label{altype3} Let $N\geq 3$, $p\in [1,p_{N})$ and $(\all, \pl)$ be a solution of {\rm $\prl$}.\\
$(a)$ If $\all\geq 0$ then there exists $\mu_+ \in (\frac{N}{2(N-1)},+\ii)$ such that,
\beq\label{n3.1}
\el\leq \frac{ (p+1)}{4N^2 (\omega_{\sscp N-1})^{\frac 2 N}}\mu_+  +
\frac{\all}{2\lm}\left(\all^{p}-\frac{1}{\mu_+}\right)\mu_+.
\eeq
$(b)$ If $\all<0$ then there exists
$\mu_- \in (\frac{N}{2(N-1)},+\ii)$ such that,
\beq\label{n3.2}
\el\leq
\frac{(p+1)}{4N^2 (\omega_{\sscp N-1})^{\frac 2 N}}\frac{\mu_-}{|\om_+|^{1-\frac{2}{N}}}+ \frac{|\all|}{2\lm}.
\eeq
\ele

\brm\label{rmmu} {\it By the claim, in principle $\mu_{+}$ and $\mu_{-}$ could depend
on the solution $(\all,\pl)$. It is tempting to try to understand whether or not there exist universal values of this sort.
This is part of the motivation behind the rest of this paper. However, as we will see in Corollary \ref{altype3.not} below, there
is no chance in dimension $N\geq 3$ to come up  with a uniform estimate as in the claim of Lemma \ref{altype2.intro}.}
\erm

Unfortunately we miss enough informations about $\mu_{-}$ suitable to derive a uniqueness result. Therefore
we argue in another way to obtain the first, at least to our knowledge, general uniqueness result
for solutions of $\prl$  with $N\geq 3$ and $p>1$.

For $N\geq 3$, let $R_N>0$ be a dimensional constant such that~$|B_{R_N}(0)|=1=|\om|$.
For~$p<p_{_N}$, pick any $s\in (p,p_{_N})$, then set $k_s\coloneqq 1-\frac{s}{p_{_N}}\in(0,1)$ and
$$
\ell_s(\om,p)\coloneqq
\left(\left(\frac{N (N-2) k_s}{R_N^{N k_s}}\right)^{\frac{p-1}{s-p}}
{\Lambda^{{p+1}}(\om,p+1)}\right)^{\frac{s-p}{p(2s-(p+1))}}.
$$

Then we have
\bte\label{thm1.5} Let $N\geq 3$, $p\in (1,p_{_N})$ and $(\all, \pl)$ be a solution of {\rm $\prl$} with
$\all\leq 0$.  Then

$$
\lm^*(\om,p)> \ell_s(\om,p).
$$
In particular for any $$
 0\leq \lm<  \min\left\{\frac{1}{p}\Lambda(\om,2p),
 \ell_s(\om,p)\,\right\}
$$
there exists one and only one solution to {\rm $\prl$}, which is a positive solution.
In particular the set of solutions in this range coincides with $\mathcal{G}(\om)$.
\ete

\bigskip

For $N\geq 3$ and $\om$ a ball, we have the value of $\lm^*(\om,p)$ in terms of the unique solution of the Emden equation, see Appendix B. Otherwise
in general it seems not so easy to improve Lemma \ref{altype3} and consequently Theorem \ref{thm1.5}.
Therefore, in dimension $N\geq 3$, we lack the sharp estimates needed to come up as well
with a neat uniqueness result.
We are thus motivated to further investigate what kind of asymptotic
behavior is really allowed for solutions of $\prl$ for large $\lm$ and $|\all|$,
with the hope that this could help in a better understanding of \rife{n3.2}.
Interestingly enough, following this route we naturally come up with the description of a new concentration-quantization-compactness
phenomenon (see Theorem \ref{thm7.2.intro} below) in a \underline{subcritical} setting (since we have $p<p_{_N}$),
the singular behavior being forced by a sort of "infinite mass" limit, in the same spirit of Theorem 4 in \cite{BM} for Liouville-type equations.\\

\medskip

It is worth to shortly illustrate this point via a comparison with mean field type problems for $N=2$ (\cite{clmp2}), where the variable $\psi$
in fact satisfies an equation of the form $-\Delta \psi=f(\al+\lm\psi)$, where $f(t)=e^t$ and $\int_{\om}f(\al+\lm\psi)=1$. Thus, putting
$u=\al+\lm \psi$ the equation takes the form of the classical Liouville equation $-\Delta u= \lm e^{u}$, $\int_{\om}e^u=1$. A subtle point arise
here since in this way the ``Lagrange multiplier'' $\al$ (see also the Appendix below) is hidden in the function $u$ and it is not anymore clear how to
single out the functional dependence of $\al$ from $\lm$ in the ``infinite mass'' limit $\lm\to+\ii$. Therefore, as far as we are concerned with the limit $\lm\to +\ii$,
it seems better to write $\tilde v=\lm \psi$ to come up with the equation $-\Delta \tilde v= \tilde \mu e^{v}$, where $\tilde \mu=\lm e^{\al}$. This has also the advantage of preserving the sign of $\tilde v$, that is if $\psi=0$ on $\pa \om$ then $\psi\geq 0$ and $\tilde v\geq 0$ as well, which allows
in turn to handle the limit $\lm\to +\ii$ by a careful use of the ``infinite mass'' blowup analysis of~\cite{BM}.
Actually nontrivial examples of singular limits of this sort have been recently pursued in \cite{LSMR} and \cite{PKI}.

\medskip

Sticking to this analogy, the equation in $\prl$ takes the same
form but with $f(t)=[t]_+^p$. Remark that, for solutions of $\prl$, $\al\leq 1$ and that if $|\al|$ and $\lm$ are bounded then it is well known
that~$\psi$ is uniformly bounded (\cite{BeBr}, see also \cite{BJ1}). Actually if $\al\to -\ii$ then necessarily $\lm\to +\ii$, see Theorem \ref{thm:7.1} below.
Also (see Theorem 1.2 in \cite{BHJY}) if $\lm\to +\ii$ then necessarily $\all<0$. Therefore, to figure out the functional relation of $\lm$ and $\al$ which yields a nontrivial, if any, ``infinite mass'' limit as $\al \to -\ii$ and
$\lm\to +\ii$, we argue as above and encode the singular behavior in $v=\frac{\lambda}{|\alpha|}\psi$, which has the advantage to solve
the model equation $-\Delta v=\mu[v-1]_+^p$ with parameter~$\mu=\lm|\al|^{p-1}\to+\ii$ and at the same time to preserve
the positivity of $v$ as far as $\psi=0$ on $\pa \om$. However, along a sequence~$\mu_n=\lambda_n |\alpha_n|^{p-1
}\to+\infty$, the solutions $\vn$ does not blow up but undergo to a condensation phenomenon
of spike type which was already pursued in \cite{FW,We}.
This is our motivation, essentially in the same spirit of \cite{BM}, to drop at first both the integral constraint and the boundary conditions and attack the local
asymptotic analysis of the ``infinite mass'' or ``singularly perturbed'' subcritical problem,
\beq\label{vn.1.intro}
\graf{-\Delta \vn =\mu_n [\vn-1]_+^p\quad \mbox{in}\;\;\om,\\ \\
\mu_n \to +\ii,\\ \\
\vn\geq 0,}
\eeq
where $\om \subset \R^N$, $N\geq 3$ is an open bounded domain and $p\in(1,p_{N})$. Remark that this general local analysis is also motivated by
applications arising from other physical problems usually equipped with different boundary conditions, see \cite{BeF}, \cite{Dolb},
\cite{Leub}, \cite{Liv} and references quoted therein.\\

We point out that refined estimates are at hand for
\rife{vn.1.intro} in the ``critical'' case $p=p_{_N}$ with $\mu_n$ uniformly bounded, see \cite{Br}.
We also refer to \cite{wy} where a full generalization of the concentration-compactness
theory in \cite{BM} has been derived directly for solutions of the equation in $\fbi$ with $p=p_{_N}$.
However for $p=p_{_N}$ there is no uniform estimates and solutions
blow up as $\lm$ converges to some specific quantized values,
see~\cite[Theorem 4]{wy}. We will also comment out, in Remark \ref{rem:suz} below, concerning the concentration-compactness theory obtained in \cite{suz} for solutions
of the equation in $\fbi$.

\bigskip

The singular limit of~\eqref{vn.1.intro} will be described in the language of spikes condensation, the basic notions are given below.

\bdf\label{spikedef} {\it Define $\eps_n^2:=\mu_n^{-1}$. The {\bf spikes set} relative to a sequence
of solutions of \eqref{vn.1.intro} is any set $\Sigma \subset \ov{\om}$ such that for any $z\in \Sigma$ there exists
$x_n\to z$, $\sg_z>0$ and $M_z>0$ such that:
\begin{itemize}
 \item $\vn(x_n)\geq 1+\sg_z$, $\forall n\in \N$;
 \item $\eps_n =o(\mbox{dist}(x_n,\pa \om))$;
 \item  $\mu_n^{\frac N 2}\int\limits_{B_{r}(z)\cap\om} [\vn-1]_+^p\leq M_z$, for some $r>0$.
\end{itemize}
Any point $z\in \Sigma$ is said to be a spike point.}
\edf

To encode the information about the spikes set we will use a formal finite combination of points in~$\ov{\Omega}$, denoted by
\begin{align}\nonumber
 \mathcal{Z}=k_1 z_1 +k_2 z_2 + \cdots + k_m z_m=\sum_{j=1}^m  k_j z_j,
 \end{align}
for some~$m\in \mathbb{N}$,~$k_j\in \mathbb{N}$ and~$z_j\in \ov{\Omega}$. Here $\mathbb{N}$ is the set of positive integers.
The degree of~$\mathcal{Z}$ is defined to be
\begin{align}\nonumber
 |\mathcal{Z}|\equiv \deg(\mathcal{Z})=k_1+\cdots+k_m,
\end{align}
while the underlying spikes set is denoted by,
\begin{align}\nonumber
 \Sigma(\mathcal{Z})= \braces{z_1, \cdots, z_m},
\end{align}
whose cardinality is $\#\mathcal{Z}\equiv\#\Sigma(\mathcal{Z})=m$. We will some time write~$\mathcal{Z}_m$ to emphasize its cardinality.
Note that we always have $|\mathcal{Z}|\geq \#\Sigma(Z)$ with equality holds iff all $k_j$ are 1.

\

In the incoming definition, the quantities~$R_0$, $w_0$ and $M_{p,0}$ are defined respectively
in \eqref{entire.Ra}, \eqref{entire.sol} and \eqref{entire.mp0} below.
In particular,~$w_0$ is the unique solution (also called a ground state in \cite{FW}) of \eqref{entire.a}, see Lemma~\ref{lem5.3}.

\bdf\label{spikeseqdef} {\it A sequence of solutions of \eqref{vn.1.intro} with a nonempty spikes set $\Sigma=\braces{z_1,\cdots, z_m} \subset \ov{\om}$ is said to be a {\bf $\mathcal{Z}$-spikes sequence} (or a spikes sequence for short) relative to~$\Sigma$ if
\begin{itemize}
 \item $\Sigma(\mathcal{Z})=\Sigma$;
 \item  there exist sequences~$\{z_{i,n}\}$, $i\in \{1,\cdots,|\mathcal{Z}|\}$, converging to~$\Sigma$;
 \item there exist $0<R_n \to +\infty$, a natural number~$n_*\in \N$ and positive constants~$C_*>0$, $t\geq 1$ and~$\sigma>0$, such that the following hold:
\end{itemize}

 \begin{enumerate}
  \item[$(i)$] $\eps_n R_n \to 0$ and, for each~$i\in\{1,\cdots,|\mathcal{Z}|\}$, $z_{i,n}\to z$
  for some~$z\in \Sigma$; and conversely, for each~$z\in\Sigma$ there is at least one sequence~$\{z_{i,n}\}$ which converges to~$z$;

  \item[$(ii)$] $B_{4\eps_nR_n}(z_{i,n})\Subset \om$ for any $i\in\{1,\cdots,|\mathcal{Z}|\}$ or equivalently
  $$
  B_{4R_n}(0)\Subset \om_{i,n}:=\frac{\Omega- z_{i,n}}{\eps_n},\quad \forall i\in\{1,\cdots,|\mathcal{Z}|\};
  $$

  \item[$(iii)$] if~$|\mathcal{Z}|\ge 2$ then
       \begin{align}\nonumber
        B_{2 \eps_n R_{n}}(z_{\!\ell,n})\cap B_{2 \eps_n R_{n}}(z_{\!i,n})=\emptyset,\;\forall\;\ell\neq i,
       \end{align}
        and
        $$
        \vn(z_{1,n})=\max\limits_{\ov{\om}}\vn\geq 1+\sigma,
        $$
       $$
\vn(z_{i,n})=
\max\limits_{\ov{\om} \setminus
\{\underset{\ell = 1,\cdots,|\mathcal{Z}|, \ell\neq i}{\cup}B_{2 \eps_n R_{n}}(z_{\!\ell,n})\}}
\vn\geq 1+\sigma,\quad i\geq 2;
$$

  \item[$(iv)$] for each~$i\in\{1,\cdots,|\mathcal{Z}|\}$, the rescaled functions
        \begin{align}\nonumber
         w_{i,n}(y)\coloneqq v_n(z_{i,n}+ \eps_n y), \quad y\in \om_{i,n},
        \end{align}
        satisfy
        \begin{align}\nonumber
          \|w_{i,n}-w_0\|_{C^2(B_{2R_n}(0))}\to 0;
        \end{align}

   \item[$(v)$] setting ${\Sigma}_{r}\coloneqq
   \underset{i =1,\cdots,m}{\cup}B_{r} (z_{i})$, for some $r>0$, and
   $\widetilde{\Sigma}_{n}=\underset{i =1,\cdots,|\mathcal{Z}|}{\cup}B_{2\eps_nR_n} (z_{i,n})$ then
$$
0\leq \max\limits_{\ov{\om} \setminus {\Sigma}_{r} } \vn\leq
C_r \eps_n^{\frac{N}{t}},\quad
0\leq \max\limits_{\ov{\om} \setminus \widetilde{\Sigma}_{n} } \vn\leq
 \frac{C_*}{R_n^{N-2}}, \quad \forall\, n>n_*;
$$
\item[$(vi)$]
the ``plasma region" $\om_{n,+}:=\{x\in \om\,:\,\vn>1\}$ consists of ``asymptotically round points''
in the sense of Caffarelli-Friedman {\rm(}\cite{CF}{\rm)} i.e., for any $r<R_0<R$,
$$
\underset{i=1,\cdots,|\mathcal{Z}|}{\cup}B_{\eps_n r}(z_{i,n})\Subset\{x\in \om\,:\,\vn>1\}\Subset
\underset{i=1,\cdots,|\mathcal{Z}|}{\cup}B_{\eps_n R}(z_{i,n}),\quad \forall\,n>n_*;
$$

  \item[$(vii)$] each sequence~$\{z_{i,n}\}$ carries the fixed mass $M_{p,0}$,
$$
\lim\limits_{n\to +\ii} \int \limits_{B_{\eps_n R_{n}}(z_{i,n})}\mu^{\frac N2}_n[\vn-1]_+^p=
\lim\limits_{n\to +\ii} \int \limits_{B_{R_{n}}(0)}[w_{i,n}-1]_+^p=
M_{p,0},\;
$$
and
$$
\lim\limits_{n\to +\ii} \int \limits_{\om}\mu^{\frac N2}_n[\vn-1]_+^p=|\mathcal{Z}| M_{p,0};
$$

 \item[$(viii)$] for any $\phi \in C^{0}(\ov{\om}\,)$
 $$
 \lim\limits_{n\to +\ii}\mu^{\frac N2}_n\int\limits_\om [\vn-1]_+^p\phi=
 M_{p,0}\sum\limits_{j=1}^m k_j\phi(z_j),
 $$
 where $k_j\geq 1$ are the local multiplicities (see also Remark \ref{multi:spike} below),
  $\sum\limits_{j=1}^m k_j=|\mathcal{Z}|\geq m$ and $k_j=1$, for any $j=1,\cdots,m$ if and only if $|\mathcal{Z}|=\#\Sigma(\mathcal{Z})$.
\end{enumerate}}
\edf

\brm\label{multi:spike}{\it Each local profile approximating $w_0$ in $B_{\eps_n R_{n}}(z_{i,n})$ is said to be a
{\bf spike}. By definition in general we only have~$m\leq |\mathcal{Z}|$.
 The equality holds iff all spike points are {\bf simple} in the following sense.
 Since~$\Sigma$ is a finite set we may pick~$r>0$ such that, for~$z\in \Sigma$,
 $B_r(z)\cap \Sigma=\braces{z}$.
 Such a spike point~$z\in \Sigma$ is said to be of multiplicity $m_z\in \N$ relative to the sequence~$\vn$
 if there are, among those in the Definition \ref{spikeseqdef},
 ~$m_z$ sequences~$\{z_{i_j,n}\}_{n\geq  1}$, $j\in \{1,\cdots,m_z\}$, that converge to~$z$.
 In this case,
 \begin{align}\nonumber
  \lim_{n\to+\infty} \mu_n^{\frac{N}{2}}\int_{B_r(z)} [v_n-1]_+^p \dd x = m_z M_{p,0}.
 \end{align}
 A spike point is {\bf simple} if it is of unit multiplicity, $m_z=1$.}
\erm
\brm\label{remspikedef}{\it We point out that, in view of $(iv)$ in the definition,
any spikes sequence also satisfies the following properties.
First of all, in view of \cite[Theorem 4.2]{DGP},
it is readily seen that, for any $n$ large enough,
$$
\max\limits_{B_{R_n}(0)}w_{i,n}=w_{i,n}(0)=\vn(z_{i,n})=\max\limits_{B_{R_n\eps_n}(z_{i,n})}\vn,
$$
where $z_{i,n}$ is the unique and non-degenerate maximum point of $\vn$ in $B_{2 R_n\eps_n}(z_{i,n})$,
see also the proof of Theorem \ref{thm5.7.intro} for further details.\\
Also, by the explicit expression of $w_0$ in \eqref{entire.sol},
 for each~$i\in\{1,\cdots,|\mathcal{Z}|\}$, on the neck domain~
 $\braces{x\in\om \mid 2\eps_n R_0 \leq |x- z_{i,n}| \leq 2\eps_n R_n}$, for $n$ large we have
        \begin{align}\label{neck}
         \parenthesis{\frac{1}{2}\frac{R_0 \eps_n}{|x-z_{i,n}|}}^{N-2}
         \leq v_n(x)
         \leq \parenthesis{2\frac{R_0 \eps_n}{|x-z_{i,n}|}}^{N-2}.
        \end{align}
}
\erm

\bigskip

A detailed asymptotic analysis is worked out in \cite{FW} of mountain-pass type solutions of \rife{vn.1.intro} with $\vn=0$ on $\pa\om$, including the existence, uniqueness and location of the unique spike point of the corresponding spike sequence.\\
Moreover, let~$\mathcal{Z}$ satisfy $|\mathcal{Z}|=\#\Sigma(\mathcal{Z})$, i.e. all $k_j$ equal 1.
It has been shown in \cite{We} (see also Theorem \ref{exist:plasma} below) that $\mathcal{Z}$-spikes sequences~$\{\vn\}$ with $\vn|_{\pa\om}=0$ exists, where $\Sigma(\mathcal{Z})=\{z_1,\cdots,z_m\}$ and all the spike points are simple.
In particular $(z_1,\cdots, z_n)$ is a critical point of
\begin{align}\label{eq:simpleKR}
 \mathcal{H}(x_1,\cdots, x_m; \mathcal{Z})
 =\sum\limits_{j=1}^m  H(x_j,x_j)+\sum\limits_{i\neq j} G(x_i,x_j),
\end{align}
where~$G(x,y)=G_\om(x,y)$ is the Green function with Dirichlet boundary condition for~$\Omega$, i.e.
\begin{align}\nonumber
 \begin{cases}
  -\Delta_x G(x,y)= \delta_{y}(x) & x\in  \Omega, \\
  G(x,y)=0 & x\in  \pa\Omega,
 \end{cases}
\end{align}
and~$H$ is the regular part of the Green function:
$$H(x,y)=G (x,y)-\frac{1}{N(N-2)\omega_N}\frac{1}{|x-y|^{N-2}}.$$
Remark that~\eqref{eq:simpleKR} is invariant under the permutation of variables, so the critical set of~$\mathcal{H}$ above is in general not unique if~$m\ge 2$.
For~$m=1$, the function~$H(x,x)$ is the Robin function, whose minimum points are by definition the harmonic centers of the domain~$\Omega$, see for example
\cite{CT}.

\

Here we make another step in the direction pursued in \cite{FW,We} and, in the same spirit of \cite{BM},
consider the asymptotic behavior of a sequence of solutions of \eqref{vn.1.intro}
with no prior boundary conditions, obtaining a result of independent interest.
A delicate point arises at this stage which is about the structure, according to Definition \ref{spikeseqdef}, of sequences of solutions whose
spikes set is not empty, see Theorem \ref{thm5.7} and Remark \ref{remnobdy} below.
Inspired by \cite{BM}, we have the following spikes condensation-vanishing alternative for solutions of \eqref{vn.1.intro}.

% IN VIEW OF THE NEWS IN COROLLARY \ref{thm7.2.intro} MAY BE WE SHOULD REFORMULATE THIS ALREADY IN
% TERMS OF $\psi_n$, HOWEVER WE SHOULD WRITE DOWN A local FORM of $\prl$ IN A MORE NATURAL FORM,
% SAY SOMETHING LIKE \eqref{vn.1.intro}.

\bte\label{thm5.7.intro}
Let $\vn$ be a sequence of solutions of \eqref{vn.1.intro} with spikes set~$\Sigma\subset\ov{\om}$ such that,
\beq\label{mass.intro.1}
\mu_n^{\frac N 2}\ino [\vn-1]_+^p\leq C_0,
\eeq
\beq\label{mass.intro.2}
\mu_n^{\frac{N}{2}}\ino \vn^{t}\leq C_t,
\eeq
for some $t\geq 1$ and $C_0>0,C_t>0$.
Then,
\begin{itemize}
 \item[{\em either}] {\rm (a)} \emph{[Vanishing]} $[v_n-1]_+\to 0$ locally uniformly in $\om$, in which case
 for any open and relatively compact subset $\om_0\Subset \om$ there exists ${n}_0\in \N$ and $C>0$,
 depending on $\om_0$, such that $[v_n-1]_+=0$ in ${\om}_0$ and in particular
\beq\label{vanishnorm}
\|\vn\|_{L^{\ii}(\om_0)}\leq C\eps_n^{\frac N t},\quad {\forall\,n>n_0},
\eeq

 \item[{\em or\quad}] {\rm (b)} \emph{[Spikes-Condensation]} up to a subsequence (still denoted $\vn$) the interior spikes set relative to $\vn$ is nonempty and finite,  that is $\Sigma_0:=\Sigma \cap \om=\{z_1,\cdots,z_m\}$ for some $m \geq 1$ and for any open and relatively compact set $\om'$ such that
$\Sigma_0\subset \om'\Subset \om$, $\{\vn|_{\Omega'}\}$ is a $\mathcal{Z}$-spikes sequence relative to $\Sigma_0$, for some $\mathcal{Z}$ such that $\Sigma(\mathcal{Z})=\Sigma_0$.

\end{itemize}
\ete

Actually there is a version of Theorem \ref{thm5.7.intro} for changing sign solutions satisfying a uniform bound on~$\|[\vn]_{-}\|_{L^1(\om)}$,
see Remark \ref{rem:local}. As a consequence of Theorem \ref{thm5.7.intro}, in the same spirit of \cite{BM,ls}, we will see that
the corresponding sequences of solutions of $\prl$ undergo to a blowup-concentration-quantization phenomenon, see Theorems~\ref{thm7.2.intro} and~\ref{thm.1.21.intro} below.

\

Some comments are in order.
The assumption~\rife{mass.intro.1} is a necessary minimal
requirement for the possible singularities to be of spike type according to Definition \ref{spikedef}.
Of course, in principle other singularities may arise if we miss this assumption, as for example
happens to be the case for classical singularly perturbed elliptic problems (\cite{AMN,DY}) or either
in the ``infinite mass'' limit of the Liouville-Dirichlet problem (\cite{PKI,LSMR}).
On the other side, as far as we miss \rife{mass.intro.2},
we don't have such a nice description of the spikes set and $w_0$ in general is not the
unique admissible spike-type profile of the limiting problem, see Theorem \ref{thm5.7}, Remark \ref{remnobdy}
and Lemma \ref{minmass} for further details. In particular we refer to Lemma \ref{minmass} as the Minimal Mass Lemma,
which plays essentially the same role as its analogue does for Liouville-type equations, see \cite{bt,ls}.

\

The proof of Theorem \ref{thm5.7.intro} is not trivial since the spikes structure arise
in a subtle competition between the divergence of $\mu_n$ and the
the vanishing of $\ino [\vn-1]^p_+$ and $\ino |\vn|^t$.
Actually it is already nontrivial to prove the boundedness of the solutions sequence, since the
uniform bound for $\mu_n^{\frac N2}\ino [\vn-1]^p_+$ just yields a uniform local
estimate in $W^{1,N}$ for $\vn$.
We succeed in the description of
the singular limit via a detailed analysis of
the ``infinite mass'' limit of a suitably defined related problem. Concerning this point
we gather first some ideas from \cite{BM}, see Theorem \ref{thm5.2}, where the assumption $p<p_{_N}$ plays a crucial role.
Then the proof relies on the structure
of the solutions of the corresponding limiting problem in $\R^N$ (described in Lemma \ref{lem5.3} below) and
blow-up type arguments (see for example the Minimal Mass Lemma, i.e. Lemma \ref{minmass} below), although there is no blow-up of the $\|\vn\|_\ii$-norm here. A crucial intermediate result of independent interest,
which we call the Vanishing Lemma (see Lemma \ref{vanishle}), says that if $[\vn-1]_+\to 0$ uniformly in some $\om_0\Subset \om$,
then in fact $[\vn-1]_+\equiv 0$ in $\om_0$ for~$n$ large enough, whence in particular \eqref{vanishnorm} holds.
An almost equivalent and useful statement, which we call the Non-Vanishing Lemma, is stated in Remark \ref{vanish-equiv}.

\brm\label{rem:suz}{\it
In \cite{suz} a very interesting generalization of the Brezis-Merle {\rm(}\cite{BM}{\rm )}
concentration-compactness theory has been obtained for classical solutions of $-\Delta {\rm v}=[{\rm v}]_+^p$ in $\om$, $p\in (1,\frac{N-2}{N+2})$ satisfying
the ``scaling invariant'' bound $\int\limits_{\om}[{\rm v}]_+^{\frac{N}{2}(p-1)}\leq C$. The fact that the integral condition is scaling invariant allows
the authors to recover a full generalization of the results in \cite{BM} for classical solutions. However, the plasma problem requires a sort of control of
the infinite mass limit of $\int\limits_{\om}[{\rm v}]_+^p$ and, as far as we could check, this particular singular limit is not directly detectable
by the results in \cite{suz}. We do not exclude of course that further elaborations about the argument in \cite{suz} could be used in our context as well.
However, compared to \cite{suz}, we adopt an entirely different argument which seems to be more natural for the singular limit described in
Theorem \ref{thm5.7.intro}. In particular the assumption $p<p_{_N}$ plays a crucial role in our case, see the proof of Theorem \ref{thm5.2}.}
\erm

\

More care is needed in the analysis of the boundary behavior of solutions of
\beq\label{w5.1.intro}
\graf{-\Delta \vn =\mu_n[\vn-1]_+^p\quad \mbox{in}\;\;\om, \\ \\
\qquad \vn= 0 \quad \mbox{on}\;\;\pa\om,   \\
\bigints\limits_{\om} {\dsp |\an|^p[\vn-1]_+^p}=1,  \\
\mu_n=\lm_n|\an|^{p-1}\to +\ii, \lm_n \to+\ii,}
\eeq
which naturally arises in the analysis of sequences of solutions of {\rm $\prl$}.
Among other things, concerning \rife{w5.1.intro} we come up with an
almost converse of existence results in Wei (\cite{We}) and also in Flucher-Wei (\cite{FW}), where the discussion was limited to the description of mountain pass solutions.
Here indeed we gather some ideas from \cite{FW} about the boundary behavior of spikes which we use
to prove a crucial boundary version of the Non-Vanishing Lemma, see Lemma \ref{vanishle-bdy}. This is a major
step in the description of the singular global (i.e. up to the boundary) limit, in the sense that, even if we do not exclude
that spikes could sit on the boundary, still we can guarantee that $(ii)$ in Definition \ref{spikeseqdef} holds. As a consequence,
spikes (if any) accumulating at the boundary must carry the same profile and the same mass of the interior spikes. Remark that boundary spikes for example
exists in the classical spike-layer construction of singularly perturbed elliptic problems~\cite{NW}.
Also, boundary spikes have been constructed in \cite{LP} for the singularly perturbed problem \eqref{w5.1.intro}, disregarding of course the constraint, by including a multiplicative weight function.

\

Let us set $\vec{k}=(k_,\cdots, k_m)\in \R^m$, $k_j\neq 0$, then, borrowing a terminology typical of the two dimensional case, we refer to the function~$\mathcal{H}(x_1,\cdots, x_m; \vec{k})$ defined by
\begin{align}\label{K-H}
 \mathcal{H}(x_1,\cdots, x_m;\vec{k}) =
 \sum\limits_{j=1}^m k_j^2 H(x_j,x_j)+\sum\limits_{i\neq j=1}^m k_i k_j G(x_i,x_j),
\end{align}
as to the $\vec{k}$-Kirchhoff-Routh Hamiltonian.
Writing~$\vec{1}_m\equiv (1,\cdots, 1)\in \mathbb{Z}^m$, then~\eqref{eq:simpleKR} is just $\mathcal{H}(\bullet;\vec{1}_m)$, where
if $m=1$ we just have $\vec{1}_1=1$.

\

\begin{theorem}\label{thm:spike-vn.intro}
 Let $v_n$ be a sequence of solutions of~\eqref{w5.1.intro} satisfying \eqref{mass.intro.1} and assume that, either
 $$
 \om \;\; \mbox{is convex},
 $$
 or
 \begin{align}\label{eq:assump-2.intro}
  \mu_n^{\frac{N}{2}}\int_{\Omega} [\vn-1]_+^{p+1} \leq C_1,
 \end{align}
 for some $C_1<+\infty$.

 Then, possibly along a subsequence,
  $\vn$ is a $\mathcal{Z}$-spikes sequence relative to $\Sigma$, for some $\mathcal{Z}=\sum_{j=1}^m k_j z_j$,~$1\leq m<+\infty$. In particular, setting
${\Sigma}_{r}\coloneqq
   \underset{j =1,\cdots,m}{\cup}B_{r} (z_{j})$ for some $r>0$, then
\beq\label{est:nm2}
\vn(x)=
\eps_n^{N-2}\sum\limits_{i=1}^{|\mathcal{Z}|}M_{p,0}G(x,z_{i,n})+o(\eps_n^{N-2}),\quad \forall\,x\in \om\setminus {\Sigma}_{r},
\eeq
for any $n$ large enough and we have the ``mass quantization identity'',
$$
\lim\limits_{n\to +\ii}\left(\frac{\lm_n}{|\an|^{1-\frac{p}{p_{_N}}}}\right)^{\frac N 2}=
\lim\limits_{n\to +\ii} \mu_n^{\frac N 2}\int\limits_{\om}[\vn-1]^p_+= |\mathcal{Z}| M_{p,0}.
$$
Moreover, if $\Sigma\subset \om$,
then, $(z_1,\cdots,z_m)$ is a critical point of the
$\vec{k}$-Kirchoff-Routh Hamiltonian, where $\vec{k}\equiv (k_1,\cdots,k_m)$.
\end{theorem}

The proof of the fact that $(z_1,\cdots,z_m)$ is a critical point of the
$\vec{k}$-Kirchoff-Routh Hamiltonian is an adaptation of an argument in \cite{MaW} based on the analysis
of the vectorial Pohozaev-identity.\\

For $\om$ convex, by a well-known moving plane argument we have that $\Sigma \subset \om$. Remark that, by a result in \cite{GrT}, if $\om$ is convex then the
$\vec{k}$-Kirchoff-Routh Hamiltonian, $k=(z_1,\cdots,z_m)$, has no critical points as far as $m\geq2$, whence the spikes set in this case is
a singleton, $z_1\in \om$ and is a critical point of the Robin function of the domain.
However, by a result in \cite{CT}, if $\om$ is convex then the Robin function is also convex and admits a unique critical point,
which is the unique harmonic center of $\om$. Therefore, as a corollary of Theorem \ref{thm:spike-vn.intro} and of the results in
\cite{GrT} and \cite{CT} we have,

\bco Let $v_n$ be a sequence of solutions of~\eqref{w5.1.intro}, \eqref{mass.intro.1} and assume that $\om$ is convex.
Then the conclusions of Theorem \ref{thm:spike-vn.intro} holds and in particular the spikes set is a singleton $\Sigma=\{z_1\}\subset \om$
and $z_1$ is the unique harmonic center of $\om$.
\eco

Actually, the almost explicit expression of the spike sequence solving \eqref{w5.1.intro} in case $\om$ is a ball is provided in Appendix B.

\brm {\it Interestingly enough, while we was completing this paper, we have been informed by P. Cosentino and F. Malizia
that they were able to show {\rm (}\cite{CM}{\rm )} that the interior spike points of Theorem \ref{thm:spike-vn.intro} are simple.}
\erm

\brm{\it
In particular we deduce that whenever $m=1$ and the spike point is an interior point, then
it must coincide with an harmonic center of $\om$. At last, remark that,
$$
 1\leq m\leq |\mathcal{Z}| \leq \frac{C_0}{M_{p,0}},
$$
(where $C_0$ has been defined in \rife{mass.intro.1}) and then in particular $C_0 < 2M_{p, 0}$ is a sufficient condition to guarantee that the spikes set is a singleton, $m=1$.}
\erm

\brm {\it
The assumption \eqref{eq:assump-2.intro} is rather natural, as can be seen by
the underlying dual formulation arising from physical arguments relative to {\rm $\prl$}.
In fact, in terms of $\pn=\frac{|\an|}{\lmn}\vn$ and in view of {\rm \rife{mass.intro.1}}
(see Theorem \ref{thm7.2.intro} below) a sufficient condition to come up with
\eqref{eq:assump-2.intro} takes the form $\ino [\an+\lm_n \pn]_+^{p+1}\leq C_2|\an|$, which
is just an assumption about the growth of the associated entropy and free energy functionals,
see Appendix A for more details.\\
Moreover all the solutions in \cite{FW} and \cite{We} satisfy \eqref{eq:assump-2.intro}.}
\erm

At this point we are ready to apply these results to the
asymptotic description of solutions of~$\prl$, which yields a new concentration-quantization
phenomenon. Recall that
$$
\om_{n,+}:=\{x\in \om\,:\,\an+\lm_n\pn>0\},
$$ and ${\Sigma}_{r}\coloneqq \underset{i =1,\cdots,m}{\cup}B_{r} (z_{i})$, for some $r>0$.
Then we have,
\bte\label{thm7.2.intro} Let $\pn$ be a sequence of solutions of {\rm $\prl$} with $\lm_n\to +\ii$ and
\beq\label{hyp.intro.pn.1}
\frac{\lm_n}{|\an|^{1-\frac{p}{p_{_N}}}}\leq (C_0)^{\frac 2 N}\quad \mbox{ for some } C_0>0.
\eeq
Assume that, either
$$
 \om \;\; \mbox{is convex},
$$
or
\beq\label{hyp.intro.pn}
\ino [\an+\lm_n \pn]_+^{p+1}\leq C_2|\an|,
\eeq
for some $C_2>0$.\\
Then $\vn=\frac{\lmn}{|\an|}\pn$ satisfies \rife{vn.1.intro}, \rife{mass.intro.1},
\rife{eq:assump-2.intro} and consequently the conclusions of Theorem \ref{thm:spike-vn.intro}
hold for $\vn$ in $\om$ with~$\mathcal{Z}=\sum_{j=1}^m k_j z_j$. In particular, possibly along a subsequence, we have:
\begin{itemize}
 \item[(i)] the ``mass quantization identity'',
\beq\label{mass:quant}
\lim\limits_{n\to +\ii}\left(\frac{\lm_n}{|\an|^{1-\frac{p}{p_{_N}}}}\right)^{\frac N 2}=|\mathcal{Z}| M_{p,0},
\eeq
whence in particular
$$
\frac{|\an|}{\lm_n}=\frac{|\an|^{\frac{p}{p_{_N}}}}{(|\mathcal{Z}| M_{p,0})^{\frac2N}}(1+o(1))\to +\ii;
$$

 \item[(ii)] for any $r<R_0<R$ and any $n$ large enough,
\beq\label{round.intro}
\underset{i=1,\cdots,N_m}{\cup}B_{\eps_n r}(z_{i,n})\Subset\om_{n,+}\Subset
\underset{i=1,\cdots,N_m}{\cup}B_{\eps_n R}(z_{i,n}),
\eeq
\beq\label{psi:est.1}
\frac{|\an|}{\lm_n}\leq \psi_n(x)\leq \left(1+ R^{\frac{2}{1-p}}_0 u(0)\right)\frac{|\an|}{\lm_n}, \quad x\in \om_{n,+},
\eeq
\beq\label{psi:est.2}
\psi_n(x)=\frac{(1+o(1))}{N_m}\sum\limits_{i=1}^{N_m}G(x,z_{n,i}),\quad x\in \om\setminus \Sigma_{r};
\eeq

\item[(iii)] for any $\phi \in C^{0}(\ov{\om}\,)$,
\beq\label{psi:est.0}
 \lim\limits_{n\to +\ii}\int\limits_\om [\an+\lm_n\pn]^p_+\phi=
 \frac{1}{N_m}\sum\limits_{j=1}^{m}k_j\phi(z_j).
\eeq

\end{itemize}

Moreover, if $\Sigma\subset \om$,
then $(z_1,\cdots,z_m)$ is a critical point of the
$\vec{k}$-Kirchoff-Routh, where $\vec{k}=(k_1,\cdots,k_m)$. In particular, if $\om$ is convex then $\Sigma=\{z_1\}$ and
$z_1\in\om$ is the unique harmonic center of $\om$.
\ete

\brm {\it
Remark that,
in view of $\frac{|\an|}{\lm_n}\to +\ii$,
$\psi_n$ blows up uniformly in $\om_{n,+}$.}
\erm

\bdf{\it  Any sequence of solutions of {\rm $\prl$} with $\an<0$, $\lm_n\to +\ii$ and such that $\vn=\frac{\lm_n}{\an}\pn$
is a $\mathcal{Z}$-spikes sequence of \eqref{w5.1.intro} is by definition said to be a {\bf $\mathcal{Z}$-blowup sequence}
(or a blowup sequence for short)
of {\rm $\prl$}. The set $\Sigma$ is said to be
the {\bf blowup set} and each spike point $z\in \Sigma(\mathcal{Z})$ is said to be a {\bf blowup point}. Any blowup point is said to be {\bf simple}
if it is simple as a spike point according to Definition \ref{multi:spike}.
Obviously any blowup sequence satisfies
{\rm \eqref{mass:quant}, \eqref{round.intro},
\eqref{psi:est.1}, \eqref{psi:est.2} and \eqref{psi:est.0}}. }
\edf

Interestingly enough, this concentration-compactness phenomenon shares some features typical of critical problems (Liouville $N=2$, Yamabe $N\geq 3$)
in a subcritical context ($p<p_{_N}$), the singularity of the limit being induced in this case by $\lm_n\to +\ii$ and \eqref{hyp.intro.pn.1}.\\

For the model case of a blow up sequence we refer to solutions of $\fbi/\prl$ on balls as described in Appendix B. However,
by using Theorem \ref{thm7.2.intro} together with those in \cite{FW,We},
we deduce the corresponding existence results for blowup sequences of $\prl$. This is immediate for \cite{We} while some care
is needed for~\cite{FW}.

\bco[\cite{FW,We}]\label{exist:plasma} Let $\om\subset \R^N$, $N\geq 3$ be a smooth and bounded domain.\\
{\rm (a)-(\cite{FW})}  Let $\vn$ be the mountain pass solutions of \eqref{vn.1.intro} with $\vn=0$ as derived in \cite{FW}.
Then $\vn$ is a  $\mathcal{Z}$-spikes sequence and $\pn=\frac{|\an|}{\lm_n}\vn$ is a
$\mathcal{Z}$-blowup sequence of {\rm $\prl$} relative to $\mathcal{Z}\equiv \{z_1\}\subset \om$ in $\om$, where $z_1$ is an harmonic center of $\om$.
In particular $\vn$ satisfies \eqref{est:nm2};\\
{\rm (b)-(\cite{We})}  Assume that $\mathcal{H}(x_1,\cdots,x_m;\vec{1}_m)$ admits a nondegenerate critical point $(z_1,\cdots,z_m)$, and let
$\vn$ be the solutions of \eqref{vn.1.intro} with $\vn=0$ as derived in \cite{We}. Then $\vn$ is a $\mathcal{Z}$-spikes sequence with~$\mathcal{Z}=\sum_{j=1}^m z_j$ and $\pn=\frac{|\an|}{\lm_n}\vn$ is a $\mathcal{Z}$-blowup sequence of {\rm $\prl$} relative to $\Sigma= \{z_1,\cdots,z_m\}\subset \om$.
\eco

We review in Appendix A the definition of variational solutions of $\prl$ introduced in~\cite{BeBr}. Recall that at least one variational solution of $\prl$ exists for any $\lm>0$.
Note that (see Theorem 1.2 in \cite{BHJY}) $\al<0$ for any solution of $\prl$ with $\lm$ large enough, which was first proved in \cite{BaMa} for variational solutions.

\bte\label{thm:varsol} Let $(\an,\pn)$ be a sequence of variational solutions such that $\lm_n\to +\ii$ and assume that $|\an|\to+\ii$.
Then $(\an,\pn)$ satisfies \eqref{hyp.intro.pn} and in particular,
$$
\limsup\limits_{n\to +\ii}|\an|^{-1}\ino [\an+\lm_n \pn]_+^{p+1}\leq \frac{p+1}{p-1}.
$$
\ete

\bigskip

Next we derive some general results about solutions of $\prl$.\\
The first one is a sort of analogue of the a priori estimates in \cite{yy} about mean field equations of Liouville type.
\bte\label{thm.1.21.intro} For any $0<\epsilon<\frac12$ let us define
$$
I_\epsilon=\left\{(0,M_{p,0}-\eps)\right\}\cup \left\{\sg>0\,:\, \sg \leq \frac{1}{\epsilon},\;M_{p,0}(k+\epsilon) <\sg<M_{p,0}(k+1-\epsilon),\;k \in\N \right\}
$$
and, for any fixed solution $(\al,\psi)$ of {\rm $\prl$} with $\al<0$,
$$
\sg=\sg(\al,\lm):=\left(\frac{\lm}{|\al|^{1-\frac{p}{p_{_N}}}}\right)^{\frac N 2}.
$$
There exists $\ov{\lm}_\eps$, $\ov{\al}_\eps$ and $\ov{C}_\eps$ such that $\lm\leq \ov{\lm}_\eps$,
$\all\geq \ov{\al}_\eps $ and $\|\pl\|_\ii\leq \ov{C}_\eps$ for any solution of {\rm $\prl$} such that $\sg\in I_\epsilon$.
\ete

At last, we can say something more about Lemmas \ref{altype2.intro}, \ref{altype3}. Here $M_{p+1,0}$ is defined in \eqref{entire.mp10} and we recall that
in view of the mass quantization \eqref{mass:quant} we have $\frac{|\an|}{\lm_n}=\frac{|\an|^{\frac{p}{p_{_N}}}}{|\mathcal{Z}| M_{p,0}}(1+o(1))\to +\ii$.
\bco\label{altype3.not} Let ${\rm }(\an,\pn)$ be a~$\mathcal{Z}$-blowup sequence with $\mathcal{Z}=z$ and~$z\in\om$.
Then
$$
\int\limits_{\om_{+,n}}|\nabla \pn|^2=\frac{|\an|}{2\lm_n}\left(\frac{M_{p+1,0}}{M_{p,0}}+o(1)\right),\; \mbox{as }n\to +\ii.
$$
In particular, for any such sequence the following estimate holds for $\mu_{n,-}$ in Lemma \ref{altype3},
$$
\mu_{n,-}\geq\frac{4N^2\omega_{\sscp N}R_0^{N-2}}{(p+1)}\left(\frac{M_{p+1,0}}{M^2_{p,0}}+o(1)\right),\; \mbox{as }n\to +\ii.
$$
\eco

As remarked above, this shows in particular that in dimension $N\geq 3$ the uniform estimate in the claim of Lemma \ref{altype2.intro} is false.
It would be interesting to find out which is the best possible universal lower bound for $\mu_{-}$ in Lemma \ref{altype3}.

\bigskip

We conclude the introduction with an open problem and a conjecture.

\medskip

{\bf OPEN PROBLEM} Let $\vn$ be a spikes sequence of either \eqref{w5.1.intro} or of \eqref{vn.1.intro} with $\vn=0$ on $\pa \om$ where
$\om$ is not convex. Prove or disprove that $\Sigma \subset \om$.

\medskip

Concerning this point, we remark that, besides the above mentioned results in \cite{LP},  when replacing $\Delta$
with the original Grad-Shafranov operator, Caffarelli-Friedman (\cite{CF}) found for $p=1$ that the spikes approach the boundary
asymptotically as $\lm\to+\ii$.

\medskip

Finally we have the following,\\
{\bf CONJECTURE} Let $(\an,\pn)$ be a sequence of variational solutions such that $\lm_n\to +\ii$.
Then $(\an,\pn)$ satisfies \eqref{hyp.intro.pn.1}.

\medskip
We support this conjecture by the explicit evaluation of $\all$ in the case where $\om$ is a ball in dimension $N\geq 3$, see Remark \ref{rem9.4}.
In fact in Theorem \ref{thm9.2} of Appendix B we generalize an argument worked out for $N=2$ in \cite{BSp} and prove that for $N\geq 3$
and $\om$ a ball, $\fbi/\prl$ admit a unique solution
for any positive $I/\lm$. In particular we prove that these solutions satisfy \eqref{hyp.intro.pn.1}. However, since a
variational solution always exists, by the uniqueness we answer the conjecture in the affirmative at least in this case. Let us shortly make some
comment about the relevance of the conjecture.
Let $(\an,\pn)$ be a sequence of variational solutions such that $\lm_n\to +\ii$. By a result
in \cite{BeBr} the plasma region of $(\an,\pn)$ is connected. Thus, if the conjecture were proved in the affirmative, by
Theorems \ref{thm7.2.intro} and \ref{thm:varsol} we would deduce that $\pn$ is a $z_1$-blowup sequence of {\rm $\prl$},
that the unique blow up point $\{z_1\}\equiv \Sigma$ is simple and that, in view of \eqref{round.intro}, the plasma region
is ``asymptotically a round point''. In particular for $\om$ convex then $z_1\in \om$ would coincide with the unique harmonic center of $\om$.
This would be a generalization/improvement of classical results concerning
the shape and the size of the plasma region for variational solutions as first pursued
in \cite{CF,Shi} for $p=1$ and then for $p>1$ in \cite{BMar2,BSp,BeBr}. For example, concerning this point, it was
shown in \cite{BMar2} that the diameter of the orthogonal projection of the plasma region vanishes as $\lm \to +\ii$.\\

\bigskip

This paper is organized as follows. In section \ref{sec2} we prove Theorem \ref{thm1.1}.
In section \ref{sec3} we prove that $\frac1p\Lambda (\om,2p)<\lm_0(\om,p)$ for any $p$ large.
In section \ref{sec5} we prove Lemma \ref{altype3}.
In section \ref{sec5.1} we prove Theorem \ref{thm1.5}.
In section \ref{sec6} we prove a Brezis-Merle type result for singularly perturbed subcritical problems in dimension $N\geq 3$, which is the starting point of the spikes analysis in Theorems \ref{thm5.7} and its refinement stated above, Theorem \ref{thm5.7.intro}.
In particular we prove the so called "Vanishing Lemma" (see Lemma \ref{vanishle}).
In section \ref{sec7} we prove a boundary version of the Vanishing Lemma (see Lemma \ref{vanishle-bdy}) and Theorems \ref{thm:spike-vn.intro} and \ref{thm7.2.intro}.
The section is concluded with the proofs of Corollary \ref{exist:plasma}, Theorem \ref{thm:varsol} and Corollary \ref{altype3.not}.
Few useful details about variational solutions (\cite{BeBr}) as well as the uniqueness of solutions and some properties of
spike sequences on balls are described in the appendices.

\medskip

\subsection*{Data availability statement}
 Data sharing not applicable to this article as no datasets were generated or analysed during the current study.
 
\medskip 

\section{The proof of Theorem \ref{thm1.1}}\label{sec2}
In this section we prove Theorem \ref{thm1.1}.
We need a preliminary estimate, namely Theorem \ref{altype2.intro}, which was proved in \cite{BJ3} for non-negative solutions of $\prl$.
Actually we will use the following Lemmas \ref{enrg0}, \ref{altype2} which, by the very definition of the energy,
$2\el=\ino |\nabla \pl|^2$, immediately implies Theorem \ref{altype2.intro}. The first Lemma is just a standard proof of \rife{enrg1.0}.
\ble\label{enrg0} Assume $\all<0$ and let
$\om_+=\{x\in\om\,:\,\all+\lm\pl>0\}$. Then
$$
2\el =\int\limits_{\om_+} |\nabla \pl|^2-\frac{\all}{\lm},
$$
and in particular setting $\om_{-}=\{x\in\om\,:\,\all+\lm\pl<0\}$ we have
$$
\int\limits_{\om_{-}} |\nabla \pl|^2=-\frac{\all}{\lm} \quad \mbox{ and }\quad \int\limits_{\{\all+\lm\pl=0\}} |\nabla \pl|^2=0.
$$

\ele
\proof
Let $t_0\coloneqq \frac{|\all|}{\lm}$, by the Sard Lemma we can take $t_n\searrow t_0^+$ such that
$\om(t_n)=\{x\in\om\,:\,\pl>t_n\}$ is smooth.
Therefore $\int\limits_{\om(t_n)}(\all+\lm\pl)^p=-\int\limits_{\pa \om(t_n)}\pa_{\nu}\pl$.
Multiplying $\prl$ by $\pl$ and integrating by parts, we find that,
$$
\int\limits_{\om(t_n)}(\all+\lm\pl)^p\pl=-t_n\int\limits_{\pa \om(t_n)}\pa_{\nu}\pl+
\int\limits_{\om(t_n)}|\nabla \pl|^2=t_n\int\limits_{\om(t_n)}(\all+\lm\pl)^p+
\int\limits_{\om(t_n)}|\nabla \pl|^2.
$$

At this point, passing to the limit as $t_n\searrow t_0^+$ we find that,
$$
\int\limits_{\om}|\nabla \pl|^2=\int\limits_{\om_+}[\all+\lm\pl]_+^p\pl\equiv
\int\limits_{\om_+}(\all+\lm\pl)^p\pl=-\frac{\all}{\lm}+\int\limits_{\om_+}|\nabla \pl|^2.
$$

Since it is well-known that
$$
\int\limits_{\{\al+\lm \pl=0\}}|\nabla \pl|^2=0,
$$
then we also deduce that $\int\limits_{\om_{-}} |\nabla \pl|^2=-\frac{\all}{\lm}$, as claimed.
\finedim

\ble\label{altype2} Let $N=2$ and $(\all, \pl)$ be a solution of {\rm $\prl$} with $\all<0$.
Then we have,
\beq\label{level1.2}
\int\limits_{\om_+} |\nabla \pl|^2\leq \frac{p+1}{8\pi},
\eeq

where the equality holds if and only if $\om_+$ is a disk and $\pl$ is radial in $\om_+$.
\ele
\proof

Let $\thl=\|\pl\|_{\ii}$ and recall that $\pl>0$ in $\om$ and $\all<0$ by assumption.
Set
$$
\om(t)=\{x\in\om\,:\,\pl>t\}, \quad \Gamma(t)=\{x\in\om\,:\,\pl=t\}, \quad t\in [0, \thl],
$$
and
$$
m(t)=\int\limits_{\om(t)}\left(\all+\lm \pl\right)^{p},\qquad \mu(t)=|\om(t)|,
\qquad {e}(t)=\int\limits_{\om(t)}|\nabla \pl|^2,
$$
where $|\om(t)|$ is the area of $\om(t)$. Let
$\om_{-}$ be as in Lemma \ref{enrg0}. Since $|\Delta\pl|$ is locally bounded below away from zero in
$\om_+$ and since $\pl$ is harmonic in $\om_{-}$, then it is not difficult to see that
the level sets have vanishing area $|\Gamma(t)|=0$ for any $t\neq t_0:=\frac{|\all|}{\lm}$.
Therefore $m(t)$ and $\mu(t)$ are continuous in $[0,t_0)$ and  $(t_0, \thl]$. However a closer inspection shows that
$m(t)$ is continuous in $[0,\thl]$ and in particular
by the co-area formula in \cite{BZ} that $m(t)$
is also locally absolutely continuous in $(0,t_0)$ and in
$(t_0,\thl)$ and that $e(t)$ is absolutely continuous in $(0,\thl)$.
We will use the fact that,
\beq\label{muem1}
m(0)=1, \quad \mu(0)=1, \quad e(0)=\int\limits_{\om}|\nabla \pl|^2,
\eeq
and
\beq\label{muem2}
m(\thl)=0, \quad \mu(\thl)=0, \quad e(\thl)=0.
\eeq
By the co-area formula and the Sard Lemma we have,
\beq\label{diffm2}
-m^{'}(t)=\int\limits_{\Gamma(t)}\frac{\left(\all+\lm \pl\right)^{p}}{|\nabla \pl |}=
\left(\all+\lm{t}\right)^{p}\int\limits_{\Gamma(t)}\frac{1}{|\nabla \pl |}=
\left(\all+\lm{t}\right)^{p}(-\mu^{'}(t)), \quad \mbox{for a.a.\,} t\in (t_0,\thl),
\eeq
while $m^{'}(t)= 0$ in $(0,t_0)$
and
\beq\label{diffem2}
m(t)=-\int\limits_{\om(t)}\Delta \pl=\int\limits_{\Gamma(t)}|\nabla \pl|=-e^{'}(t),
\eeq
for a.a. $t\in (t_0, \thl)$, while $m(t)=m(0)=1$ for $t\in [0,t_0]$.
By the Schwarz inequality and the isoperimetric inequality we find that,
\begin{align}\nonumber
 -m^{'}(t)m(t)
 =&
\int\limits_{\Gamma(t)}\frac{\left(\all+\lm{\pl}\right)^{p}}{|\nabla \pl |}\int\limits_{\Gamma(t)}|\nabla \pl|=
\left(\all+\lm{t}\right)^{p}\int\limits_{\Gamma(t)}\frac{1}{|\nabla \pl |}\int\limits_{\Gamma(t)}|\nabla \pl| \\ \nonumber
\geq & \left(\all+\lm {t}\right)^{p}\left(|{\Gamma(t)}|_{1}\right)^2 \\ \nonumber
\geq & \left(\all+\lm {t}\right)^{p}
4\pi \mu(t),
\mbox{ for a.a. } t\in (t_0, \thl),
\end{align}
where $|{\Gamma(t)}|_{1}$ denotes the length of $\Gamma(t)$.
It follows that
\beq\label{diff12}
\frac{(m^2(t))^{'}}{8\pi }+
\left(\all+\lm{t}\right)^{p}\mu(t)\leq 0,\mbox{ for a.a. } t\in (t_0, \thl).
\eeq
By using the following identity,
$$
\left(\all+\lm {t}\right)^{p}\mu(t)=\frac{1}{\lm(p+1)}
\left(\left( \all+\lm {t}\right)^{p+1}\mu(t)\right)^{'}
-\frac{1}{\lm (p+1)}\left(\all+\lm {t}\right)^{p+1}\mu^{'}(t),
$$
$\mbox{ for a.a. } t\in (t_0, \thl)$, together with \rife{diff12} and \rife{diffm2} we conclude that,
$$
\left(\frac{m^2(t)}{8 \pi}+
\frac{\left( \all+\lm {t}\right)^{p+1}}{\lm (p+1)}\mu(t) \right)^{'}
-\frac{1}{\lm (p+1)}\left( \all+\lm{t}\right)m^{'}(t)\leq 0,
$$
$\mbox{ for a.a. } t\in (t_0, \thl)$.
Thus we see that,
\beq\label{diff22}
-\frac{m^2(t)}{8\pi}-\frac{\left( \all+\lm {t}\right)^{p+1}}{\lm (p+1)}\mu(t)
-\frac{1}{\lm (p+1)}\int\limits_{t}^{\thl}ds \,\left( \all+\lm {s}\right)
m^{'}(s)\leq 0,\; \forall\,t\in (t_0,\thl).
\eeq
In view of \rife{diffem2} and $-\lm t_0=\all$, recalling that
$m(t_0)=m_+(t_0)=\lim\limits_{t\to t_0^+}m(t)$ and $e(t_0)=e_+(t_0)=\lim\limits_{t\to t_0^+}e(t)$, we have
\begin{align}\nonumber
 \int\limits_{t_0}^{\thl}ds \,\left( \all+\lm {s}\right)m^{'}(s)
 =&-\all m(t_0)+
\lm\int\limits_{t_0}^{\thl}ds\,s\, m^{'}(s)=-\all m(t_0)-\lm t_0m(t_0)-
\lm\int\limits_{t_0}^{\thl}ds m(s)\\ \nonumber
=& -\lm\int\limits_{t_0}^{\thl}ds m(s)=\lm\int\limits_{t_0}^{\thl}ds e^{'}(s)=-\lm e_+(t_0)=-\lm e(t_0). \nonumber
\end{align}
Taking the limit of \rife{diff22} as $t\to t_0^+$ we find that
$$
-\frac{m^2_+(t_0)}{8\pi}+\frac{e_+(t_0)}{p+1}=-\frac{m^2(t_0)}{8\pi}+\frac{e(t_0)}{p+1}\leq 0,
$$
that is,
$$
e(t_0)\leq \frac{p+1}{8\pi}m^2(t_0)=\frac{p+1}{8\pi},
$$
which proves \rife{level1.2}.
The equality holds in \rife{level1.2} if and only if $\Gamma(t)$ is a circle
for a.a. $t$, whence if and only if $\om_+$ is a disk as well and $\pl$ is radial in $\om_+$.
\finedim

\

Let us set
$R_{p+1}(w)=\dfrac{\ino |\nabla w|^2}{\left(\ino |w|^{p+1}\right)^{\frac{2}{p+1}}}$, for
$w\in H^1_0(\om)\setminus\{0\}$.
Then
$$
\Lambda(\om,p+1)=\inf\limits_{w\in H^{1}_0(\om)\setminus\{0\}}R_{p+1}(w).
$$
Recall also the definition of~$\lm_0(\om,p)$ from~\eqref{eq:lambda0}.
The following result implies Theorem~\ref{thm1.1}, generalizing
an estimate first proved in \cite{BJ3} for non-negative solutions as well as
the celebrated result by Temam in \cite{Te2} for $p=1$.
\bte\label{thm4.3}
Let $N=2$, $p\in [1,+\ii)$, $|\om|=1$ and $(\all, \pl)$ be a solution of {\rm $\prl$} with
$\all\leq 0$. Then  $\lm\geq \lm_{0}(\om,p)$ and the equality holds if and only if:\\
$(i)$ $\all=0$;\\
$(ii)$ either $p=1$ or $p>1$ and $\om$ is a disk of unit area;\\
$(iii)$ in both situations occurring in $(ii)$, $\xil:=\lm \pl$ is a minimizer of $R_{p+1}$.
\ete
\proof
The Dirichlet energy~$\el=\frac12\ino |\nabla \pl|^2$ can be written as follows,
$$
\el=\frac{1}{2\lm}\ino[\all+\xil]_+^p\xil.
$$
Hence we have
\begin{align}
 \lm=&\dfrac{\ino |\nabla \xil|^2}{\int_{\om_+}[\all+\xil]^p_+\xil} \nonumber \\
 =& \dfrac{\int_{\om} |\nabla \xil|^2}{\int_{\om_+} |\nabla \xil|^2}
 \dfrac{\int_{\om_+} |\nabla \xil|^2}{\left(\int_{\om_+} \xil^{p+1}\right)^{\frac{2}{p+1}}}
 \dfrac{\left(\int_{\om_+} \xil^{p+1}\right)^{\frac{2}{p+1}}}
{\left(\int_\om [\all+\xil]_+^{p}\xil\right)^{\frac{2}{p+1}}}
\frac{1}{\left(\int_\om [\all+\xil]_+^{p}\xil\right)^{\frac{p-1}{p+1}}} \nonumber\\ \nonumber
\geq& \dfrac{\int_{\om} |\nabla \xil|^2}{\int_{\om_+} |\nabla \xil|^2}
\Lambda(\om_+,p+1)
\left(\dfrac{\int_{\om_+} \xil^{p+1}}{\int_\om [\all+\xil]_+^{p}\xil}\right)^{\frac{2}{p+1}}
\frac{1}{\left(2\lm \el\right)^{\frac{p-1}{p+1}}} \\
=&\frac{\Lambda(\om_+,p+1)}
{\left(\lm \int_{\om_+} |\nabla \pl|^2\right)^{\frac{p-1}{p+1}}}
\dfrac{\lm^2\int_{\om} |\nabla \pl|^2}{\lm^2\int_{\om_+} |\nabla \pl|^2}
\left(\dfrac{\int_{\om_+} \xil^{p+1}}{\int_\om [\all+\xil]_+^{p}\xil}\right)^{\frac{2}{p+1}}
\frac{\left(\lm \int_{\om_+} |\nabla \pl|^2\right)^{\frac{p-1}{p+1}}}
{\left(\lm \int_{\om} |\nabla \pl|^2\right)^{\frac{p-1}{p+1}}} \nonumber \\
 = &\frac{\Lambda(\om_+,p+1)}{\left(\lm \int_{\om_+} |\nabla \pl|^2\right)^{\frac{p-1}{p+1}}}
(\mathcal{A}(\lm))^{\frac{1}{p+1}}, \nonumber
\end{align}
where
\begin{align}
 \mathcal{A}(\lm)
 \coloneqq &
\left(\dfrac{\lm^2\int_{\om} |\nabla \pl|^2}{\lm^2\int_{\om_+} |\nabla \pl|^2}\right)^{p+1}
\left(\dfrac{\int_{\om_+} \xil^{p+1}}{\int_\om [\all+\xil]_+^{p}\xil}\right)^{2}
\left(\frac{\lm\int_{\om_+} |\nabla \pl|^2}
{\lm\int_{\om} |\nabla \pl|^2}\right)^{p-1} \nonumber \\
=& \left(\dfrac{\int_{\om} |\nabla \pl|^2}{\int_{\om_+} |\nabla \pl|^2}\right)^{2}
\left(\dfrac{\int_{\om_+} \xil^{p+1}}{\int_\om [\all+\xil]_+^{p}\xil}\right)^{2}
\geq 1. \nonumber
\end{align}
Using the estimate in Lemma \ref{altype2}, we readily deduce that,
$$
\lm^{2p}\geq \left(\frac{8\pi}{p+1}\right)^{p-1}\Lambda^{p+1}(\om_+,p+1),
$$
that is, in particular by the well-known
monotonicity of $\Lambda(\om,p)$ w.r.t. $\om$ (\cite{CRa1}), we have,
\beq\label{lmstar}
\lm\geq \lm_0(\om,p),
\eeq
where the equality holds if and only if  $\all=0$, $\xil$ is a minimizer of $R_{p+1}$ and
$\om$ is a disk of unit area, as desired.
\finedim

\

{\bf Remark} We remark that, for~$p=1$, there is a simpler argument which goes as follows.
Using the definition of~$R_{p+1}$ and a similar strategy we have
\begin{align}
 \lm
 =&\dfrac{\ino |\nabla \xil|^2}{\int_{\om}[\all+\xil]^p_+\xil}
 =\dfrac{\ino |\nabla \xil|^2}{\left(\int_{\om_+} [\all+\xil]_+^{p}\xil\right)^{\frac{2}{p+1}}}
\dfrac{1}{\left(\int_\om [\all+\xil]_+^{p}\xil\right)^{\frac{p-1}{p+1}}}
   \nonumber \\
\geq& \dfrac{\ino |\nabla \xil|^2}{\left(\int_{\om_+} \xil^{p+1}\right)^{\frac{2}{p+1}}}
\dfrac{1}{\left(2\lm\el\right)^{\frac{p-1}{p+1}}} \nonumber \\
\geq &
\dfrac{\int_{\om} |\nabla \xil|^2}{\left(\int_{\om} \xil^{p+1}\right)^{\frac{2}{p+1}}}
\dfrac{1}{\left(2\lm\el\right)^{\frac{p-1}{p+1}}} =
\dfrac{R_{p+1}(\xil)}{\left(2\lm\el\right)^{\frac{p-1}{p+1}}} \nonumber \\
\geq & \dfrac{\Lambda(\om,p+1)}{\left(2\lm\el\right)^{\frac{p-1}{p+1}}}. \nonumber
\end{align}
In particular, following the inequalities used so far it is readily seen that,
\beq\label{lmest.1}
\lm^{2p}\geq \dfrac{\Lambda^{p+1}(\om,p+1)}{\left(2\el\right)^{{p-1}}},
\eeq
where the equality holds if and only if $\all=0$ and $\xil$ is a minimizer of $R_{p+1}$.
For $p=1$ this already implies the claim, which was well-known since the pioneering result
of Temam~\cite{Te2}.

\

\section{An estimate about best Sobolev constants}\label{sec3}

For~$N=2$ and $p\in [1,+\infty)$, we wish to compare
$$
 \mu(\Omega,p)\equiv \frac{1}{p}\Lambda(\Omega, 2p) \quad \mbox{ and } \quad \lambda_0(\Omega,p)= \parenthesis{\frac{8\pi}{p+1}}^{\frac{p-1}{2p}}\Lambda(\Omega, p+1)^{\frac{p+1}{2p}}.
$$

We will need the following result from \cite{CTZ}: if~$p>2$ it holds,
\begin{align}\label{ctz.1}
 \frac{4\pi}{|\Omega|^{\frac{2}{p}}}
 \frac{1}{\parenthesis{ \prod_{k=0}^{[\frac{p}{2}]-1} (\frac{p}{2}-k)}^{\frac{2}{p}}}
 \leq \Lambda(\Omega,p)
 \leq \frac{8\pi e}{p}\parenthesis{\pi d_{\Omega}^2}^{-\frac{2}{p}},
\end{align}
where
$$
 d_{\Omega}\coloneqq \sup\braces{R: B_R(x_0)\subset \Omega \mbox{ for some } x_0\in\Omega}.
$$

\bpr\label{lme:31}
There exists $p_0>1$ depending on $\om$ such that
$$
 \mu(\Omega, p) < \lambda_0(\Omega,p), \qquad \forall p>p_0.
$$
\epr
\proof
Note that for~$p>1$, we have~$2p>2$ and~$p+1>2$, so we can use \rife{ctz.1} in our setting.
Thus we have
\begin{align}\nonumber
 \parenthesis{\frac{\mu(\Omega,p)}{\lambda_0(\Omega,p)}}^{p}
 =&\frac{\Lambda(\Omega, 2p)^p}{p^p} \parenthesis{\frac{p+1}{8\pi}}^{\frac{p-1}{2}}
 \frac{1}{\Lambda(\Omega, p+1)^{\frac{p+1}{2}}}  \\
 \nonumber
 \leq & \parenthesis{\frac{8\pi e}{2p}}^p \frac{1}{\pi d_\Omega^2}
 \frac{1}{p^p} \parenthesis{\frac{p+1}{8\pi}}^{\frac{p-1}{2}}
 \parenthesis{\frac{|\Omega|^{\frac{2}{p+1}}}{4\pi}
 \parenthesis{\prod_{k=0}^{[\frac{p+1}{2}]-1}(\frac{p+1}{2}-k)}^{\frac{2}{p+1}}}^{\frac{p+1}{2}} \\
 \nonumber
 =& \frac{|\Omega|}{\pi d_\Omega^2}\frac{(4\pi e)^p}{p^{2p}}
 \frac{(p+1)^{\frac{p-1}{2}}}{2^{\frac{p-1}{2}} (4\pi)^p }
 \parenthesis{\prod_{k=0}^{[\frac{p+1}{2}]-1}(\frac{p+1}{2}-k)}\\
 \nonumber
 =&\frac{|\Omega|}{\pi d_\Omega^2} \frac{e^p}{2^{\frac{p-1}{2}}}
 \frac{(p+1)^{\frac{p-1}{2}}}{p^{2p}}\parenthesis{\prod_{k=0}^{[\frac{p+1}{2}]-1}(\frac{p+1}{2}-k)}\\
 \nonumber
 =& \frac{|\Omega|}{\pi d_\Omega^2} \cdot
 \sqrt{\frac{2}{p+1}}\cdot
 \parenthesis{ \frac{e}{\sqrt{p}}\sqrt{\frac{p+1}{2p}} }^p
 \cdot
 \frac{1}{p^p} \parenthesis{\prod_{k=0}^{[\frac{p+1}{2}]-1}(\frac{p+1}{2}-k)}
 \eqqcolon F(p)
\end{align}

Remark that (since $\min\limits_{x\in[1,2]}\Gamma(x)> \frac12$)
$$
\parenthesis{\prod_{k=0}^{[\frac{p+1}{2}]-1}(\frac{p+1}{2}-k)}=
\frac{\Gamma(\frac{p+1}{2}+1)}{\Gamma([\frac{p+1}{2}]-\frac{p+1}{2}+1)}\leq 2\Gamma\left(\frac{p+1}{2}+1\right),
$$
which readily implies
$$
 \lim_{p\to+\infty} F(p)=0.
$$
The conclusion follows.
\finedim

\

\section{Energy estimates for \texorpdfstring{$N\geq 3$}{N>=3}.}\label{sec5}
In this section we prove Lemma \ref{altype3}.\\
{\bf The proof of Lemma \ref{altype3}.}
Let $\thl=\|\pl\|_{\ii}$ and recall that $\pl>0$ in $\om$. We take the same notations as in Lemma \ref{altype2},
where if $\all=0$ we just set $\om_+\equiv \om$ and
$\om_{-}=\emptyset$. Here we set
$$t_0:=\graf{\frac{|\all|}{\lm} \mbox{ if } \all<0,\\ 0 \mbox{ if } \all\geq 0.}$$
As in Lemma \ref{altype2}, we see that
$m(t)$ and $\mu(t)$ are continuous in $[0, \thl]$ and locally absolutely continuous in $(0,t_0)$ and in
$(t_0,\thl)$. We will need \rife{muem1}, \rife{muem2} as well.

By the co-area formula and the Sard Lemma we have,
\beq\label{diffm}
-m^{'}(t)=\int\limits_{\Gamma(t)}\frac{\left(\all+\lm \pl\right)^{p}}{|\nabla \pl |}=
\left(\all+\lm{t}\right)^{p}\int\limits_{\Gamma(t)}\frac{1}{|\nabla \pl |}=
\left(\all+\lm{t}\right)^{p}(-\mu^{'}(t)), \quad \mbox{for a.a.\,} t\in (t_0,\thl),
\eeq
while obviously $m^{'}(t)= 0$ in $(0,t_0)$
and
\beq\label{diffem}
m(t)=-\int\limits_{\om(t)}\Delta \pl=\int\limits_{\Gamma(t)}|\nabla \pl|=-e^{'}(t),
\eeq
for a.a. $t\in (t_0, \thl)$, while $m(t)=m(0)=1$ for $t\in [0,t_0]$.
By the Schwarz inequality and the isoperimetric inequality we find that,
\begin{align}
 -m^{'}(t)m(t)
 =&
\int\limits_{\Gamma(t)}\frac{\left(\all+\lm{\pl}\right)^{p}}{|\nabla \pl |}\int\limits_{\Gamma(t)}|\nabla \pl|=
\left(\all+\lm{t}\right)^{p}\int\limits_{\Gamma(t)}\frac{1}{|\nabla \pl |}\int\limits_{\Gamma(t)}|\nabla \pl| \nonumber \\
\geq &
\left(\all+\lm {t}\right)^{p}\left(|{\Gamma(t)}|_{1}\right)^2  \nonumber \\
\geq &\left(\all+\lm {t}\right)^{p}
\left(N (\omega_{\sscp N})^{\frac 1 N} (\mu(t))^{1-\frac{1}{N}}\right)^2,
\mbox{ for a.a. } t\in (t_0, \thl), \nonumber
\end{align}
where $|{\Gamma(t)}|_{1}$ denotes the length of $\Gamma(t)$ and we used the isoperimetric inequality.
%\todoRuijun[inline]{What is~$\omega_{N-1}$ here? Should it be~$\omega_N$?}
Therefore we deduce that,
\beq\label{diff1}
\frac{(m^2(t))^{'}}{2N^2 (\omega_{\sscp N})^{\frac 2 N}}+
\left(\all+\lm{t}\right)^{p}(\mu(t))^{2-\frac{2}{N}}\leq 0,\mbox{ for a.a. } t\in (t_0, \thl).
\eeq
By using the following identity,
\begin{align}
 \left(\all+\lm {t}\right)^{p}&(\mu(t))^{2-\frac{2}{N}} \nonumber \\
 = &\frac{1}{\lm(p+1)}
\left(\left( \all+\lm {t}\right)^{p+1}(\mu(t))^{2-\frac{2}{N}}\right)^{'}
-\frac{2(N-1)}{\lm N(p+1)}\left(\all+\lm {t}\right)^{p+1}(\mu(t))^{1-\frac{2}{N}}\mu^{'}(t), \nonumber
\end{align}
$\mbox{ for a.a. } t\in (t_0, \thl)$, together with \rife{diff1} and \rife{diffm} we conclude that,
$$
\left(\frac{m^2(t)}{2N^2 (\omega_{\sscp N})^{\frac 2 N}}+
\frac{\left( \all+\lm {t}\right)^{p+1}}{\lm (p+1)} (\mu(t))^{2-\frac{2}{N}}\right)^{'}
-\frac{2(N-1)}{\lm N(p+1)}\left( \all+\lm{t}\right)(\mu(t))^{1-\frac{2}{N}}m^{'}(t)\leq 0,
$$
$\mbox{ for a.a. } t\in (t_0, \thl)$.
Therefore we see that, for any~$t\in (t_0,\thl)$,
\beq\label{diff2}
-\frac{m^2(t)}{2N^2 (\omega_{\sscp N})^{\frac 2 N}}-\frac{\left( \all+\lm {t}\right)^{p+1}}{\lm (p+1)}
(\mu(t))^{2-\frac{2}{N}}
-\frac{2(N-1)}{\lm N(p+1)}\int\limits_{t}^{\thl}ds \,\left( \all+\lm {s}\right)(\mu(s))^{1-\frac{2}{N}}
m^{'}(s)\leq 0.
\eeq

Next observe that
\begin{align}
 e(t)
 =&\int\limits_{\om(t)}|\nabla \pl|^2=-t \int\limits_{\Gamma(t)}|\nabla \pl|-\int\limits_{\om(t)}\pl\Delta \pl  \nonumber \\
 =& te^{'}(t)+\int\limits_{\om(t)} \pl(\all+\lm \pl)^p=-tm(t)+\int\limits_{\om(t)} \pl(\all+\lm \pl)^p, \nonumber
\end{align}
from which we conclude that
\begin{align}\label{etil1}
 (\all +\lm t)m(t)+\lm e(t)
 =&(\all +\lm t)m(t)-\lm tm(t)+\int\limits_{\om(t)} \lm \pl(\all+\lm\pl)^p \nonumber \\
 =&\all m(t)+\lm\int\limits_{\om(t)} \pl(\all+\lm \pl)^p \\
 =& \int\limits_{\om(t)}(\all+\lm \pl)^{p+1}. \nonumber
\end{align}

\bigskip

Let $\epsilon>0$ be any positive number and denote by $\opep$ an infinitesimal function of $\epsilon$.
For the last term in \rife{diff2}, after some integration by parts and using \rife{etil1}, we deduce that
\begin{align}
 &-\int\limits_{t}^{\thl}ds \,\left( \all+\lm{s}\right)(\mu(s))^{1-\frac{2}{N}} m^{'}(s)
=\opep
-\int\limits_{t}^{\thl-\epsilon}ds \,\left( \all+\lm{s}\right)(\mu(s))^{1-\frac{2}{N}} m^{'}(s) \nonumber \\
=& \opep +
(\all+\lm t)(\mu(t))^{1-\frac{2}{N}}m(t)
-\lm\int\limits_{t}^{\thl-\epsilon}ds(\mu(s))^{1-\frac{2}{N}}e^{'}(s) \nonumber \\
 & \qquad \qquad +
\frac{N-2}{N}\int\limits_{t}^{\thl-\epsilon}ds(\all+\lm s)(\mu(s))^{-\frac{2}{N}}\mu^{'}(s)m(s) \nonumber \\
=&\opep +(\all+\lm t)(\mu(t))^{1-\frac{2}{N}}m(t)+
\lm(\mu(t))^{1-\frac{2}{N}}e(t) \nonumber \\
& \qquad \qquad +
\frac{N-2}{N}\int\limits_{t}^{\thl-\epsilon}ds(\mu(s))^{-\frac{2}{N}}\mu^{'}(s)
\left(\lm e(s)+(\all+\lm s)m(s)\right)\nonumber \\
=&\opep +(\all+\lm t)(\mu(t))^{1-\frac{2}{N}}m(t)+
\lm(\mu(t))^{1-\frac{2}{N}}e(t) \nonumber \\
 &\qquad\qquad +
\frac{N-2}{N}\int\limits_{t}^{\thl-\epsilon}ds(\mu(s))^{-\frac{2}{N}}\mu^{'}(s)
\int\limits_{\om(s)}(\all+\lm \pl)^{p+1}. \nonumber
\end{align}
Therefore we can pass to the limit as $\eps\to 0$ and substitute in \rife{diff2} to deduce that
\begin{align}\label{eqenrg3.1}
& -\frac{m^2(t)}{2N^2 (\omega_{\sscp N})^{\frac 2 N}}-\frac{\left( \all+\lm {t}\right)^{p+1}}{\lm (p+1)}
(\mu(t))^{2-\frac{2}{N}} \nonumber \\
& +
\frac{2(N-1)}{\lm N(p+1)}\left(
(\all+\lm t)(\mu(t))^{1-\frac{2}{N}}m(t)+
\lm(\mu(t))^{1-\frac{2}{N}}e(t)\right) \nonumber \\
&+\frac{2(N-1)}{\lm N(p+1)}\left(
\frac{N-2}{N}\int\limits_{t}^{\thl}ds(\mu(s))^{-\frac{2}{N}}\mu^{'}(s)
\int\limits_{\om(s)}(\all+\lm \pl)^{p+1}\right)\leq 0,\; \forall\,t\in (t_0,\thl)
\end{align}

\bigskip
\bigskip

At this point we split the proof in two cases.

\

{\bf Proof of $(a)$.}\\ In this case we have $\all\geq 0$, whence $t_0=0$, $\om_+\equiv \om$ and
$\om_{-}=\emptyset$. Clearly $\ino (\all+\lm\pl)^{p+1}=\all+2\lm \el$ and, since $\mu^{'}(s)\leq 0$, we deduce that,
\begin{align}
  \frac{N-2}{N}\int\limits_{0}^{\thl}ds(\mu(s))^{-\frac{2}{N}}\mu^{'}(s)
\int\limits_{\om(s)}(\all+\lm \pl)^{p+1}
\ge & (\all+2\lm \el)\frac{N-2}{N}
\int\limits_{0}^{\thl}ds(\mu(s))^{-\frac{2}{N}}\mu^{'}(s)
\nonumber \\
=& -(\all+2\lm \el). \nonumber
\end{align}
Therefore there exists $\widetilde{\beta_+} \in (0,1)$ such that
$$
0>\frac{N-2}{N}\int\limits_{0}^{\thl}ds(\mu(s))^{-\frac{2}{N}}\mu^{'}(s)
\int\limits_{\om(s)}(\all+\lm \pl)^{p+1}=-\widetilde{\beta_+} (\all+2\lm \el).
$$
Thus, since $t_0=0$, then  putting $\beta_+=1-\widetilde{\beta_+}\in (0,1)$, we deduce
from \rife{eqenrg3.1} that
\begin{align}
-\frac{1}{2N^2 (\omega_{\sscp N})^{\frac 2 N}}-\frac{\all^{p+1}}{\lm (p+1)}+
\frac{2(N-1)}{\lm N(p+1)}\left(\all+2\lm \el\right)
\leq
\widetilde{\beta_+}
\frac{2(N-1)}{\lm N(p+1)}\left(\all+2\lm \el\right),
\end{align}
% $$
% -\frac{2(N-1)}{\lm N(p+1)}\left(
% \frac{N-2}{N}\int\limits_{0}^{\thl}ds(\mu(s))^{-\frac{2}{N}}\mu^{'}(s)
% \int\limits_{\om(s)}(\all+\lm \pl)^{p+1}\right)=
% $$
that is,
$$
-\frac{1}{2N^2 (\omega_{\sscp N})^{\frac 2 N}}-\frac{\all^{p+1}}{\lm (p+1)}+
\beta_+\frac{2(N-1)}{\lm N(p+1)}\left(\all+2\lm \el\right)\leq 0,
$$
which is the same as
$$
2\lm \el\leq \frac{N}{2\beta_+ (N-1)}\frac{\lm (p+1)}{2N^2 (\omega_{\sscp N-1})^{\frac 2 N}}+
\frac{N}{2\beta_+ (N-1)}\all\left(\all^{p}-\beta_+\frac{2(N-1)}{ N}\right)
$$
which immediately implies \rife{n3.1}.

\bigskip

{\bf Proof of $(b)$.}\\
In this case we have $\all<0$, whence $t_0=\frac{|\all|}{\lm}$ and we recall that
$\mu(t_0)=|\om_+|, m(t_0)=1, \all+\lm {t_0}=0$ and, by definition of $\el$,
$$
\int\limits_{\om_+} (\all+\lm\pl)^{p+1}=\all+2\lm \el.
$$
Since $\mu^{'}(s)\leq 0$, we have
$$
\frac{N-2}{N}\int\limits_{t_0}^{\thl}ds(\mu(s))^{-\frac{2}{N}}\mu^{'}(s)
\int\limits_{\om(s)}(\all+\lm \pl)^{p+1}>
$$
$$(\all+2\lm \el)\frac{N-2}{N}
\int\limits_{t_0}^{\thl}ds(\mu(s))^{-\frac{2}{N}}\mu^{'}(s)=-|\om_+|^{1-\frac2N}(\all+2\lm \el),
$$
and consequently there exists $\widetilde{\beta_-}\in (0,1)$ such that
$$
0>\frac{N-2}{N}\int\limits_{t_0}^{\thl}ds(\mu(s))^{-\frac{2}{N}}\mu^{'}(s)
\int\limits_{\om(s)}(\all+\lm \pl)^{p+1}=-\widetilde{\beta_-}|\om_+|^{1-\frac2N}(\all+2\lm \el).
$$
Therefore we deduce that
$$
-\frac{1}{2N^2 (\omega_{\sscp N})^{\frac 2 N}}+
\frac{2(N-1)}{\lm N(p+1)}\left(\lm |\om_+|^{1-\frac{2}{N}}e(t_0)\right)\leq
\frac{2(N-1)}{\lm N(p+1)}\widetilde{\beta_-}|\om_+|^{1-\frac2N}(\all+2\lm \el)
$$

which, in view of Lemma \ref{enrg0}, i.e. $\lm e(t_0)=\all+2\lm \el$, and putting
$\beta_-=1-\widetilde{\beta_-}\in (0,1)$, implies that
$$
-\frac{1}{2N^2 (\omega_{\sscp N})^{\frac 2 N}}+
\beta_- \frac{2(N-1)}{\lm N(p+1)}|\om_+|^{1-\frac{2}{N}}\left(\all+2\lm \el\right)\leq 0
$$

which is the same as

$$
\el\leq
\frac{(p+1)}{4N^2 (\omega_{\sscp N})^{\frac 2 N}}\frac{N}{2(N-1)}\frac{1}{|\om_+|^{1-\frac{2}{N}}\beta_-}
+ \frac{|\all|}{2\lm},
$$
which immediately implies \rife{n3.2}.
\finedim
\bigskip

We remark that, in the case~$\alpha_\lambda<0$, combining the estimate for~$\|\psi_\lambda\|_{L^\infty}$ from next section and that
\begin{align}
 1=\int_{\Omega_+ } [\alpha_\lambda + \lambda\psi_\lambda]_+^p
 < \lambda^p \|\psi_\lambda\|_{L^\infty}^p |\Omega_+|,
\end{align}
we can obtain a lower bound of~$|\Omega_+|$ and hence a uniform upper bound for~$E_\lambda$
independent of~$|\Omega_+|$. But this is far from being sharp.

\

\section{Uniqueness for \texorpdfstring{$N\geq 3$}{N>=3}.}\label{sec5.1}
In view of Theorem A, the next result immediately implies Theorem \ref{thm1.5}.
\bte\label{thm4.1} Let $N\geq 3$, $p\in [1,p_{_N})$, $|\om|=1$ and $(\all, \pl)$ be a solution of {\rm $\prl$} with
$\all\leq 0$.
Fix any $s\in(p,p_{_N})$. Then
$$
\left\|\pl\right\|_{\ii}\leq  \frac{\lambda^{\frac{p}{s-p}}}{[N(N-2)\omega_N^{2/N}]^{\frac{s}{s-p}} [N(1-\frac{s}{p_{_N}})]^{\frac{1}{s-p}}},
$$
and
$$
\lm> \nu_s(\om,p)=
\left(N(N-2)\omega_N\right)^{\frac{s(p-1)}{p(2s-p-1)}}
\left(N(1-\frac{s}{p_{_N}})\right)^{\frac{p-1}{p(2s-p-1)}}
\Lambda(\Omega,p+1)^{(p+1)\frac{s-p}{p(2s-p-1)}}.
$$
\ete
\proof For fixed $\lm>0$, since $p<p_{_N}$, then it is well-known (\cite{BeBr,BJ1}) that $\left\|\pl\right\|_{\ii}\leq C_\lm$,
whence in particular $\rl=[\alpha_\lambda+\lambda\psi_\lambda(y)]_+^p\in L^{\ii}(\om)$.
Let~$G(x,y)$ denote the Green function on $\om$ with Dirichlet boundary conditions.
Then
\begin{align}\nonumber
 \psi_\lambda(x)
 =\int_\Omega G(x,y)[\alpha_\lambda+\lambda\psi_\lambda(y)]_+^p dy
\end{align}
and consequently for any $s\in (1,p_{_N})$ we have
\begin{align}\label{eq:pointwise estimate for psi}
 |\psi_\lambda(x)|
 \leq & \|G(x,y)\|_{L^s(\Omega, dy)} \|[\alpha_\lambda+\lambda\psi_\lambda(y)]_+^p\|_{L^t(\Omega)}
\end{align}
where~$\frac{1}{s}+\frac{1}{t}=1$.
Note that
\begin{align}
 G(x,y)=\frac{1}{N(N-2)\omega_N} \frac{1}{|x-y|^{N-2}} + H(x,y)\ge 0, \nonumber
\end{align}
and~$H(x,y)$ is smooth and nonpositive, and
\begin{align}
 |G(x,y)|\leq \frac{1}{N(N-2)\omega_N} \frac{1}{|x-y|^{N-2}}.
\end{align}
Thus, with~$|B_{R_N}|=|\Omega|=1$ (namely~$R_N=\omega_N^{-1/N}$), by a well known rearrangement argument (\cite{Bandle}) we have
\begin{align}
 \int_\Omega |G(x,y)|^s dy
 \leq &\frac{1}{(N(N-2)\omega_N)^s} \int_{\Omega} \frac{1}{|x-y|^{s(N-2)}} dy  \nonumber  \\
 \leq & \frac{1}{(N(N-2)\omega_N)^s} \int_{B_{R_N}(0)} \frac{1}{|y|^{s(N-2)}} dy \\
 =&\frac{|\Omega|^{1-\frac{s}{p_{_N}}}}{\omega_N^{s-\frac{s}{p_{_N}}} N^{s+1} (N-2)^s (1-\frac{s}{p_{_N}})}\nonumber
\end{align}
Recall that~$|\Omega|=1$.
Therefore, for any~$x\in\Omega$,
\begin{align}\label{gs}
 \|G(x,y)\|_{L^s(\Omega)}
 \leq \frac{1}{\omega^{2/N} N(N-2)} \frac{1}{N^{1/s}(1-\frac{s}{p_{_N}})^{1/s}}\eqqcolon g(s).
\end{align}

As for the second factor in the RHS of~\eqref{eq:pointwise estimate for psi}, with
$\om_+=\{x\in\om\,:\,\all+\lm\pl>0\}$ and $\theta_\lm=\left\|\pl\right\|_{\ii}$, we observe that,
\begin{align}\nonumber
 \|[\alpha_\lambda+\lambda\psi_\lambda(y)]_+^p\|_{L^t(\Omega)}
  =&\left(\int_\Omega [\alpha_\lambda+\psi_\lambda]_+^{pt} \right)^{\frac{1}{t}}  \nonumber\\
  \leq& \left(\int_{\Omega_+} (\alpha_\lambda+\lambda\psi_\lambda)^p (\alpha_\lambda+\lambda\psi_\lambda)^{pt-p}\right)^{\frac{1}{t}} \nonumber\\
  \leq & (\lambda\theta_\lambda)^{p(1-\frac{1}{t})}
  \left(\int_{\Omega_+} (\alpha_\lambda+\lambda\psi_\lambda)_+^p dx\right)^{\frac 1 t}
  \nonumber\\
  =&(\lambda\theta_\lambda)^{\frac{p}{s}}.\nonumber
\end{align}
Therefore
$
 |\psi_\lambda(x)|\leq g(s) (\lambda\theta_\lambda)^{\frac{p}{s}}
$
and taking supremum over~$x$ we obtain
\begin{align}\nonumber
 \theta_\lambda\leq g(s)(\lambda\theta_\lambda)^{\frac{p}{s}},
\end{align}
which, together with \rife{gs} implies, for $p<s<p_{_N}$,
\begin{align}\nonumber
 \left\|\pl\right\|_{\ii}
 =\theta_\lambda \leq g(s)^{\frac{s}{s-p}}\lm^{\frac{p}{s-p}}
 = \frac{\lambda^{\frac{p}{s-p}}}{[N(N-2)\omega_N^{2/N}]^{\frac{s}{s-p}} [N(1-\frac{s}{p_{_N}})]^{\frac{1}{s-p}}},
\end{align}
as claimed.
Therefore we also have an upper bound for the energy,
\begin{align}\nonumber
2 \el
 =\ino [\alpha_\lambda+\lambda\psi_\lambda]_+^p \psi_\lambda dx
 < \left\|\pl\right\|_{\ii}\int_\Omega [\alpha_\lambda+\lambda\psi_\lambda]_+^p dx =\left\|\pl\right\|_{\ii}
  \leq g(s)^{\frac{s}{s-p}}\lm^{\frac{p}{s-p}}.
\end{align}
At this point we argue as in the Remark right after the proof of Theorem \ref{thm4.3} to come up with
\rife{lmest.1} and then deduce that
$$
\lm^{2p}\geq \dfrac{\Lambda^{p+1}(\om,p+1)}{\left(2\el\right)^{{p-1}}}>
\left(N(N-2)\omega_N^{2/N}\right)^{\frac{s(p-1)}{s-p}}\left(N(1-\frac{s}{p_{_N}})\right)^{\frac{p-1}{s-p}}\lm^{-\frac{p}{s-p}(p-1)}{\Lambda^{p+1}(\om,p+1)},
$$
or equivalently
$$
\lm^{\frac{p(2s-(p+1))}{s-p}}
>
\left( N(N-2)\omega_N^{2/N}\right)^{\frac{s(p-1)}{s-p}}
\left(N(1-\frac{s}{p_{_N}})\right)^{\frac{p-1}{s-p}}
\Lambda^{p+1}(\om,p+1),
$$
which immediately implies the conclusion.\finedim

\

To illustrate the above estimate, one may extremize the right hand sides with respect to~$s\in (p, p_{_N})$.
This might be tedious and not very enlightening.
Here instead, taking~$s=\frac{1}{2}(p+p_{_N})$, we obtain that
\begin{align}\nonumber
 \|\psi_\lambda\|_{L^\infty}
 \leq \frac{ \lambda^{2\frac{\frac{p}{p_{_N}}}{1-\frac{p}{p_{_N}}}}}
 { \parenthesis{N(N-2)\omega_N^{2/N}}^{\frac{1+\frac{p}{p_{_N}}}{1-\frac{p}{p_{_N}}}}
 \parenthesis{\frac{N}{2}(1-\frac{p}{p_{_N}})^{\frac{2}{p_{_N}} \frac{1}{1-\frac{p}{p_{_N}}} }} },
\end{align}
and
\begin{align}\nonumber
 \lambda > \parenthesis{N(N-2)\omega_N}^{\frac{p-1}{p_{_N} -1} \frac{p_{_N}+p}{2p}}
 \parenthesis{\frac{N}{2}(1-\frac{p}{p_{_N}})}^{\frac{1}{p}\frac{p-1}{p_{_N} -1}}
 \Lambda(\Omega, p+1)^{\frac{p+1}{2p} \frac{p_{_N}-p}{p_{_N}-1}}.
\end{align}

\

\section{A Brezis-Merle type result for subcritical problems in dimension \texorpdfstring{$N\geq 3$}{N>2}}\label{sec6}
We deduce some properties of sequences of solutions of equations
arising in the study of $\prl$. A preliminary analysis, based on some arguments concerning the
``infinite mass" limit of mean field type equations as in \cite{BM},
 is needed concerning sequences of solutions of
\beq\label{wn.1}
\graf{-\Delta w_n = [w_n-1]_+^p\quad \mbox{in}\;\;\om,\\ \\
\|[w_n]_-\|_{1}\leq C.
}
\eeq
Here it is crucial that $p<p_{_N}$.
Unlike the remaining results in this section which are derived for $N\geq 3$ only, Theorem \ref{thm5.2} below holds in dimension $N\geq 2$.

\bte\label{thm5.2} Let $\om\subset \R^N$, $N\geq 2$ be an open set and $w_n$ be a sequence of solutions of \rife{wn.1} with $p\in (1,p_{_N})$.
Then $[w_n]_-$ is locally uniformly bounded and there exists a subsequence, still denoted $w_{n}$,
such that
\begin{itemize}
 \item[{\em either}] {\rm (i)} $ w_n \to +\ii$, locally uniformly in $\om$,

 \item[{\em or}~$\quad $] {\rm (ii)}  $w_n \to w$  in $C^{2}_{\rm loc}(\om)$.
\end{itemize}
More exactly {\rm (ii)} occurs if and only if $[w_n-1]_+^p$ is bounded in $L^{1}_{\rm loc}(\om)$.
\ete

\proof We argue as in Theorem 4 of \cite{BM}. By the Kato inequality (\cite{Ka}), we have
$$
\Delta [w_n]_-\geq -\Delta w_n \chi(\{w_n\leq 0\})= [w_n-1]_+^p\chi(\{w_n\leq 0\})=0,
$$
whence $[w_n]_-$ is a weakly subharmonic, nonnegative function, which is also uniformly bounded in $L^1(\om)$ and then it readily follows by the mean value inequality that it is also bounded in~$L^\ii_{\rm loc}(\om)$.
Let us peak any smaller open set $\om_0\Subset \om$
and consider the function
$u_n=w_n-c$, where $c=\inf\limits_{n\in\N, x \in \ov{\om_0}}w_n(x)>-\ii$. Then $u_n$ satisfies
$$
\graf{-\Delta u_n = [u_n-1+c]_+^p\quad \mbox{in}\;\;\om_0,\\ \\
u_n\geq 0 \mbox{ in }\om_0.
}
$$

{\bf CASE 1}. There exists a compact set $K\subset \om_0$ and a subsequence $u_n$
such that $\int\limits_K [u_n-1+c]_+^p\to +\ii$. Then (i) holds.\\
Take any compact set $K_0 \subset \om_0$ and let $G$ be the Green function with Dirichlet boundary
conditions in $\om_0$. Clearly $G(x,y)\geq d>0$ for any $(x,y)\in K \times K_0$, whence we deduce that
$$
u_n(x)\geq \int\limits_{\om_0} G(x,y)[u_n-1+c]_+^p\geq d\int\limits_{K} [u_n-1+c]_+^p\to +\ii,
\;\forall\,x\in K.
$$
Obviously, by the arbitrariness of $\om_0$, the same is true for $w_n(x)$ where $x$ is any point in a
compact subset of $\om$.

\bigskip

{\bf CASE 2}. $[u_n-1+c]_+^p$ and $u_n$ are bounded in $L^1_{\rm loc}(\om_0)$. Then (ii) holds.\\
Let $K\Subset \om_0$ be any compact subset and pick any open set $\omega\Subset \om_0$ such that
$K\Subset \omega$. Let $u_n=u_{1,n}+u_{2,n}$ where
$$
\graf{-\Delta u_{1,n} =[u_n-1+c]_+^p\quad \mbox{in}\;\;\omega,\\ \\
u_{1,n}= 0 \quad \mbox{ on } \pa\omega,
}
$$

$$
\graf{-\Delta u_{2,n} =0\quad \mbox{in}\;\;\omega,\\ \\
u_{2,n}\geq  0 \quad \mbox{ on }\pa \omega.
}
$$
Remark that by the maximum principle we have $u_{1,n}\geq 0$ and $u_{2,n}\geq 0$.
By classical estimates (\cite{stam}) $u_{1,n}$ is bounded in $W_0^{1,r}(\omega)$ for any
$r<\frac{N}{N-1}$, and then by the Sobolev embedding in $L^q(\omega)$ for any $q<p_{_N}$.
Since by assumption $u_n$ is bounded in $L^1(\omega)$, then $u_{2,n}=u_n- u_{1,n}$ is bounded in $L^1(\omega)$ and
by the mean value theorem also locally bounded in $\omega$, whence
passing to a subsequence $u_{2,n_k}\to u_2$ in $C^{2}_{\rm loc}(\omega)$ for some nonnegative and harmonic $u_2$ in $\omega$.\\
Clearly for any $q<p_{_N}$ there exists $u_1\in L^q(\omega)$ and a sub-subsequence $u_{1,n_k}$
such that $u_{1,n_k}\to u_1$ in $L^q(\omega)$.
Therefore $u_{n_k}\to u$ in $L^{q}_{\rm loc}(\omega)$ and since $p<p_{_N}$ then by classical elliptic estimates and a bootstrap argument
we deduce that $u_{1,n_k}\to u_1$ in $C^{2}_{\rm loc}(\omega)$.
Therefore $u_{n_k}\to u$ and obviously $w_{n_k}\to w$ in $C^{2}_{\rm loc}(\omega)$.

\bigskip

{\bf CASE 3}. $[u_n-1+c]_+^p$ is bounded in $L^1_{\rm loc}(\om_0)$. Then $u_n$ is bounded
in $L^1_{\rm loc}(\om_0)$ and by CASE 2 we have that (ii) holds.\\
In fact, since $u_n\geq 0$, for any compact set $K\Subset \om_0$, we obviously have
$$
\int\limits_{K}|u_n|=\int\limits_{K}u_n\leq |K||c-1|+\int\limits_{K}[u_n-1+c]_+\leq
 |K||c-1|+C_K\left(\int\limits_{K}[u_n-1+c]^p_+\right)^{\frac1p}.
$$
\finedim

\bigskip

From now on in the rest of this section we will be concerned with open sets $\om\subset \R^N$, $N\geq 3$.
\ble\label{lem5.3} Let $N\geq 3$, $p\in (1,p_{_N})$ and $w$ be a solution of
$$
\graf{-\Delta w =[w-1]_+^p\quad \mbox{in}\;\;\R^N,\\ \\
\bigints\limits_{\R^N} {\dsp [w-1]_+^p}<+\ii,\\ \\
w\geq 0 \mbox{ in } \R^N.}
$$
Then either $w\equiv c$ for some constant $c\in [0,1]$, or up to a translation $w= w_a$, for some $a\in [0,1)$
where $w_a$ is the unique solution of
\beq\label{entire.a}
\graf{-\Delta w =[w-1]_+^p\quad \mbox{in}\;\;\R^N,\\ \\
w(x)>a, \quad w(0)=\max\limits_{R^N}w>1,\\ \\
w(x)\to a \quad \mbox{ as}\quad |x|\to +\ii,}
\eeq
which takes the form \rife{entire.sol} below.
\ele
\proof
If $w\leq 1$ in $\R^N$ then $w$ would be harmonic and bounded, whence $w\equiv c$ for some $c\in[0,1]$. Otherwise
$w>1$ somewhere and we can argue as \cite[Theorem 2.1]{cli1} via a moving plane argument
to deduce that $w$ is radial and radially decreasing around some point $x_0\in\R^N$. After a translation we
may assume $x_0=0$. Remark that $w<1$ somewhere, otherwise $w_0=w-1$ would be a solution of
$\Delta v = v^p$ in $\R^N$, $v\geq 0$ in $\R^N$ which admits only the trivial solution (\cite{cli1}),
implying a contradiction to $\{w>1\}\neq \emptyset$. Let $R_a>0$ be defined by $w(R_a)=1$, then we readily deduce that
$$
w(r)=a+\left(\frac{R_a}{r}\right)^{N-2}(1-a), \;\; r\in [R_a,+\ii)\mbox{ and }
R_a w^{'}(R_a)=(2-N)(1-a).
$$

Next let us define $u(\varrho)$, $\varrho \in [0,1]$ as follows:
$$
w(r)-1=R_a^{\frac{2}{1-p}}u(R_a^{-1}r),\;\; r\in [0,R_a],
$$
then $u$ is the unique (\cite{gnn}) solution of
\beq\label{emden}
\graf{-\Delta u =u^p\quad \mbox{in}\;\;B_1(0),\\ \\
u>0 \mbox{ in } B_1(0),\\ \\
u= 0 \mbox{ on } \pa B_1(0).}
\eeq
The universal value $u^{'}(1)$ uniquely determines $R_a$ via the identity
$$
(2-N)(1-a)R_a^{-1}
= w^{'}(R_a)
=R_a^{\frac{2}{1-p}-1}u^{'}(1)<0,
$$
i.e.
\beq\label{entire.Ra}
R_a=\left(\frac{-u^{'}(1)}{(N-2)(1-a)}\right)^{\frac{p-1}{2}}.
\eeq
In other words for any $a\in [0,1)$ there exists a unique solution $w_a$ of \rife{entire.a}
which takes the form
\beq\label{entire.sol}
\graf{w_a(r)=1+R_a^{\frac{2}{1-p}}u(R_a^{-1}r),\;\; r\in [0,R_a],\\ \\
w_a(r)=a+\left(\frac{R_a}{r}\right)^{N-2}(1-a),\;\; r\in [R_a,+\ii).}
\eeq
\finedim

\bigskip

\brm\label{rmdgp}{\it
In view of Theorem 4.2 in \cite{DGP},
it is readily seen by the proof of Lemma \ref{lem5.3} that $x=0$ is a nondegenerate maximum point of $w_a(x)$.}
\erm

\

\brm\label{rmd12} {\it Since~$w_a$ is asymptotic to the constant~$a\geq 0$ at infinity, we see that $w_a\in D^{1,2}(\R^N)=\{w\colon \R^N\to \R \,\mid\,w\in L^{2^*}(\R^N),
|\nabla w|\in L^{2}(\R^N) \}$, $2^*=\frac{2N}{N-2}$,
if and only if $a=0$, in which case $w_0$ is said to be a ground state {\rm (\cite{FW})}.}
\erm

\

\brm\label{rem5.4}
{\it Let $w_a$ be the unique solution of \rife{entire.a},
and let us define,
$$
M_{p,a}:=\int\limits_{\R^N}[w_a-1]_+^p\equiv \int\limits_{B_{R_a}(0)}[w_a-1]_+^p \in (0,+\ii).
$$
Note that for the solution~$u$ of~\eqref{emden},
$$I_p=\int\limits_{B_{1}(0)}u^p=\int\limits_{B_1(0)} -\Delta u
=-\int\limits_{\pa B_1(0)} \frac{\pa u}{\pa \nu} = N\omega_N (-u'(1)). $$
Then, by a straightforward evaluation we find that
$$
M_{p,a}
=\frac{1}{R_a^{\frac{2}{p-1}-N+2}}I_p
=
\left(\frac{(N-2)(1-a)}{-u^{'}(1)}\right)^{{\frac{N}{2}(1-\frac{p}{p_{_N}})}} N\omega_N (-u'(1)),
$$
which is monotonic decreasing for $a\in [0,1)$ and, for any $1<p<p_{_N}$,
\begin{align}\label{entire.mp0}
M_{p,0}=N\omega_N (N-2)^{\frac{N}{2}({1-\frac{p}{p_{_N}}})} (-u'(1))^{1-\frac{N}{2}(1-\frac{p}{p_{_N}})}
\end{align}\
with~$1-\frac{N}{2}(1-\frac{p}{p_{_N}})>0$, while
$$ \lim_{a\to 1^-}M_{p,a}= 0.
$$
}
\erm

\brm{\it
  We also remark that, for the solutions~$w_a$ as above, the analogous quantity
  \begin{align}\label{entire.mp10}
   M_{p+1,a}:=\int\limits_{\R^N}[w_a-1]_+^{p+1} = \int\limits_{B_{R_a}(0)}[w_a-1]_+^{p+1},
  \end{align}
  can also be evaluated in terms of $u'(1)$.
  Indeed, a Pohozaev type argument shows that
  \begin{align}\nonumber
    I_{p+1}\coloneqq \int\limits_{B_1(0)} u^{p+1} = \frac{p+1}{(N+2)-p(N-2)}N\omega_N (u'(1))^2
  \end{align}
  and
  \begin{align}\nonumber
   M_{p+1, a}= \frac{1}{R_a^{\frac{4}{p-1}-N+2}} I_{p+1}.
  \end{align}
 }\erm

\bigskip

We will analyze the asymptotic behavior of solutions of
\beq\label{vn.1}
\graf{-\Delta \vn =\mu_n [\vn-1]_+^p\quad \mbox{in}\;\;\om,\\ \\
\mu_n \to +\ii,\\ \\
\vn\geq 0,}
\eeq
where $\om \subset \R^N$, $N\geq 3$, is any open set and $p\in(1,p_{_N})$.
Remark that the constant functions~$v_n=c_n\in [0,1]$ are trivial solutions of~\eqref{vn.1}.
These trivial solutions justify some of our later results.
Nonnegative harmonic functions taking values in~$[0,1]$ are also solutions, showing that~$v_n$ may not achieve local maximum in interior points.
More interesting solutions are those of the plasma problem, see next section.

\brm\label{rem5.2} {\it
By definition we assume that for each $n$, $[\vn-1]_+^p\in L^1(\om)$.
Therefore, since $p<p_{_N}$,
by classical elliptic estimates and a bootstrap argument
we have $\vn\in C^{2}(\om)$.}
\erm

\bigskip

{\bf Example.} Consider the unit ball~$\Omega=B_1\subset \R^N$.
With~$a$ and~$R_a$ as in \rife{entire.sol} and~$\eps_n\coloneqq \frac{1}{\sqrt{\mu_n}}\to 0^+$, the radial function~$v_n\colon B_1\to \R^+$ defined by
\beq\label{spike.sol}
\graf{v_n(x)=1+R_a^{\frac{2}{1-p}}u(\frac{|x|}{R_a\eps_n}),\;\; |x|\in [0,R_a\eps_n],\\ \\
v_n(x)=a+\left(\frac{R_a\eps_n}{|x|}\right)^{N-2}(1-a),\;\; |x|\in [R_a\eps_n,1],}
\eeq
is a solution of \rife{vn.1} in $\om=B_1$ with spike set $\Sigma=\{0\}$,
\beq\label{example1}
\lim\limits_{n\to +\ii}\mu_n^{\frac{N}{2}}\int\limits_{B_1}[\vn-1]_+^p=M_{p,a},
\eeq
and in particular, for $t\geq 1$,
\beq\label{example2}
\limsup\limits_{n\to +\ii}\mu_n^{\frac{N}{2}}\int\limits_{B_1}\vn^t\leq C
\mbox{ if and only if } a=0 \mbox{ and } t> \frac{N}{N-2}.
\eeq

\bigskip

These examples motivates Theorem \ref{thm5.7} below and its refinement Theorem \ref{thm5.7.intro} already
stated in the introduction. Recall that the spikes set for solutions of
\rife{vn.1} was introduced in Definition \ref{spikedef}.

\bte\label{thm5.7}
Let $\vn$ be a sequence of solutions of \rife{vn.1}, such that
\beq\label{mass}
\mu_n^{\frac N 2}\ino [\vn-1]_+^p\leq C_0.
\eeq
Then, $\vn$ is bounded in $L^{\infty}_{\rm loc}(\om)$ and there exists a subsequence, still denoted $\vn$, such that
\begin{itemize}
 \item[{\em either}] {\rm (i)} $[v_n-1]_+\to 0$ locally uniformly in $\om$, $\vn\to v$  in $W^{1,N}_{\rm loc}(\om)\cap L^{t}_{\rm loc}(\om)$ for any $t\geq 1$, for some harmonic function $v\leq 1$ in $\om$ and
$\mu_n [v_n-1]_+^p\to 0$ in $L^{\frac{N}{2}}_{\rm loc}(\om)$;

 \item[{\em or}~$\quad$] {\rm (ii)} the spikes set
relative to $\vn$ is not empty, i.e.
there exists a nonempty set $\Sigma_0\subset \om$ (the interior spikes set)
such that for any open and relatively compact set $\om_0\Subset \om$
there exists a subsequence such that
for any $z\in\Sigma_0\cap \om_0$ there exists $z_n\to z$, $\sg>0$ and $a=a_z\in [0,1)$ such that
$\vn(z_n)\geq 1+\sg$ and
$w_n(x):=\vn(z_n+\eps_n x)\to w_a(x)$ in $C^{2}_{\rm loc}(\R^N)$, where $\eps_n=\mu_n^{-\frac12}$, $w_a$ is the unique solution
of \rife{entire.a} and for any $R\geq 1$ we
have
$$
w_n(0)=\max\limits_{B_R}w_n,\quad \vn(z_n)=\max\limits_{B_{R\eps_n}(z_n)}\vn,
$$
where $z_n$ is the unique maximum point of $\vn$ in $B_{R\eps_n}(z_n)$. Moreover for any $z\in \Sigma$, we have
$$
\lim\limits_{R\to +\ii}\lim\limits_{n \to +\ii}\mu_n^{\frac N 2}
\int\limits_{B_{R\eps_n}(z_n)} [\vn-1]_+^p= M_{p,a_z},
$$
and for any~$r>0$,
$$
\lim\limits_{n \to +\ii}\|\mu_n[\vn-1]_+^p\|_{L^{\frac{N}{2}}(B_{r}(z))}> 0.
$$
\end{itemize}

\ete

\proof \, {\bf STEP 1:} {\it $\vn$ is bounded in $L^{\infty}_{\rm loc}(\om)$.}\\
By contradiction there exist open and relatively compact subset $\om_0\Subset \om$,
and $\om_0\ni x_n\to x_0\in \ov{\om}_0$ such that
\beq\label{blowup}
\vn(x_n)\to +\ii.
\eeq

We define the scaling parameter $\eps_n= \frac{1}{\sqrt{\mu_n}}$,
to obtain functions~$\widetilde{w}_n\colon \Omega_n\to \R^+$ where
\begin{align}\nonumber
\widetilde{w}_n(x)=\vn(x_n + \eps_n x),\qquad x\in \om_n:=\{x\in\R^N\,:\,x_n +\eps_n x\in \om\}.
\end{align}
Note that for any $R\geq 1$, $B_R\Subset \om_n$ for $n$ large enough.
Remark that
$$
\int\limits_{\om_n} {\dsp [\widetilde{w}_n-1]_+^p}=\eps_n^{-N} \int\limits_{\om} {\dsp [\vn-1]_+^p}\equiv
\mu_n^{\frac N 2} \int\limits_{\om} {\dsp [\vn-1]_+^p}\leq {C_0},
$$
whence $\widetilde{w}_n$ satisfies,
$$
\graf{-\Delta \widetilde{w}_n =[\widetilde{w}_n-1]_+^p\quad \mbox{in}\;\;\om_n,\\ \\
\bigints\limits_{\om_n} {\dsp [\widetilde{w}_n-1]_+^p}\leq{C_0},\\ \\
\widetilde{w}_n\geq 0 \mbox{ in } \om_n.
}
$$
Therefore we deduce from Theorem \ref{thm5.2} that, possibly passing to a subsequence,
$\widetilde{w}_n\to w$ in $C^{2}_{\rm loc}(\R^N)$ where $w\in C^{2}(\R^N)$ is a solution of

$$
\graf{-\Delta w =[w-1]_+^p\quad \mbox{in}\;\;\R^N,\\ \\
\bigints\limits_{\R^N} {\dsp [w-1]_+^p}\leq{C_0},\\ \\
w\geq 0 \mbox{ in } \R^N.
}
$$
As a consequence we have
$\vn(x_n)=\widetilde{w}_n(0)\to w(0)$, which is a contradiction to \rife{blowup}.

\

{\bf STEP 2:} {\it Assume that $[\vn -1]_+\to 0$ locally uniformly in $\om$. Then
$\vn\to v$  in $W^{1,N}_{\rm loc}(\om)\cap L^{t}_{\rm loc}(\om)$ for any $t\geq 1$,
for some harmonic function $v\leq 1$ in $\om$ and
$\mu_n [v_n-1]_+^p\to 0$ in $L^{\frac{N}{2}}_{\rm loc}(\om)$.}

\

Let $B_{4R}\Subset \om$ be any relatively compact ball and for each~$n$ let us write $\vn=v_{1,n}+v_{2,n}$ where
$$
\graf{-\Delta v_{1,n} =\mu_n [v_n-1]_+^p\quad \mbox{in}\;\;B_{4R}\\ \\
v_{1,n}=0 \mbox{ on }\pa B_{4R}}
$$
while $\Delta v_{2,n}=0$ in $B_{4R}$.
%{\color{red}The following lines can be simplified, since now~$v_n, v_{1,n}, v_{2,n}$ are all positive functions. Hence~$v_{i,n}\le v_n$ for~$i=1,2$, and we already know~$(v_n)$ is locally bounded. Thus we can directly conclude the convergence of~$v_{i,n}$ for~$i=1,2$, in the desired spaces.}
By the maximum principle,~$v_{1,n}$ and~$v_{2,n}$ are both positive, and
\begin{align}\nonumber
 0\leq v_{1,n}< v_n, \quad  0< v_{2,n}\leq v_n \quad \mbox{ in } B_{4R}.
\end{align}
Being a sequence of bounded positive harmonic functions, $v_{2,n}$ sub-converges to a
nonnegative harmonic function $v_2$, uniformly on any compact subset of $B_{4R}$, with
$v_2\leq 1$ because the assumption $[\vn -1]_+\to 0$ locally uniformly in $\om$ implies
$$
\limsup\limits_{n\to +\ii}\sup\limits_{B_R}v_{2,n}\leq
\limsup\limits_{n\to +\ii}\sup\limits_{\pa B_R}v_{n}\leq 1.
$$
%
%
% First of all by the mean value formula we have
% $$
% \|v_{2,n}\|_{L^{\ii}(B_{2R})}\leq C_R\int\limits_{B_{4R}}|v_{2,n}|\leq
% C_R\int\limits_{B_{4R}}v_{1,n}+C_R\int\limits_{B_{4R}}\vn,
% $$
% where we use that $v_{1,n}\geq 0$. Recall that by assumption
%  $[\vn -1]_+\to 0$ uniformly in $B_{4R}$, whence we can write
% $$
% \int\limits_{B_{4R}}\vn
% \leq \int\limits_{B_{4R}}\vn\chi(\{0\leq \vn\leq1\})+
% \int\limits_{B_{4R}}\vn\chi(\{\vn>1\}) \leq
% $$
% $$
% C_{1,R}+\int\limits_{B_{4R}}[\vn-1]_+\leq C_{2,R}.
% $$
% Also $\int\limits_{B_{4R}}v_{1,n}\leq C$ by classical results \cite{stam}.
% Therefore in particular $v_{2,n}$ is a sequence of harmonic functions bounded in $L^{\ii}(B_{2R})$
% and consequently possibly along a subsequence converges uniformly to an harmonic function $v_{2}$ in $B_R$.
% Observe also that by the maximum principle
% $$
% \limsup\limits_{n\to +\ii}\sup\limits_{B_R}v_{2,n}\leq
% \limsup\limits_{n\to +\ii}\sup\limits_{\pa B_R}v_{n}\leq 1,
% $$
% because of the assumption $[\vn -1]_+\to 0$ locally uniformly in $\om$. Therefore $v_2\leq 1$ in $B_R$.\\
Concerning $v_{1,n}$ we have,
\beq\label{holdrq0}
\int\limits_{B_{4R}} (\mu_n [v_n-1]_+^p)^{\frac N 2}
=\mu_n^{\frac N 2}\int\limits_{B_{4R}} [v_n-1]_+^p [v_n-1]_+^{p\frac{N-2}{2}}
\leq C_0\|[v_n-1]_+\|^{p\frac{N-2}{2}}_{L^{\ii}(B_{4R})},
\eeq
implying that $\mu_n [v_n-1]_+^p\to 0$ in $L^{\frac N 2}(B_{4R})$ and then by standard
elliptic regularity theory and the Sobolev embedding that $v_{1,n}\to v_1=0$ in
$W^{1,N}_0(B_{4R})\cap L^t(B_{4R})$ for any $t\in [1,+\ii)$.

\

In particular we deduce
that
$\vn=v_{1,n}+v_{2,n}\to v_2\leq 1$ in $W^{1,N}_{\rm loc}(\om)\cap L^t_{\rm loc}(\om)$ for any $t\in [1,+\ii)$.

\bigskip

{\bf STEP 3:} In view of STEP 2 we can assume w.l.o.g. that
there exists an open and relatively compact subset $\om_0\Subset \om$ and a sequence $\{x_n\}\subset \om_0$,
such that, possibly along a subsequence, $x_n\to z_0\in {\om}_0$ and $\vn(x_n)\geq 1+\sigma$ for any $n\in \N$
for some $\sigma>0$.\\
As in STEP 1, the rescaled functions
$\widetilde{w}_n(x)\coloneqq \vn(x_n+\eps_n x)$, with~$\eps_n=\frac{1}{\sqrt{\mu_n}}$,
defined on $\om_n=\{x\in\R^N\,:\,x_n+\eps_n x\in \om\}\nearrow \R^N$, satisfy
\beq\label{tildwn}
\graf{-\Delta \widetilde{w}_n =[\widetilde{w}_n-1]_+^p\quad \mbox{in}\;\;\om_n,\\ \\
\bigints\limits_{\om_n} {\dsp [\widetilde{w}_n-1]_+^p}=
\mu_n^{\frac N 2} \int\limits_{\om} {\dsp [\vn-1]_+^p}\leq C_0,\\ \\
\widetilde{w}_n\geq 0 \mbox{ in } \om_n,\\ \\
\widetilde{w}_n(0)=\vn(0)\geq 1+\sg>1.
}
\eeq

Therefore we deduce from Theorem \ref{thm5.2} that, possibly passing to a subsequence,
$\widetilde{w}_n\to w$ in $C^{2}_{\rm loc}(\R^2)$ where $w$ is a solution of

$$
\graf{-\Delta w =[w-1]_+^p\quad \mbox{in}\;\;\R^N,\\ \\
\bigints\limits_{\R^N} {\dsp [w-1]_+^p}\leq C_0,\\ \\
w\geq 0 \mbox{ in } \R^N\\ \\
w(0)\geq 1+\sg>1.
}
$$
As a consequence of Lemma \ref{lem5.3} we see that there exists $y_0$ and $a\in [0,1)$
such that $w(x)=w_a(x-y_0)$ where $w_a$ is the unique solution of \rife{entire.a}.
By Remark \ref{rmdgp} we see that $x=0$ is a nondegenerate maximum point of $w_a(x)$.
Therefore, since $\widetilde{w}_n\to w$ in $C^{2}_{\rm loc}(\R^N)$,
there exists $y_n\to y_0$
such that, for any $R\geq 1$ and any $n$ large enough, we have
$$
\widetilde{w}_n(y_n)=\vn(x_n+\eps_n y_n)=\max\limits_{B_{R\eps_n}(x_n+\eps_n y_0)}\vn=
\max\limits_{B_R(y_0)}\widetilde{w}_n.
$$

At this point it is easy to see that $z_n=x_n+\eps_ny_n\to z_0$ and the
functions rescaled with respect to the centers~$z_n$ (replacing the centers $x_n$)
$$
w_n(x)\coloneqq \vn(z_n+\eps_n x),
$$
indeed satisfy $w_n\to w_a$ in $C^{2}_{\rm loc}(\R^N)$ and in particular for any $R\geq 1$ it holds,
$$
w_n(0)=\max\limits_{B_R}w_n,\quad \vn(z_n)=\max\limits_{B_{R\eps_n}(z_n)}\vn,
$$
for any $n$ large enough, where $z_n$ is the unique (still by Remark \ref{rmdgp}) maximum point of $\vn$ in $B_{R\eps_n}(z_n)$.
Since by assumption we have
$$
\mu_n^{\frac N 2}\int\limits_{B_{r}(z_0)} {\dsp [v_n-1]_+^p}\leq C_{0},
$$
then $z_0$ is a spike point, as claimed.
In particular we have,
$$
\lim\limits_{R\to +\ii}\lim\limits_{n +\ii}\mu_n^{\frac N 2}\int\limits_{B_{R\eps_n}(z_n)} [\vn-1]_+^p= M_{p,a_z}.
$$

At last remark that, for any fixed $z\in \Sigma$, denoting by $z_n\to z$ the sequence of local maximum points,
$a=a_z$, $R_a=R_{a_z}$ and $r>0$ any small enough radius, we have,
\begin{align}\nonumber
 \lim\limits_{n \to +\ii}
\mu_n^{\frac N 2}\int\limits_{B_r(z)} [\vn-1]_+^{p\frac{N}{2}}
\geq &
\lim\limits_{R\to +\ii}\lim\limits_{n \to +\ii}
\mu_n^{\frac N 2}\int\limits_{B_{R\eps_n}(z_n)} [\vn-1]_+^{p\frac{N}{2}} \\
=& \lim\limits_{n \to +\ii}\int\limits_{B_{2R_a}(0)} [\widetilde{w}_n-1]_+^{p\frac{N}{2}} \nonumber \\
=&
\lim\limits_{n \to +\ii}\int\limits_{B_{2R_a}(0)}
[\widetilde{w}_n-1]_+^{p}[\widetilde{w}_n-1]_+^{p\frac{N-2}{2}}
\nonumber \\
=& \int\limits_{B_{R_a}}
[w_a-1]_+^{p}[w_a-1]_+^{p\frac{N-2}{2}}>0, \nonumber
\end{align}
as claimed.
\finedim

\brm\label{remnobdy}{\it
As far as we miss further assumptions, as for example  in \cite{FW} or either in Theorem \ref{thm5.7.intro},
it seems not easy to say much more about the spikes set, mainly due to the fact that $M_{p,a}\searrow 0^+$ as
$a\nearrow 1^-$. In principle this could allow the existence of a spikes set of local unbounded cardinality.
More in general, in view of \rife{example2}, for $a\neq 0$ the $\vn$ as defined in \rife{spike.sol} provides an example
of a sequence which is not, according to Definition \ref{spikeseqdef}, a regular spikes sequence. A natural
assumption to possibly rule out non-regular spikes turns out to be \eqref{mass.intro.2} in Theorem \ref{thm5.7.intro}, which
allows one to deduce the so called Minimal Mass Lemma, see Lemma \ref{minmass} below.}
\erm

As above we set $\eps_n=\mu_n^{-\frac12}$ while $R_0$, $M_{p,0}$ are defined in \rife{entire.Ra}
and \rife{entire.mp0}.
Let us also remark that the case (a) in Theorem \ref{thm5.7.intro}
happens for example for constant solutions and for solutions as \eqref{spike.sol} when the spike
point seats on the boundary. A preliminary crucial Lemma will be needed during the proof. Recall
that by definition a (sub)domain is an open and connected set.
\ble[The Vanishing Lemma]\label{vanishle}
Let $\vn$ satisfy all the assumptions of Theorem \ref{thm5.7.intro} and
assume that $[\vn-1]_+\to 0$ locally uniformly in a subdomain $\om'\subseteq \om$.
Then for any open and relatively compact subset $\om_0\Subset \om'$ there exists ${n}_0\in \N$ and
$C>0$, depending on $\om_0$,
such that $[v_n-1]_+=0$ in ${\om}_0$ and in particular,
\beq\label{vanish:decay}
\|\vn\|_{L^{\ii}(\om_0)}\leq C\eps_n^{\frac N t},
\eeq
for any $n>{n}_0$.
\ele
\proof
We argue by contradiction and assume that
there exists an open and relatively compact subset $\om_0\Subset \om'$ and a sequence $\{x_n\}\subset \om_0$,
such that, possibly along a subsequence, $x_n\to z_0\in {\om}_0$, $v_n(x_n)\to 1$ and $\vn(x_n)>1$
for any $n\in \N$. At this point we follow step by step the construction of STEP 3 in the
proof of Theorem \ref{thm5.7}
and define $\widetilde{w}_n$ and $\om_n$ in the same way. Therefore $\widetilde{w}_n$ satisfies
\rife{tildwn} and we observe that,
due to the assumption \rife{mass.intro.2},
in this case we also have,
\beq\label{goodbd1}
\int\limits_{\om_n}|\widetilde{w}_n|^{t}\leq C_1.
\eeq
As a consequence, we deduce from Theorem \ref{thm5.2} that, possibly passing to a subsequence,
$\widetilde{w}_n\to w$ in $C^{2}_{\rm loc}(\R^2)$ where $w\geq 0$ in $\R^N$ is a solution of

\beq\label{eqinfty}
\graf{-\Delta w =[w-1]_+^p\quad \mbox{in}\;\;\R^N,\\ \\
\bigints\limits_{\R^N} {\dsp [w-1]_+^p}\leq{C_0},\\
\bigints\limits_{\R^N} |{w}|^{t}\leq C_1.
}
\eeq

The condition $\int\limits_{\R^N}|{w}|^{t}\leq C_1$ rules out positive constants and solutions
which do not vanish at infinity.
Therefore, by Lemma \ref{lem5.3} and Remark \ref{rmd12}, we see that either $w\equiv 0$ or there exists $y_0$
such that $w(x)=w_0(x-y_0)$ where $w_0$ is the unique solution of
\rife{entire.a} with $a=0$.
Actually $w$ cannot vanish identically, otherwise
$\widetilde{w}_n\to 0$ in $C^{2}_{\rm loc}(\R^2)$ would imply $\vn(x_n)<1$ for $n$ large.
Since by Remark \ref{rmdgp} we have that  $x=0$ is a nondegenerate maximum point of $w_0$, then $\widetilde{w}_n\to w$ in $C^{2}_{\rm loc}(\R^N)$
implies that  there exists $y_n\to y_0$
such that, for any $R\geq 1$ and any $n$ large enough, we have
$$
\widetilde{w}_n(y_n)=
\max\limits_{B_R(y_0)}\widetilde{w}_n.
$$
As a consequence, by setting $z_n=x_n+\eps_n y_n\to z_0$, then
$$
w_n(x)=\vn(z_n+\eps_n x),
$$
satisfies $w_n\to w_0$ in $C^{2}_{\rm loc}(\R^N)$ and in particular for any $R\geq 1$ it holds,
$$
\max\limits_{B_{2R}}w_n=w_n(0)=\vn(z_n)=\max\limits_{B_{2R\eps_n}(z_n)}\vn,
$$
for any $n$ large enough, where $z_n$ is the unique maximum point of $\vn$ in $B_{R\eps_n}(z_n)$.\\
In particular since by construction $w_n(0)\to w_0(0)>1$, and since by assumption $[\vn-1]_+\to 0$ uniformly near
$z_0$, we would also have,
$$
1<\lim\limits_{n\to +\ii}w_n(0)=\lim\limits_{n\to +\ii}\vn(z_n)\leq 1,
$$
which is the desired contradiction. Therefore $\vn$ is harmonic for $n$ large in $\om_0$ and by the mean value theorem
and \rife{mass.intro.2} we deduce that for any ball $B_{2R}(x_0)\Subset \om_0$ we have,
$$
\|\vn\|_{L^{\ii}(B_R(x_0))}\leq C_R\eps_n^{\frac N t}.
$$
The proof of the Vanishing Lemma is completed by taking a finite cover of the relatively compact set
$\Omega_0\Subset \Omega^{'}$. $\square$\\

\brm\label{vanish-equiv}{\it  There is a useful, almost equivalent formulation of the Vanishing Lemma,
which is just what the most part of its proof shows which we call the Non-Vanishing Lemma.\\

{\bf Lemma}{\rm (The Non-Vanishing Lemma)}. If there exists a sequence
$\{x_n\}\Subset \om$ such that $\vn(x_n)>1$, then
there exists $z_0\in \om$ and $\{z_n\}\Subset \om$ such that,
possibly along a subsequence still denoted $\vn$, we have $z_n\to z_0$, $|z_n-x_n|\leq C\eps_n$
and $w_n(x)=\vn(z_n+\eps_nx)\to w_0(x)$
in $C^{2}_{\rm loc}(\R^N)$. In particular $\vn(z_n)\to w_0(0)>1$, i.e.
$\vn(z_n)\geq 1+ \sg$ for some $\sg>0$.$\square$}
\erm

\bigskip

Next we prove Theorem \ref{thm5.7.intro}.

{\bf The proof of Theorem \ref{thm5.7.intro}.}\\
The claim in $(a)$ follows just by taking $\om'=\om$ in the Vanishing Lemma,
whence we are just left with the proof of $(b)$.

{\bf STEP 1}: Assume that $[\vn-1]_+$ does not converge to $0$ locally uniformly, then
by either Theorem \ref{thm5.7} or the Non-Vanishing Lemma (Remark \ref{vanish-equiv}),
we necessarily have that the interior spikes set $\Sigma_0=\Sigma\cap\om$
relative to a suitable subsequence is not empty. Let us choose any point $z_1\in \Sigma_0$ so that, by definition,
there exists $x_{n}\to z_1$ such that
$\vn(x_{n})\geq 1+ \sg$, for some $\sg>0$. In view of Remark \ref{vanish-equiv} we just follow
the argument of the Vanishing Lemma and define $\widetilde{w}_{n}$, $y_n$ and ${w}_{n}$
in the same way. Therefore
$z_n=x_n+\eps_n y_n \to z_1$ and
$w_{n}(x)=\vn(z_n+\eps_n x)$ satisfies, possibly passing to a subsequence, $w_n\to w_0$
in $C^{2}_{\rm loc}(\R^N)$,
where $w_0$ is the unique solution of \eqref{entire.a} and in particular for any $R\geq 1$ it holds,
\beq\label{max.1}
\max\limits_{B_{2R}(0)}w_n=w_n(0)=\vn(z_n)=\max\limits_{B_{2R\eps_n}(z_n)}\vn,
\eeq
for any $n$ large enough, where $z_n$ is the unique and non-degenerate (see Remark \ref{rmdgp})
maximum point of $\vn$ in $B_{R\eps_n}(z_n)$.
Obviously no contradiction arise in this case and in particular
$$
\lim\limits_{n\to +\ii} \int \limits_{B_{2R_{0}}(0)}[w_n-1]_+^p=\lim\limits_{n\to +\ii} \int \limits_{B_{R_{0}}(0)}[w_n-1]_+^p=M_{p,0},
$$
which is, scaling back, the same as,
\beq\label{mass-limit}
\lim\limits_{n\to +\ii} \mu_n^{\frac N 2}\int \limits_{B_{2\eps_n R_{0}}(z_n)}[\vn-1]_+^p=\lim\limits_{n\to +\ii} \mu_n^{\frac N 2}\int \limits_{B_{\eps_n R_{0}}(z_n)}[\vn-1]_+^p=M_{p,0}.
\eeq
By a diagonal argument, along a subsequence which we will not relabel,
we can find $R_{1,n}\to +\ii$ such that,
$$
\eps_n R_{1,n}\to 0, \quad \|w_{n}-w_0\|_{C^2(B_{2 R_{1,n}}(0))}\to 0,
$$
$$
\left(\frac{R_0}{2|x|}\right)^{N-2}\leq w_{n}(x)\leq \left(\frac{2R_0}{|x|}\right)^{N-2},\quad
2 R_0 \leq|x|\leq 2 R_{1,n},
$$
\beq\label{decay3}
\left(\frac{R_0\eps_n}{2|x-z_n|}\right)^{N-2}\leq \vn(x)\leq \left(\frac{2R_0\eps_n}{|x-z_n|}\right)^{N-2},\quad
2 \eps_n R_0 \leq|x-z_n|\leq 2 \eps_n R_{1,n},
\eeq
and for any $r<R_0<R$,
\beq\label{profile1}
B_{\eps_n r}(z_n)\Subset \{x\in B_{2\eps_n R_{1,n}}(z_n)\,:\vn(x)>1\,\}\Subset B_{\eps_n R}(z_n),
\eeq
for any $n$ large enough which, in view of \rife{mass-limit}, implies that
\beq\label{mass1}
\lim\limits_{n\to +\ii} \mu_n^{\frac N 2}\int \limits_{B_{2\eps_n R_{1,n}}(z_n)}[\vn-1]_+^p=M_{p,0}.
\eeq

There is of course no loss of generality in assuming that,
\beq\label{eprn}
\eps_n^{\frac{N}{t}}\leq \frac{1}{R_{1,n}^{N-2}}.
\eeq

In particular, what the argument right above shows is the following,
\ble[Minimal Mass Lemma]\label{minmass} Let $z_1\in \Sigma_0$ be any spike point so that, by definition,
there exists $x_{n}\to z_1$ such that
$\vn(x_{n})\geq 1+ \sg$, for some $\sg>0$. Then there exists $z_{n}\to z_1$ such that
$$
\lim\limits_{R\to +\ii}\lim\limits_{n +\ii}\mu_n^{\frac N 2}\int\limits_{B_{R\eps_n}(z_n)} [\vn-1]_+^p= M_{p,0}.
$$
\ele

{\bf STEP 2}: By STEP 1 $\Sigma_0 \neq \emptyset$ while by the Minimal Mass Lemma and \eqref{mass.intro.1} we have
$$
\#(\Sigma_0)\leq \frac{C_0}{M_{p,0}}.
$$

Therefore $\Sigma_0$ is finite and by definition we have $[\vn-1]_+\to 0$
locally uniformly in $\om\setminus \Sigma_0$. As a consequence
of (a) we see that for any open and relatively compact set
${\om}_0\Subset \om\setminus \Sigma_0$
there exists ${n}_0\in \N$ such that $[v_n-1]_+=0$ in ${\om}_0$ and in particular
\beq\label{decay1}
\|\vn\|_{L^{\ii}(\om_0)}\leq C\eps_n^{\frac N t},
\eeq
for any $n>{n}_0$.\\
At this point, for $z_n\equiv z_{1,n}$ and $R_{1,n}$ defined as in STEP 1 and for any domain $\om_m\Subset \om$ such that $\Sigma_0\subset \om_m$, we define
$$
z_{2,n}:\vn(z_{2,n})=\max\limits_{\ov{\om}_m \setminus B_{\eps_n R_{1,n}}(z_{1,n})\}}\vn.
$$
There are only two possibilities:
\begin{itemize}
 \item[either] (j) $[\vn(z_{2,n})-1]_+\to 0$,
 \item[or~$\quad$] (jj) $\vn(z_{2,n})\geq 1+ \sg$, for some $\sg>0$.
\end{itemize}

We discuss the two cases separately. If (j) holds we define $m=1$, $z_{1,n}=z_n$ and $R_n=R_{1,n}$ and claim that
\beq\label{noresmass0}
\sup\limits_{\om_m\setminus B_{\eps_n R_{n}}(z_{1,n})}[\vn-1]_+=0,
\eeq
for any $n$ large enough and
\beq\label{noresmass1}
\lim\limits_{n\to +\ii} \mu_n^{\frac N 2}\int\limits_{\om_m}[\vn-1]^p_+\to M_{p,0}.
\eeq

In view of \rife{mass1}, clearly \rife{noresmass0} implies \rife{noresmass1} which is why we just prove
the former. By the Vanishing Lemma or either Remark \ref{vanish-equiv} it is enough to rule out $z_{2,n}\to z_1$.
By contradiction, if \rife{noresmass0} where false we would have
$\vn(z_{2,n})>1$ for any $n$ and we could define ${w}_{2,n}(x)=\vn(z_{2,n}+\eps_n x)$ for
$|x|\leq \frac{R_{n}}{2}$. Remark that
$$
\left\{z_{2,n}+\eps_n x,|x|< \frac{R_{n}}{2}\right\}\bigcap
\left\{z_{1,n}+\eps_n x,|x|< \frac{R_{n}}{2}\right\}=\emptyset,
$$
whence we deduce as in the proof of the Vanishing Lemma that $w_{2,n}\to w_0$ in $C^{2}_{\rm loc}(\R^N)$
and in particular that for any $R\geq 1$ we have,
\beq\label{max2}
\max\limits_{B_R}w_{2,n}=w_{2,n}(0)=\vn(z_{2,n})=\max\limits_{B_{R\eps_n}(z_{2,n})}\vn,
\eeq
for any $n$ large enough, where $z_{2,n}$ is the unique maximum point of $\vn$ in $B_{R\eps_n}(z_{2,n})$.\\
In particular since by construction $w_{2,n}(0)\to w_0(0)>1$, and since by assumption
$[\vn(z_{2,n})-1]_+\to 0$, we would also have,
$$
1<\lim\limits_{n\to +\ii}w_{2,n}(0)=\lim\limits_{n\to +\ii}\vn(z_{2,n})\leq 1,
$$
which is the desired contradiction. Therefore both \eqref{noresmass0} and \eqref{noresmass1}
holds true. Actually,
we deduce by \rife{noresmass0} that $\vn$ is harmonic in $\om_m \setminus B_{R_n\eps_n}(z_{1,n})$
and since by
\rife{decay3},\rife{decay1},  we have
$$
0\leq \sup\limits_{\pa \om_m \cup \pa B_{2 R_n\eps_n}(z_{1,n})}\vn\leq
C_*\max\left\{\eps_n^{\frac N t},\frac{1}{R_n^{N-2}}\right\}
$$
for any $n$ large enough, we also deduce from \eqref{eprn} that,
$$
0\leq \sup\limits_{\om_m \setminus B_{2 R_n\eps_n}(z_{1,n})}\vn\leq
\frac{C_*}{R_n^{N-2}},
$$
for any such $n$. In particular, in view of \rife{max.1}, possibly choosing a larger $n_1$, we have
$$
\vn(z_{1,n})=\max\limits_{\om_m}\vn,\quad \forall\,n>n_1,
$$
and by \rife{profile1}, for any $r<R_0<R$,
\beq\label{profile1.1}
B_{\eps_n r}(z_{1,n})\Subset \{x\in \om_m\,:\vn(x)>1\,\}\Subset B_{\eps_n R}(z_{1,n}),\quad \forall\,n>n_1.
\eeq

This completes the proof as far as (j) holds just putting $n_*=n_1$ and then we are left with (jj).
Assuming that (jj) holds, we define ${w}_{2,n}(x)=\vn(z_{2,n}+\eps_n x)$ and
again as in the proof of the Vanishing Lemma we deduce that $z_{2,n}\to z$ for some $z\in \Sigma_0$, $w_{2,n}\to w_0$ in $C^{2}_{\rm loc}(\R^N)$
and in particular that, for any $R\geq 1$, \rife{max2} holds. Remark that we do not exclude that $z=z_1$. Moreover we have,
$$
\lim\limits_{n\to +\ii} \int \limits_{B_{2R_{0}}(0)}[w_{2,n}-1]_+^p=
\lim\limits_{n\to +\ii} \int \limits_{B_{R_{0}}(0)}[w_{2,n}-1]_+^p=M_{p,0},
$$
which is, scaling back, the same as,
\beq\label{mass9}
\lim\limits_{n\to +\ii} \mu_n^{\frac N 2}\int \limits_{B_{2\eps_n R_{0}}(z_{2,n})}[\vn-1]_+^p=
\lim\limits_{n\to +\ii} \mu_n^{\frac N 2}\int \limits_{B_{\eps_n R_{0}}(z_{2,n})}[\vn-1]_+^p=M_{p,0}.
\eeq
Actually, again by a diagonal argument we can find $\{R_{2,n}\}\subset \{R_{1,n}\}$, $R_{2,n}\to +\ii$,
such that, along a subsequence which we will not relabel, we have,
$$
\eps_n R_{2,n}\to 0, \quad \|w_{2,n}-w_0\|_{C^2(B_{2 R_{2,n}}(0))}\to 0,
$$
$$
\left(\frac{R_0}{2|x|}\right)^{N-2}\leq w_{i,n}(x)\leq \left(\frac{2R_0}{|x|}\right)^{N-2},\quad
2 R_0 \leq|x|\leq 2 R_{2,n},\;i=1,2,
$$
$$
\left(\frac{R_0\eps_n}{2|x-z_{i,n}|}\right)^{N-2}\leq \vn(x)\leq \left(\frac{2R_0\eps_n}{|x-z_{i,n}|}\right)^{N-2},\quad
2 \eps_n R_0 \leq|x-z_{i,n}|\leq 2 \eps_n R_{2,n},\;i=1,2,
$$
$$
B_{2 \eps_n R_{2,n}}(z_{1,n})\cap B_{2 \eps_n R_{2,n}}(z_{2,n})=\emptyset,
$$
and for any $r<R_0<R$,
$$
\underset{i=1,2}{\cup}B_{\eps_n r}(z_{i,n})\Subset\{x\in \underset{i=1,2}{\cup}B_{2\eps_n R_n}(z_{i,n})\,:\,\vn>1\}\Subset
\underset{i=1,2}{\cup}B_{\eps_n R}(z_{i,n}),
$$
for any $n$ large enough which in view of \rife{mass9} implies that,
\beq\label{mass2}
\lim\limits_{n\to +\ii} \mu_n^{\frac N 2}\int \limits_{B_{2\eps_n R_{2,n}}(z_{2,n})}[\vn-1]_+^p=M_{p,0}.
\eeq

At this point, setting
$$
z_{3,n}:\vn(z_{3,n})=\max\limits_{\ov{\om}_m
\setminus B_{2 \eps_n R_{2,n}}(z_{n})\cup B_{2 \eps_n R_{2,n}}(z_{2,n})\}}\vn,
$$
and if $[\vn(z_{3,n})-1]_+\to 0$, we
argue as in (j) above to deduce that,
\beq\label{noresmass2}
\sup\limits_{\om_m\setminus \{\underset{i = 1,2}
{\cup}B_{2 \eps_n R_{n}}(z_{\!i,n})\}}[\vn-1]_+=0,
\eeq
for any $n$ large enough and, in view of \rife{mass1}, \rife{mass2},
\beq\label{noresmass3}
\lim\limits_{n\to +\ii} \mu_n^{\frac N 2}\int\limits_{\om_m}[\vn-1]^p_+= 2 M_{p,0}.
\eeq
In particular, putting $m=2$, we deduce as above that, possibly choosing a larger $n_*=n_1$,
$$
0\leq \sup\limits_{\om_m \setminus \{\underset{i = 1,2}
{\cup}B_{2 \eps_n R_{n}}(z_{i,n})\}}\vn\leq
C_*\max\left\{\eps_n^{\frac N t},\frac{1}{R_n^{N-2}}\right\}\leq \frac{C_*}{R_n^{N-2}},
\quad \forall\,n>n_*,
$$
and, for any $r<R_0<R$,
$$
\underset{i = 1,2}
{\cup}B_{2 \eps_n r}(z_{i,n})\Subset \{x\in \om_m\,:\vn(x)>1\,\}
\Subset \underset{i = 1,2}
{\cup}B_{2 \eps_n R}(z_{i,n}),\quad \forall\,n>n_*,
$$
which obviously concludes the proof in this case as well. Otherwise, by induction we can find $\sg>0$
such that there exists $\{R_n\}\subset \{R_{2,n}\}$, $R_{n}\to +\ii$ and
$N_m$ sequences $\{z_{i,n}\}\subset \om$, $i=1,\cdots,N_m$ such that, along a sub-subsequence
which we will not relabel we have,
$$
\vn(z_{i,n})=\max\limits_{\ov{\om}_m \setminus \{\underset{\ell = 1,\cdots,m,\,\ell\neq i}
{\cup}B_{2 \eps_n R_{n}}(z_{\!\ell,n})\}}\vn\geq 1+\sg,
$$
and defining $w_{i,n}(x)=\vn(z_{i,n}+\eps_n x)$, then
$$
\eps_n R_{n}\to 0,\quad \|w_{i,n}-w_0\|_{C^2(B_{2 R_{n}}(0))}\to 0,\; \forall i\geq 1,
$$
$$
\left(\frac{R_0}{2|x|}\right)^{N-2}\leq w_{i,n}(x)\leq \left(\frac{2R_0}{|x|}\right)^{N-2},\quad
2 R_0 \leq|x|\leq 2 R_{2,n},\;\forall i\geq 1,
$$
$$
\left(\frac{R_0\eps_n}{2|x-z_{i,n}|}\right)^{N-2}\leq \vn(x)\leq \left(\frac{2R_0\eps_n}{|x-z_{i,n}|}\right)^{N-2},\quad
2 \eps_n R_0 \leq|x-z_{i,n}|\leq 2 \eps_n R_{2,n},\;\forall i\geq 1,
$$
$$
B_{2 \eps_n R_{n}}(z_{\!\ell,n})\cap B_{2 \eps_n R_{n}}(z_{\!i,n})=\emptyset,\;\ell\neq i,
$$
and
$$
\lim\limits_{n\to +\ii} \int \limits_{B_{2R_{0}}(0)}[w_{i,n}-1]_+^p=
\lim\limits_{n\to +\ii} \int \limits_{B_{R_{0}}(0)}[w_{i,n}-1]_+^p=M_{p,0},\; \forall i\geq 1,
$$
that is, scaling back,
$$
\lim\limits_{n\to +\ii} \mu_n^{\frac N 2}\int \limits_{B_{2\eps_n R_{0}}(z_{i,n})}[\vn-1]_+^p=
\lim\limits_{n\to +\ii} \mu_n^{\frac N 2}\int \limits_{B_{\eps_n R_{0}}(z_{i,n})}[\vn-1]_+^p=M_{p,0},\;
\forall i\geq 1.
$$
Recall that since $N_m M_{p,0}\leq C_0$, then $N_m$ must be finite and then in particular
as in (j) above we deduce that, possibly choosing a larger $n_*\geq n_1$, and $C_*\geq C$,
$$
\sup\limits_{\om_m\setminus \{\underset{i = 1,\cdots,N_m}
{\cup}B_{2 \eps_n R_{n}}(z_{i,n})\}}[\vn-1]_+=0,
$$

$$
0\leq \sup\limits_{\om_m \setminus \{\underset{i = 1,\cdots,N_m}
{\cup}B_{2 \eps_n R_{n}}(z_{i,n})\}}\vn\leq
C_*\max\left\{\eps_n^{\frac N t},\frac{1}{R_n^{N-2}}\right\}\leq \frac{C_*}{R_n^{N-2}},\quad \forall\,n>n_*,
$$
and for any $r<R_0<R$,
$$
\underset{i=1,\cdots,N_m}{\cup}B_{\eps_n r}(z_{i,n})\Subset\{x\in \om_m\,:\,\vn>1\}\Subset
\underset{i=1,\cdots,N_m}{\cup}B_{\eps_n R}(z_{i,n}),\quad \forall\,n>n_*,
$$
and
$$
\lim\limits_{n\to +\ii} \mu_n^{\frac N 2}\int\limits_{\om_1}[\vn-1]^p_+= N_m M_{p,0}.
$$
Obviously we have,
\begin{align}\nonumber
 m\equiv \#\Sigma_0 \leq N_m \leq \frac{C_0}{M_{p,0}},
\end{align}
where $\Sigma_0=\{z_1,\cdots,z_m\}$. The proof is completed putting $\mathcal{Z}=\sum_{i=1}^mk_jz_j$, where $k_j$ is the multiplicity of $z_j$.
\finedim

\brm\label{rem:local} {\it
There is a version of Theorems \ref{thm5.7} and \ref{thm5.7.intro} for functions which are not
assumed to be nonnegative satisfying,
\beq\label{vn.1.neg}
\graf{-\Delta \vn =\mu_n [\vn-1]_+^p\quad \mbox{in}\;\;\om,\\ \\
\mu_n \to +\ii,\\ \\
\|[\vn]_-\|_1\leq C.}
\eeq

Indeed, as in the proof of Theorem \ref{thm5.2}, by the Kato inequality it is readily seen that
$[\vn]_-$ is bounded in $L^{\ii}_{\rm loc}(\om)$.  Also, we obviously take the same
definition of spikes set (see Definition \ref{spikedef}) and by a straightforward evaluation
obtain the following,\it }
\erm
\ble
Let $\vn$ be a sequence of solutions of \rife{vn.1.neg} and assume
$$
I_{n,p}:=\mu_n^{\frac N 2}\ino [\vn-1]_+^p<+\ii.
$$

For any fixed open set $\om_0\subseteq \om$ let
$-c:=\min\{0,\inf\limits_{n\in\N, x \in \ov{\om_0}}v_n(x)\}$ and~$\delta_c\coloneqq (1+c)^{-\frac{p-1}{2}}$.
Then the functions~$u_n\colon \Omega_{n,c}\to \R_+$ defined by
$$
 u_n(x)\coloneqq\frac{1}{1+c}(\vn(x_n+\delta_c x)+c),\quad \mbox{ for }
x\in \om_{n,c}:=\{x\in \R^N\,:\,x_n+\delta_c x\in \om_0\}
$$
satisfy
\beq\label{scaling1}
\graf{-\Delta u_n =\mu_n [u_n-1]_+^p\quad \mbox{in}\;\;\om_{n,c}, \\ \\
\mu_n \to +\ii,\\ \\
\mu_n^{\frac N 2}\int\limits_{\om_{n,c}} [u_n-1]_+^p= (1+c)^{-\tau_p}I_{n,p},\quad\mbox{ with }
\tau_p=\frac{N}{2}\left(1-\frac{p}{p_{_N}}\right)>0,\\ \\
\inf\limits_{\om_{n,c}}u_n\geq  0.
}
\eeq
\ele
{\it
Therefore the Definition \ref{spikeseqdef} of regular $\mathcal{Z}$-spikes sequence  and
Theorem \ref{thm5.7} can be applied to $u_n$. Also, by further assuming that
$\mu_n^{\frac N 2}\int\limits_{\om_{n,c}} u_n^t\leq C_t$ for some $t\geq 1$,
then Theorem \ref{thm5.7.intro} applies
to $u_n$ as well.} $\square$

\section{\bf Asymptotics of free boundary problems}\label{sec7}
Let $(\lmn,\pn)$ be a sequence of solutions of $\prl$ with $\an<0$.
We aim to study the case where $\alpha_n\to -\ii$ and first point out that in this case necessarily $\lambda_n\to +\ii$.

\bte\label{thm:7.1} Let $(\lmn,\pn)$ be a sequence of solutions of {\rm $\prl$} such that $\an<0$ and $|\an|\to +\ii$. Then $\lmn\to +\ii$
and in particular $\lmn\|\pn\|_\ii\to +\ii$.
\ete
\proof Concerning the first claim, by contradiction along a subsequence we would have $\an\to -\ii$ and $\lmn\leq \ov{\lm}$.
Recall by Theorem \ref{thm4.1} that then we would also have $\|\pn\|_{\ii}\leq C(\om,\ov{\lm})$ whence
$\sup\limits_{\om}(\an+\lmn\pn)\leq \an +\ov{\lm}C(\om,\ov{\lm})<0$ for $n$ large enough. This is impossible
in view of the integral constraint in $\prl$.
As for the second claim, again by contradiction
along a subsequence we
would have $\an\to -\ii$ and $\lmn\pn\leq C$, whence the contradiction arise in the same way.
\finedim

\bigskip

The asymptotic behavior is better understood via the functions~$\vn=\frac{\lmn}{|\an|}\pn$ which
satisfy \eqref{w5.1.intro}. As above we set~$\eps_n=\mu_n^{-1/2}\to 0^+$.

\medskip
\medskip

We first state an energy estimate based on ~\eqref{mass.intro.1} and~\eqref{eq:assump-2.intro}.
Remarks that~$\mathcal{C}_S$ below is the critical Sobolev constant,
\begin{align}
 \|u\|_{L^{2^*}(\Omega)} \leq \mathcal{C}_S\|\nabla u\|_{L^2(\Omega)},
\end{align}
which does not depend on $\om$.
\ble\label{lem:2star} Let $\vn$ satisfy all the assumptions of Theorem \ref{thm:spike-vn.intro}.
Then
\begin{align}\label{eq:Dirichlet bound-vn}
  \int_{\Omega} |\nabla v_n|^2 \dd x
  \leq  (C_0 + C_1)\eps_n^{N-2},
 \end{align}
 and
 \begin{align}\label{eq:L2*estimate}
  \mu_n^{\frac{N}{2}}\int_\Omega |v_n|^{2^*}\dd x \leq \mathcal{C}_S^{2^*}(C_0+C_1)^{\frac{N}{N-2}},
 \end{align}
 where $2^*=\frac{2N}{N-2}$ and $\mathcal{C}_S$ is the Sobolev constant.
\ele
\proof
 First of all, because of $\int\limits_{\om} {|\an|^p[\vn-1]_+^p}=1$, we have
\beq\label{unifb}
\mu_n^{\frac N 2}\int\limits_{\om} {[\vn-1]_+^p}=
\left(\frac{\lm_n}{|\an|^{1-\frac{p}{p_{_N}}}}\right)^{\frac N 2}
\leq C_0.
\eeq
Since~$v_n$ is positive in~$\Omega$, Theorem \ref{thm5.7} implies that $\vn$
is locally uniformly bounded. For each~$n\geq 1$, we denote the plasma region for~$v_n$
by~$\Omega_{+,n}\equiv\{v_n > 1\}$ and its vacuum region by~$\Omega_{-,n}\equiv \{v_n < 1\}$.
 In case both~$\Omega_{+,n}$ and~$\Omega_{-,n}$ are smooth enough to integrate by parts,
 we just test the equations in~$\Omega_{+,n}$ by~$(v_n-1)$ and in~$\Omega_{-,n}$ by~$v_n$ respectively,
 to deduce that,
 \begin{align}\nonumber
  \int_{\Omega_{+,n}} |\nabla v_n|^2 \dd x
  =\mu_n\int_{\Omega} [v_n-1]_+^{p+1} \leq C_1 \mu_{n}^{\frac{2-N}{2}}= C_1 \eps_n^{N-2},
 \end{align}
 and
 \begin{align}\nonumber
  \int_{\Omega_{-,n}} |\nabla v_n|^2\dd x
  =&-\oint_{\pa \Omega_{+,n} } \frac{\pa v_n}{\pa\nu}
  =\int_{\Omega_{+,n}} -\Delta v_n \dd x\\\nonumber
  =&\int_{\Omega} \mu_n [v_n-1]_+^p
  =\frac{\lambda_n}{|\alpha_n| } \int_{\Omega} |\alpha_n|^p [v_n-1]_+^p  \\\nonumber
  =& \frac{\lambda_n}{|\alpha_n|}.
 \end{align}
 Although we don't know in general whether or not such a regularity property about
 the free boundary holds true, since $\vn$ is smooth away from $\{\vn=1\}$,
 we can still use an approximation argument together with the Sard lemma
 as in the proof of Lemma~\ref{enrg0} to come up with the same formulas.
 Next, observe that
 \begin{align}\nonumber
  1= 1+ (p-1)(1-\frac{N}{2}) - (p-1)(1-\frac{N}{2})
   = \frac{N}{2}(1-\frac{p}{p_{_N}}) - (p-1)(1-\frac{N}{2}),
 \end{align}
 whence, recalling $\mu_n=\lambda_n |\alpha_n|^{p-1}$ and \eqref{unifb}, we find that,
 \begin{align}\nonumber
  \frac{\lambda_n}{|\alpha_n|}
  =& \frac{\lambda_n^{\frac{N}{2}}}{|\alpha_n|^{\frac{N}{2}(1-\frac{p}{p_{_N}})}}
  \frac{\lambda_n^{1-\frac{N}{2}}}{|\alpha_n|^{-(p-1)(1-\frac{N}{2})}}
  =\parenthesis{\frac{\lambda_n}{|\alpha_n|^{1-\frac{p}{p_{_N}}}}}^{\frac{N}{2}}
   \parenthesis{\lambda_n |\alpha_n|^{p-1}}^{\frac{2-N}{2}}
   \leq C_0 \eps_n^{N-2},
 \end{align}
which immediately implies that \eqref{eq:Dirichlet bound-vn} holds.
At last, by the Sobolev embedding, we have,
 \begin{align}\nonumber
  \int_{\Omega} v_n^{2^*} \leq \mathcal{C}_S^{2^*} \parenthesis{\int_{\Omega} |\nabla v_n|^2}^{\frac{N}{N-2}}
  \leq \mathcal{C}_S^{2^*}\parenthesis{C_0+C_1}^{\frac{N}{N-2}} \eps_n^N,
 \end{align}
 that is the same as to say that \eqref{eq:L2*estimate} holds.\finedim

\bigskip

We will need the following crucial boundary version of the Non Vanishing Lemma.
Here we gather some ideas from \cite{FW}.
\ble[The Boundary Non Vanishing Lemma]\label{vanishle-bdy}
Let $\vn$ satisfy all the assumptions of Theorem \ref{thm:spike-vn.intro} and
assume that there exists a sequence
$\{x_n\}\subset \om$ such that $\vn(x_n)>1$ and dist$(x_n,\pa \om)\to 0$.
Then
\beq\label{bdy:dist}
\frac{\mbox{\rm dist}(x_n,\pa \om)}{\eps_n}\to +\ii,
\eeq
and, possibly along a subsequence,
there exists $z_0\in \pa \om$ and $\{z_n\}\subset \om$ such that
$z_n\to z_0$, $|z_n-x_n|\leq C\eps_n$
and $w_n(x)=\vn(z_n+\eps_nx)\to w_0(x)$
in $C^{2}_{\rm loc}(\R^N)$. In particular $\vn(z_n)\to w_0(0)>1$, i.e.
$\vn(z_n)$ stays uniformly bounded far away from $1$.
\ele
\proof We split the proof in four steps. \\
{\bf STEP 1}: We first show that the rescaled functions
        \begin{align}\nonumber
         w_n(y)\coloneqq \vn(x_n+\eps_n y), \quad y\in \Omega_n:=
         \frac{\Omega-x_n}{\eps_n},
        \end{align}
        is uniformly bounded in~$C^{1,s}$ norm, for some~$s\in (0,1)$.

        \bigskip

We adapt an argument in \cite{FW}. Note that $w_n$ satisfies
\beq\label{eq:wn.1}
        \graf{-\Delta w_n =[w_{n}-1]_+^p\quad \mbox{in}\;\;\om_n,\\ \\
        w_{n}= 0 \quad \mbox{on}\;\;\pa\om_n, \\ \\
        w_n >0 \quad \mbox{ in } \;\;\Omega_n,
        }
        \eeq
and
        $$
         \int_{\Omega_n} [w_n-1]_+^p \dd  y= \frac{1}{\eps_n^N} \int_{\Omega}[v_n-1]_+^p\dd x
         = \frac{1}{\eps_n^N |\alpha_n|^p }
         =\parenthesis{\frac{\lambda_n}{|\alpha_n|^{1-\frac{p}{p_{_N}}}}}^{\frac{N}{2}} \leq C_0.
        $$
 In view of \eqref{eq:Dirichlet bound-vn} we have,
        \begin{align}\nonumber
         \int_{\Omega_n} |\nabla w_n|^2 \dd y
         =\eps_n^{2-N}\int_{\Omega}  |\nabla v_n|^2 \dd x \leq C_0+C_1,
        \end{align}
        and consequently by either the Sobolev embedding or by \eqref{eq:L2*estimate},
        \begin{align}\nonumber
         \int_{\Omega_n } |w_n|^{2^*}\dd y
         \leq \mathcal{C}_S^{2^*}\parenthesis{\int_{\Omega_n} |\nabla w_n|^2\dd y}^{\frac{2^*}{2}}
         \leq \mathcal{C}_S^{2^*}\parenthesis{C_0+C_1}^{\frac{N}{N-2}}.
        \end{align}
         Let $$q_0=q_0(N,p):=\frac{2^*}{p},$$ and observe that, because of $p<p_{_N}$, we have $q_0>2$.
        Since $w_n\ge 0$, the Calderon-Zygmund inequality implies that,
        \begin{align}\label{CZ}\nonumber
         \|D^2 w_n\|_{L^{q_0}(\Omega_n)}
         \leq& C_{CZ}(N,q_0)\|\Delta w_n\|_{L^{q_0}(\Omega_n)}
         =C_{CZ}(N,q_0)\|[w_n-1]_+^p \|_{L^{q_0} (\Omega_n)} \\\nonumber
         \leq&C_{CZ}(N,q_0)\|w_n^p \|_{L^{q_0} (\Omega_n)} = C_{CZ}(N,q_0)\| w_n\|_{L^{2^*}(\Omega_n)}^p \\
         \leq & C_{CZ}(N,q_0)\parenthesis{\mathcal{C}_S\sqrt{C_0+C_1}}^p.
        \end{align}
        Remark that all the constants in the estimates do not depend on $n$ and assume for the moment that
        \beq\label{q0}
        q_0>N,
        \eeq
        whence in particular $2-\frac{N}{q_0}=1+ s$ for some $s\in (0,1)$.  We will use the symbol $C$ to denote various constants that do not depend on $n$
        and which may change from line to line. We claim that the $C^{1,s}(\om_n)$ norm of $w_n$ is uniformly bounded. In fact, setting $f_i=\frac{\pa w_n}{\pa x_i}$, by standard techniques one can prove that
        $f_i$ can be extended to $\tilde f_i$ on $\R^N$ in such a way that
        $$\|\tilde f_i\|_{L^{2}(\R^N)}\leq C(N)\|f_i\|_{L^{2}(\Omega_n)}\leq C(N)\|\nabla w_n\|_{L^{2}(\Omega_n)},$$ and
        $$\|\nabla \tilde f_i\|_{L^{q_0}(\R^N)}\leq C(N,q_0) \|\nabla f_i\|_{L^{q_0}(\Omega_n)}\leq C(N,q_0)\|D^2 w_n\|_{L^{q_0}(\Omega_n)}, $$
        where $C(N),C(N,q_0)$ depend in general also on the $L^{\ii}$ norm of the gradient and the Hessian of the functions defining the local charts on $\pa \om_n$,
        which however are readily seen to provide a uniformly bounded contribution. Now by a classical argument due to Morrey
        we have,
        $$
        |f_i(x)-f_i(y)|\leq C\|\nabla \tilde f_i\|_{L^{q_0}(\R^N)}|x-y|^{1-\frac{N}{q_0}},\forall |x-y|\leq 1,
        $$
        and, for $x\in \R^N$ and $Q$ any cube of side length $1$ containing $x$,
        $$
        |f_i(x)|\leq \frac{1}{|Q|}\int\limits_{Q}|f_i| + C\|\nabla \tilde f_i\|_{L^{q_0}(\R^N)}\leq
        \|\tilde f_i\|_{L^{2}(\R^N)}+C\|\nabla \tilde f_i\|_{L^{q_0}(\R^N)},
        $$
        implying that $\|\nabla w_n\|_{C^{0,s}(\om_n)}\leq C_s$, for some $C_s$ depending only on $N,p$.
        Therefore $w_n\in C^{1,s}(\om_n)$ and we have, for any $x\in \R^N$,
        $$
        |w_n(x)|\leq \frac{1}{|B_1(x)|}\int\limits_{B_1(x)}|w_n|+\int\limits_{B_1(x)}|w_n(x)-w_n(y)|dy\leq C\|w_n\|_{L^{2^*}(\Omega_n)}+2CC_s,
        $$
        implying that $\|w_n\|_{C^{1,s}(\om_n)}\leq C$, for some $C$ depending only on $N$ and~$p$, as claimed.\\
        On the other side, if $q_0\leq N$, assume first that,
        \beq\label{q03}
        \frac{N}{2}\leq q_0< N.
        \eeq
        Thus, for $f_i$ and $\tilde{f}_i$ defined as above, since $+\ii>q_0^*=\frac{N q_0}{N-q_0}\geq  N\geq 3$,
        by the Sobolev embedding and the Young inequality we have
        \begin{align}\label{q02}
         \|\tilde f_i\|_{L^{t}(\R^N)}
         \leq & \|\tilde f_i\|_{L^{2}(\R^N)} + \|\tilde f_i\|_{L^{q_0^*}(\R^N)}
         \leq \|\tilde f_i\|_{L^{2}(\R^N)}+ C(N,q_0)\|\nabla \tilde f_i\|_{L^{q_0}(\R^N)} \\  \nonumber
         \leq & C(N)\|\nabla w_n\|_{L^{2}(\Omega_n)}+C(N,q_0)\|D^2 w_n\|_{L^{q_0}(\Omega_n)},\quad t\in [N,q_0^*].
        \end{align}
        In particular we have that $\|\nabla w_n\|_{L^{N}(\om_n)}$ is uniformly bounded.
        Since~$w_n\in W^{1,N}_0(\om_n)$, extending $w_n$ by zero to $\tilde w_n\in W^{1,N}(\R^N)$, we have that, for any $s\in [N,+\ii)$, $\|w_n\|_{L^{s}(\om_n)}=\|\tilde w_n\|_{L^{s}(\R^N)}$ is uniformly bounded. Therefore, for any
        $q_1>Np$ we deduce as in \rife{CZ} that
        $$
        \|D^2 w_n\|_{L^{q_1}(\Omega_n)}
         \leq C_{CZ}(N,q_1)\|\Delta w_n\|_{L^{q_1}(\Omega_n)}\leq C_{CZ}(N,q_1)\| w_n\|_{L^{q_1 p}(\Omega_n)}^p \leq C,
        $$
        for some $C$ depending only on $N,q_1,p$. Therefore we are reduced to the same case as in \rife{q0}. Obviously, if $q=N$,
        the only difference we have is that $q_0^*=+\ii$, whence the same argument with minor changes applies. Therefore we are left with the case
        $$
        2<q_0<\frac{N}{2},
        $$
        where we have $2^*<q_0^*<N$. Since $2^*>2$ we can still use \rife{q02} to deduce that for any $t\in[2^*,q_0^*]$
        $\|\tilde f_i\|_{L^{t}(\R^N)}$ is uniformly bounded, that is, $\|\nabla w_n\|_{L^{t}(\Omega_n)}$ is uniformly bounded and then, again since
        $w_n\in W^{1,t}_0(\om_n)$, extending $w_n$ by zero to $\tilde w_n\in W^{1,t}(\R^N)$, we have that $\|w_n\|_{L^{q_0^{**}}(\om_n)}$ is uniformly bounded, where
        $$
        \frac{2N}{N-4}=2^{**}<q_0^{**}=\frac{N q_0}{N-2q_0}<+\ii.
        $$
        Remark that $2<q_0<\frac{N}{2}$ implies $N\geq 5$, whence the left hand sided inequality is well defined.
        At this point we start a Moser-type iteration putting,
        $$
        q_{k+1}=\frac{1}{p}\frac{N q_k}{N-2q_k}>\left(\frac{N-2}{N-2q_k}\right)q_k,\;k\geq 0,
        $$
        where we used $p<p_{_N}$, which defines a strictly increasing sequence such that $q_{k+1}<+\ii$ as far as $q_k<\frac{N}{2}$.
        It is readily seen that after a finite number of steps we will have $q_k\geq \frac{N}{2}$ and we are reduced to the same case as in
        \rife{q03}. Therefore we deduce that the $C^{1,s}(\om_n)$ norm of $w_n$ is uniformly bounded for some $s\in(0,1)$.

        \

        {\bf STEP 2}: There exists a constant~$\gamma_0>0$ such that~
        $\gamma_n:=\frac{\dist(x_n, \pa\om_n)}{\eps_n}\ge \gamma_0$.

        \

        As~$\pa\Omega$ is compact and of class~$C^{2,\beta}$, there exists~$\hat{x}_n\in\pa\om$
        such that~$|x_n-\hat{x}_n|=\dist(x_n,\pa\om)$.
        Therefore, putting $\hat{x}_n=x_n+\eps_n\hat{y}_n$, we have $|\hat{y}_n|=\gamma_n$ and by STEP 1,
        \begin{align}\nonumber
         1\leq v_n(x_n)- v_n(\hat{x}_n) = w_n(0)- w_n(\hat{y}_n) \leq
         \|\nabla w_n\|_{L^\infty(\Omega_n)} |\hat{y}_n|\leq C_{N,p} \gamma_n,
        \end{align}
        that is,~$\gamma_n\geq C_{N,p}^{-1}=: \gamma_0$.

        \

        {\bf Step 3} $\gamma_n\to +\ii$.

        The proof goes exactly as in \cite{FW} p.693 and we include it here
        for reader's convenience.\\
        Suppose on the contrary that~$\gamma_n$ (sub)converges to~$\gamma \in (0,+\infty)$.
        Note that in this case, after suitably defined rigid motions,
        $\Omega_n$ converges to a translated half space of the form,
        \begin{align}\nonumber
         \R^N_{\gamma}\coloneqq \braces{ y\in \R^N\;\mid \; y^N < \gamma}.
        \end{align}
        By STEP 1 ~$w_n$ converges in~$C^{2}_{\rm loc}$ to a function~
        $w\in \mathcal{D}^{1,2}(\R^N_\gamma)\cap C^2_{\rm loc}(\R^N_\gamma)$ which satisfies
         \beq\nonumber
        \graf{-\Delta w =[w-1]_+^p\quad \mbox{in}\;\;\R^N_\gamma,\\ \\
        w= 0 \quad \mbox{on}\;\;\pa\R^N_\gamma.
        }
        \eeq
        According to Theorem I.1 in \cite{EL} we have ~$w\equiv 0$ which is a contradiction to
        $\liminf\limits_{n\to +\ii}w_n(0)=\liminf\limits_{n\to +\ii}v_n(x_n)\geq 1$.\\

        {\bf STEP 4} In view of STEP 3 we have that \eqref{bdy:dist} holds and then we can argue
        as in the proof of the Vanishing Lemma. Indeed, since \eqref{mass.intro.1} and
        \eqref{eq:L2*estimate} hold true, then we have all the assumptions needed to follow
        step by step in that proof, we skip the details to avoid repetitions.
        \finedim

\bigskip

Of course there is a useful almost equivalent version of the Boundary Non Vanishing Lemma whose
proof is obtained by the same argument.
\ble[The Boundary Vanishing Lemma]
Let $\vn$ satisfy all the assumptions of Theorem \ref{thm:spike-vn.intro} and
assume that there exists $z_0\in \pa \om$ and $r>0$ such that
$[\vn-1]_+\to 0$ uniformly in $\om\cap B_{2r}(z_0)$.
Then there exists ${n}_r\in \N$ and
$C_r>0$,
such that $[v_n-1]_+=0$ in $\om\cap \ov{B_{r}(z_0)}$
for any $n>n_r$.
\ele
\proof By contradiction there exists $\{x_n\}\subset \om\cap B_{2r}(z_0)$ such that $v_n(x_n)>1$ and then along
a subsequence we can assume w.l.o.g. that $x_n \to x_0 \in \ov{\om}\cap\ov{B_{r}(z_0)}$.
If $x_0\notin \pa \om$
we obtain a contradiction by the Non Vanishing Lemma, see Remark \ref{vanish-equiv}, otherwise,
if $x_0\in \pa \om$ we obtain a contradiction by the Boundary Non Vanishing Lemma. Therefore
$[v_n-1]_+=0$ in $\om\cap {B_{r}(z_0)}$ for any $n$ large enough.\finedim

\medskip

At this point we are ready to provide the proof of Theorem \ref{thm:spike-vn.intro}.\\
{\bf The proof of Theorem \ref{thm:spike-vn.intro}}\\
{\bf STEP 1} Remark that, in view of Lemma \ref{lem:2star}, \eqref{eq:L2*estimate} holds, which is
\eqref{mass.intro.2} with $t=2^*$. Assume first that $[\vn-1]_+ \to 0 \mbox{ uniformly in } \Omega$. By a covering argument and
the Vanishing Lemma and the Boundary
Vanishing Lemma, we deduce that $[\vn-1]_+=0$ in $\om$ for $n$ large enough, which contradicts
$\ino [\vn-1]^p_+=|\an|^{-p}$. Therefore there exists $x_n\subset \om$ such that, possibly
along a subsequence which we will not relabel, $x_n \to z \in \ov{\om}$,
$\vn(x_n)=\max\limits_{\om}\vn$ and $\vn(x_n)\geq 1+\sg$, for some $\sg>0$.
In view of \eqref{mass.intro.1}, \eqref{mass.intro.2} (with $t=2^*$), the Non Vanishing Lemma
and the Boundary Non Vanishing Lemma we can argue as in the proof of Theorem \ref{thm5.7.intro}
to deduce that, possibly along a further subsequence, $\vn$ is an $\mathcal{Z}$-spikes sequence for some positive 0-chain~$\mathcal{Z}=\sum_{j=1}^m k_j z_j$ in~$\ov{\om}$, with spikes set~$\Sigma=\Sigma(\mathcal{Z})$.
We skip details to avoid repetitions.
Next we show that \eqref{est:nm2} holds and that if $\Sigma\subset \om$, then $(z_1,\cdots,z_m)$ is a critical point of the $\vec{k}$-Kirchoff-Routh Hamiltonian, with~$\vec{k}=(k_1,\cdots, k_m)$.

\

We split the proof in two steps. We first prove \eqref{est:nm2}.\\
{\bf STEP 2}: Recall that ${\Sigma}_{r}\coloneqq\underset{j =1,\cdots,m}{\cup}B_{r} (z_{i})$.
By the Green representation formula, for any $x\in \om\setminus \Sigma_r$ and for any $n$ large
enough we have,
\begin{align}
 \eps_n^{-(N-2)}\vn(x)-& \sum\limits_{i=1}^{N_m}M_{p,0}G(x,z_{n,i})
 = \ino G(x,y)\mu_n^{\frac N 2}[\vn-1]^p_+-\sum\limits_{i=1}^{N_m}M_{p,0}G(x,z_{n,i}) \nonumber  \\
 =& \sum\limits_{i=1}^{N_m}\left(\int\limits_{B_{r/2}(z_{n,i})\cap \om} G(x,y)\mu_n^{\frac N 2}[\vn(y)-1]^p_+dy-M_{p,0}G(x,z_{n,i})\right) \nonumber \\
 =& \sum\limits_{i=1}^{N_m}
\left(\int\limits_{B_{2R_0\eps_n}(z_{n,i})} G(x,y)\mu_n^{\frac N 2}[\vn(y)-1]^p_+dy-M_{p,0}G(x,z_{n,i})\right),
\nonumber
\end{align}
where we used \eqref{bdy:dist} and for fixed $i\in\{1,\cdots,N_m\}$, $z_{n,i}\to z_j$ for some $j\in \{1,\cdots,m\}$.
It is obviously enough to estimate the argument in the round parenthesis here above. Thus, for any fixed $i\in \{1,\cdots,N_m\}$, we see that
\begin{align}
 \int\limits_{B_{2R_0\eps_n}(z_{n,i})}
  & G(x,y)\mu_n^{\frac N 2}[\vn(y)-1]^p_+dy
    - M_{p,0}G(x,z_{n,i}) \nonumber \\
 = & \int\limits_{B_{2R_0}(0)} G(x,z_{n,i}+\eps_n z)[w_n(z)-1]^p_+dz-M_{p,0}G(x,z_{n,i}) \nonumber \\
 = &\int\limits_{B_{2R_0}(0)} \left(G(x,z_{n,i}+\eps_n z)[w_n(z)-1]^p_+-G(x,z_{n,i})[w_0(z)-1]^p_+\right)dz \nonumber \\
 \leq & \int\limits_{B_{2R_0}(0)} \left|G(x,z_{n,i}+\eps_n z)-G(x,z_{n,i})\right|[w_n(z)-1]^p_+dz \nonumber \\
      &\qquad\quad   +\int\limits_{B_{2R_0}(0)} G(x,z_{n,i})\left|[w_n(z)-1]^p_+-[w_0(z)-1]^p_+\right|dz \nonumber \\
 = & O(\eps_n)+O(\|[w_n(z)-1]^p_+-[w_0(z)-1]^p_+\|_{\infty})=o(1), \nonumber
\end{align}
which readily implies \eqref{est:nm2}.\\

{\bf STEP 3}: Here we adapt the argument in \cite{MaW}. By assumption
$\Sigma\subset \om$, whence,
from the first part of the proof and $(viii)$ of Definition \ref{spikeseqdef},
we have that,
\begin{align}\nonumber
 \mu_n^{\frac{N}{2}} [ v_n-1]_+^p \rightharpoonup M_{p,0}\sum_{j=1}^m m_j \delta_{z_j}
\end{align}
weakly in the sense of measures in $\om$,
where~$m_j$ denotes the multiplicity of the spike point~$z_j$.
Let us set,
\begin{align}\nonumber
 u_n(x) \coloneqq \mu_n^{\frac{N-2}{2}}v_n(x),
\end{align}
which satisfies,
\begin{align}\label{eq:for un}
 -\Delta u_n=\mu_n^{\frac{N}{2}}[v_n-1]_+^p \; \mbox{ in }\om
\end{align}
and~$u_n=0$ on~$\pa\Omega$.
It is well-known that
\begin{align}\label{eq:local convergence for un.1}
 u_n \to M_{p,0}\sum_{j=1}^m m_j G(x,z_j)\eqqcolon \bar{G}(x) \quad \mbox{ in } \; C^{2,\beta}_{\rm loc}(\overline{\Omega}\setminus \Sigma).
\end{align}

This observation helps to find the locations of the spike points.

Testing~\eqref{eq:for un} against the vector~$\nabla u_n$, we get a Pohozaev identity for vectors: for any~$\Omega'\subseteq \Omega$,
\begin{align}\label{eq:Pohozaev}
 \int_{\pa\Omega'} -(\pa_\nu u_n)\nabla u_n +\frac{1}{2}|\nabla u_n|^2 \nu
 =\int_{\pa\Omega'} \frac{\mu_n^{N-1}}{p+1}[v_n-1]_+^{p+1}\nu.
\end{align}
Now without loss of generality we fix the spike point~$z_1\in \Sigma$, and consider~$\Omega'= B_t(z_1)$
for sufficiently small~$t>0$ (so that~$\overline{B}_{t}(z_1) \cap \Sigma = \braces{z_1}$ ).
Due to~\eqref{eq:local convergence for un.1} the LHS of \eqref{eq:Pohozaev} satisfies:
\begin{align}\nonumber
 \int_{\pa B_t(z_1)} -(\pa_\nu u_n)\nabla u_n +\frac{1}{2}|\nabla u_n|^2 \nu
 \to \int_{\pa B_t(z_1)} -(\pa_\nu \bar{G})\nabla \bar{G} +\frac{1}{2}|\nabla \bar{G}|^2 \nu
 \quad \mbox{as } n\to +\infty,
\end{align}
meanwhile the RHS in~\eqref{eq:Pohozaev}, because of \eqref{est:nm2},
is readily seen to vanish in the limit:
\begin{align}\nonumber
 \int_{\pa B_t(z_1)} \frac{\mu^{N-1}_n}{p+1}[v_n-1]_+^{p+1}\nu=0, \quad \mbox{for any $n$ large enough}.
\end{align}
Hence, passing to the limit (along a subsequence if necessary) we deduce that,
\begin{align}\nonumber
  \int_{\pa B_t(z_1)} -(\pa_\nu \bar{G})\nabla \bar{G} +\frac{1}{2}|\nabla \bar{G}|^2 \nu=0.
\end{align}

\

Observe that, putting $r(x)=|x-z_1|$,
\begin{align}\nonumber
 \bar{G}(x)
 =& M_{p,0} \sum_{j=1}^m m_m G(x, z_j)
 =M_{p,0}\parenthesis{ m_1 G(x, z_1) + \sum_{j\neq 1} m_j G(x, z_j)} \\ \nonumber
 =&M_{p,0}\parenthesis{\frac{m_1}{N(N-2)\omega_N} \frac{1}{r^{N-2} } + m_1 H(x, z_1) +
 \sum_{j\neq 1} m_j G(x, z_j)} \\ \nonumber
 =& M_{p,0}\parenthesis{\frac{m_1}{N(N-2)\omega_N} \frac{1}{r^{N-2} } + F_1(x) }
\end{align}
where in the last step we have introduced
\begin{align}\nonumber
F_1(x):= m_1 H(x, z_1) + \sum_{j\neq 1} m_j G(x, z_j).
\end{align}
To simplify the notations let us also write,
\begin{align}\nonumber
 \widetilde{G}(x)\coloneqq M_{p,0}^{-1} \bar{G}(x)= \frac{m_1}{N(N-2)\omega_N} \frac{1}{r^{N-2} } + F_1(x),
\end{align}
which clearly satisfies
\begin{align}\label{eq:tildeG integral.1}
  \int_{\pa B_t(z_1)} -(\pa_\nu \widetilde{G})\nabla \widetilde{G} +\frac{1}{2}|\nabla \widetilde{G}|^2 \nu=0.
\end{align}
At this point, since~$\nu=\nabla r$, we can compute this sum term by term as follows,
\begin{align}\nonumber
 \nabla\widetilde{G}(x)
 = -\frac{m_1}{N\omega_N}\frac{1}{r^{N-1}}\nabla r + \nabla F_1(x), & &
 \nabla_\nu \widetilde{G}
 =-\frac{m_1}{N\omega_N}\frac{1}{r^{N-1}} + \nabla_{\nu} F_1(x),
\end{align}
\begin{align}\nonumber
 |\nabla\widetilde{G}|^2
 =& \frac{m_1^2}{(N\omega_N)^2}\frac{1}{r^{2(N-1)}}
    -\frac{2m_1}{N\omega_N}\frac{\nabla_r F_1}{r^{N-1}}
    + |\nabla F_1|^2, \\ \nonumber
 (\pa_\nu \widetilde{G})\nabla\widetilde{G}
 =&\frac{m_1^2}{(N\omega_N)^2}\frac{\nabla r}{r^{2(N-1)}}
   -\frac{m_1}{N\omega_N}\frac{\nabla F_1+ (\nabla_\nu F_1)\nabla r }{r^{N-1}}
   + (\nabla_\nu F_1)\nabla F_1
\end{align}
Thus,
\begin{align}
  -(\pa_\nu \widetilde{G})\nabla \widetilde{G} +\frac{1}{2}|\nabla \widetilde{G}|^2 \nu
  =& -\frac{m_1^2}{(N\omega_N)^2}\frac{\nabla r}{r^{2(N-1)}}
   +\frac{m_1}{N\omega_N}\frac{\nabla F_1+ (\nabla_\nu F_1)\nabla r }{r^{N-1}}
   - (\nabla_\nu F_1)\nabla F_1  \nonumber\\
   & +\frac{1}{2} \frac{m_1^2}{(N\omega_N)^2}\frac{\nabla r}{r^{2(N-1)}}
    -\frac{m_1}{N\omega_N}\frac{\nabla_r F_1}{r^{N-1}}\nabla r
    + \frac{1}{2}|\nabla F_1|^2\nabla r \nonumber\\
 =& -\frac{1}{2}\frac{m_1^2}{(N\omega_N)^2}\frac{\nabla r}{r^{2(N-1)}}
    +\frac{m_1}{N\omega_N}\frac{\nabla F_1}{r^{N-1}}
    +\frac{1}{2}|\nabla F_1|^2\nabla r - (\nabla_r F_1)\nabla F_1.\nonumber
\end{align}
Inserting this expression into the identity~\eqref{eq:tildeG integral.1}, we see that
\begin{align}
 \nabla F_1(z_1)=0. \nonumber
\end{align}
Therefore, the spike point~$z_1$ is contained in the critical set of the function~$F_1$
and this obviously holds for any~$j=1,\cdots, m$.\\

{\bf STEP 4} At last we prove the claim about the case where $\om$ is convex.
We first prove that if $\om$ is convex then we can drop \rife{eq:assump-2.intro} while all the properties
proved so far in STEPS 1,2,3 hold true and moreover $\Sigma=\{z_1\},z_1\in \om$.\\
It is useful to recall the notion of interior parallel
sets in~\cite{Bandle}: for each~$t>0$, the interior~$t$-parallel set of~$\Omega$ is
\begin{align}\nonumber
 \Omega_{-t} \coloneqq \braces{x\in \Omega\; \mid \; B_t(x)\subset \Omega}
 =\braces{x\in \Omega \; \mid \; \dist(x, \Omega^c)>t}.
\end{align}
In view of $\int\limits_{\om} {|\an|^p[\vn-1]_+^p}=1$ we have,
$$
\mu_n^{\frac N 2}\int\limits_{\om} {[\vn-1]_+^p}=
\left(\frac{\lm_n}{|\an|^{1-\frac{p}{p_{_N}}}}\right)^{\frac N 2}
\leq C_0.
$$
Therefore Theorem \ref{thm5.7} implies that $\vn$ is locally uniformly bounded.
However since $\om$ is a convex set of class $C^{2,\beta}$, then it is a well known consequence of the moving
plane method (see \cite{gnn} and also \cite{DLN} p.45) that there exists $t_0>0$ do not
depending on $n$, such that $\vn(x-t\nu(x))$ is nondecreasing for $t\in[0,t_0]$ for any
$x\in \pa\om$, where $\nu(x)$ is the unit outer normal at $x$.
In particular we have,
\beq\label{bdy.h}
\sup \limits_{\om\setminus \om_{-t_0/2}}\vn\leq \sup \limits_{\pa \om_{-t_0/2}}\vn.
\eeq
Thus, from Theorem \ref{thm5.7}, we have in particular that $\vn$ is uniformly bounded in $L^{\ii}(\om)$.
 As a consequence we deduce that
$$
\ino |\nabla \vn|^2=\mu_n \ino [\vn-1]_+^p\vn=
\frac{\mu_n^{\frac N 2}}{\mu_n^{\frac{N-2}{2}}}\ino [v_n-1]_+^p\vn\leq
\frac{\|\vn\|_{L^{\ii}(\om)}}{\mu_n^{\frac{N-2}{2}}}C_0\leq C\mu_n^{-\frac{N-2}{2}}.
$$
Therefore the conclusions of Lemma \ref{lem:2star} hold, since by the Sobolev embedding we have
that also \eqref{eq:L2*estimate} holds true.
As a consequence, also \rife{mass.intro.2} is satisfied with $t=2^*$ and then we can argue as in {STEP 1}
to deduce that $\vn$ is a $\mathcal{Z}$-spikes sequence for some $\mathcal{Z}=\sum_{j} k_j z_j$ in $\ov{\om}$. However, since $\om$ is convex,
by a classical argument based on the moving plane method (\cite{gnn}), there can be no critical points
of $\vn$ in a sufficiently small uniform neighborhood of the boundary, say $\om\setminus \om_{-t_0/2}$, implying
by the Boundary Non Vanishing Lemma that there can be no spikes either.
Therefore $\Sigma(\mathcal{Z})\subset \om$ and
by the result in \cite{GrT}, there is no critical point of $\mathcal{H}(x_1,\cdots,x_m;\vec{k})$ as far
as $\# \mathcal{Z}=m\geq 2$, whence $\Sigma=\{z_1\}$ and $z_1\in \om$ is by STEP 3 a critical point of $\mathcal{H}_1(x_1;(1)_1)$, which
is just the defining equation of an harmonic center of $\om$. However, according to \cite{CT}, there is only one harmonic center as far as
$\om$ is convex. The proof of Theorem~\ref{thm:spike-vn.intro} is complete.
\finedim

\

Next we have,\\
{\bf The Proof of Theorem \ref{thm7.2.intro}}\\
Let $\vn=\frac{\lmn}{|\an|}\pn$ which obviously satisfies \rife{vn.1.intro}.
Also, since
$$
\mu_n^{\frac N 2}\int\limits_{\om}[\vn-1]^p_+=\frac{\mu_n^{\frac N 2}}{|\an|^p}=
\left(\frac{\lm_n}{|\an|^{1-\frac{p}{p_{_N}}}}\right)^{\frac N 2},
$$
then \eqref{hyp.intro.pn.1} implies that
\rife{mass.intro.1} holds true. Next assume \eqref{hyp.intro.pn},  then we also have,
$$
\mu_n^{\frac N 2}\int\limits_{\om}[\vn-1]^{p+1}_+=
\frac{\mu_n^{\frac N 2}}{|\an|^{p+1}}\int\limits_{\om} [\an+\lm \pn]_+^{p+1}\leq C_2 \frac{\mu_n^{\frac N 2}}{|\an|^{p}}\leq C_2C_0.
$$
Therefore
\rife{eq:assump-2.intro} holds true as well and consequently the conclusions of Theorem \ref{thm:spike-vn.intro}
hold for $\vn$ in $\om$, as claimed. This shows at once that the mass quantization identity holds as well as,
 in view of $(vi)$ in Definition \ref{spikeseqdef}, the asymptotic ``round'' form of $\om_{n,+}$.
Also observe that \eqref{psi:est.1} just follows from the definition of $\om_{n,+}$ and
$\pn=\frac{|\an|}{\lmn}\vn$. Concerning \eqref{psi:est.2}, observe that, for $x\in \om\setminus \Sigma_r$,
from \eqref{est:nm2} we have
$$
\pn(x)=\frac{|\an|}{\lmn}\vn(x)=\frac{|\an|}{\lmn \mu_n^{(N-2)/2}}\sum\limits_{i=1}^{|\mathcal{Z}|}M_{p,0}G(x,z_{n,i})+o(\frac{|\an|}{\lmn}\mu_n^{-(N-2)/2})
$$
and then the conclusion follows from the mass quantization since we have,
\beq\label{psiconc.1}
\frac{|\an|}{\lmn \mu_n^{(N-2)/2}}=\frac{|\an|}{\lmn (\lm_n |\an|^{p-1})^{(N-2)/2}}=\left(\frac{|\an|^{1-\frac{p}{p_{_N}}}}{\lm_n}\right)^{\frac N 2}=\frac{(1+o(1))}{|\mathcal{Z}| M_{p,0}},\;n\to +\infty.
\eeq

At this point \eqref{psi:est.0} follows similarly by \eqref{psiconc.1} and
$(viii)$ of Definition \ref{spikeseqdef} whenever we observe that,

$$
[\an+\lm_n\pn]_+^p=|\an|^p[\vn-1]_+^p=\frac{|\an|}{\lm_n}\lm_n|\an|^{p-1}[\vn-1]_+^p=
\frac{|\an|}{\lm_n\mu_n^{(N-2)/2}}\mu_n^{N/2}[\vn-1]_+^p.
$$

The assertions about the case where $\om$ is convex follows immediately from Theorem \ref{thm:spike-vn.intro}.
\finedim

\bigskip
\bigskip

Next we have,\\
{\bf The proof of Corollary \ref{exist:plasma}.}\\
\emph{Proof of} {\rm (a)}: It has been proved in \cite{FW} that there exists a family of solutions $v_\eps$ of
$$
\graf{-\Delta v =\mu [v-1]_+^p\quad \mbox{in}\;\;\om,\\
\eps^{-2}:=\mu\to +\ii,\;\\
v=0 \quad \mbox{on}\;\;\pa \om
}
$$

such that (see Step 8 in \cite{FW})
$$
\ino |\nabla v_\eps|^2\leq C \eps^{N-2},
$$

and (see Step 12 in \cite{FW}), for any fixed $R\geq R_0$, $v_\eps(z_{\eps}+\eps y)=w_0(y)+o(1)$ for $|y|\leq R$, where $z_\eps$ is the unique critical point of $v_\eps$ and
$z_\eps\to z_1$ where $z_1$ is an harmonic center of $\om$. On the other side, by using the Sard Lemma as in Lemma \ref{enrg0} and
testing the equation with $v_\eps-1$ on $\om_+=\{v_\eps-1>0\}$,
$$
\int\limits_{\om_+}|\nabla v_\eps|^2=\mu \int\limits_{\om_+}[v_\eps-1]_+^{p+1},
$$
while testing the equation with $v_\eps$ on $\om_-=\{v_\eps-1<0\}$,
$$
0=\int\limits_{\om_-}|\nabla v_\eps|^2-\int\limits_{\pa \om_-}v_\eps \pa_\nu v_\eps=
\int\limits_{\om_-}|\nabla v_\eps|^2+\int\limits_{\pa \om_+}\pa_\nu v_\eps=\int\limits_{\om_-}|\nabla v_\eps|^2-\mu\int\limits_{\om_+}[v_\eps-1]_+^p.
$$

Therefore we have,

\beq\label{2hyp}
\int\limits_{\om_+}[v_\eps-1]_+^{p+1}+\int\limits_{\om_+}[v_\eps-1]_+^{p}=\eps^{2}\int\limits_{\om}|\nabla v_\eps|^2\leq C\eps^N.
\eeq

At this point we define $\al_\eps<0$ and $\lm_\eps$ as follows,
$$
|\al_\eps|^p\int\limits_{\om_+}[v_\eps-1]^p=1,\qquad \mu=\eps^{-2}=\lm_\eps |\al_\eps|^{p-1},
$$
so that, in view of \eqref{2hyp}, we have,
$$
\frac{1}{|\al_\eps|^p}=\int\limits_{\om_+}[v_\eps-1]^p\leq C\eps^N.
$$
As a consequence we see that $|\al_\eps|\to +\ii$ and $\lm_\eps|\al_\eps|^{-1}\leq C\eps^{N-2}\to 0$. At this point, defining $\psi_\eps$ as
$v_\eps=\frac{\lm_\eps}{|\al_\eps|}\psi_\eps$, it is readily seen that $(\al_\eps,\psi_\eps)$ is a solution of $\prl$ with $\lm=\lm_\eps$ and we deduce from Theorem \ref{thm:7.1}
that necessarily $\lm_\eps \to +\ii$. In particular, in view of \eqref{2hyp} we see that both \eqref{mass.intro.1} and
\eqref{eq:assump-2.intro} are satisfied. It follows from Theorem \ref{thm:spike-vn.intro} that in fact as $\eps_n\to 0$, $v_{\eps_n}$ is a
$\mathcal{Z}$-spikes sequence where $\mathcal{Z}=z_1$ and $z_1$ is an harmonic center of $\om$.

\medskip

\emph{Proof of} {\rm (b)}. Under the assumptions in the claim, it is shown in \cite{We} that, for $\gamma$ and $I$ as in $\fbi$,
there exists a family $v_\eps$ of solutions of
\beq
\graf{-\Delta v =\mu [v-1]_+^p\quad \mbox{in}\;\;\om,\nonumber \\
\eps^{-2}:=\mu=|\gamma|^{p-1},\mu \to +\ii,\; \nonumber\\
v=0, \quad \mbox{on}\;\;\pa \om \nonumber\\
\mu^{\frac{p}{p-1}}\ino [v-1]_+^p=I
}
\eeq
which satisfy $v_\eps(x)=w_\eps(x)+O(\eps^{N-2})$, where
$$
w_\eps(x)=\sum\limits_{j=1}^m w_0(\eps^{-1}(x-z_j)).
$$
In particular one has
$$
I=\mu^{\frac{p}{p-1}}\ino [v_\eps-1]_+^p=\eps^{-\frac{2p}{p-1}}\ino [v_\eps-1]_+^p=\eps^{-\frac{2p}{p-1}}\eps^N (mM_{p,0}+O(\eps^{N-2})),
$$
as $\mu\to +\ii$. Recalling that $\gamma=\lm^{\frac{1}{p-1}}\al$, $I=\lm^{\frac{p}{p-1}}$ we see that in fact $v_\eps$ solves \eqref{w5.1.intro} and
that in fact $\vn:=v_{\eps_n}$, $\eps_n^{-2}=\mu_n\to +\ii$ is a $\mathcal{Z}$-spikes sequence where $\mathcal{Z}=\sum\limits_{j=1}^mz_j$. Therefore $\pn$ defined
as in the claim is a $\mathcal{Z}$-blow up sequence.
\finedim

\bigskip
\bigskip

Next we have,\\
{\bf The proof of Theorem \ref{thm:varsol}.}\\
We take the same notations about $\mathcal{F}_{\lm}(\rh)$ and $\mathcal{J}_{\lm}(u)$ as in the Appendix below.
If the claim were false we would have, possibly along a subsequence, $\tau_n:=|\an|^{-1}\ino [\an+\lm_n \pn]_+^{p+1}\to \tau >\frac{p+1}{p-1}$.
However, since for a solution of $\prl$ we have
$\rh=\rh_n=[\an+\lm_n\pn]_+^p$, by \eqref{Fl:sol} below we would have that,
$$
\inf \{2\mathcal{F}_{\lm_n}(\rh), \rho \in \mathcal{P}\}=2\mathcal{F}_{\lm_n}(\rh_n)=\frac{p-1}{p+1}\ino [\an+\lm_n \pn]_+^{p+1}+{\an}=
|\an|\left(\tau_n\frac{p-1}{p+1}-1\right)\to +\ii.
$$

On the other side, it readily follows from Lemma 2.2 in \cite{BMar2} (see in particular (2.6) in \cite{BMar2}) that \\
$\inf \{\mathcal{J}_{\lm_n}(u), u \in H\}\leq C$ for some positive constant~$C>0$. Remark that the functional used in \cite{BMar2}, say $\widetilde{\mathcal{J}}_{I}({\rm v})$,
deals directly with solutions of $\fbi$ and is related to ours as follows $\widetilde{\mathcal{J}}_{I}({\rm v})=I^{1+\frac1p}\mathcal{J}_{\lm}(u)$, where
${\rm v}=I^{\frac1p}u$. However,
by \eqref{equiv} below we would also have $\inf \{\mathcal{F}_{\lm_n}(\rh), \rho \in \mathcal{P}\}\leq C$, which is the desired
contradiction.  \finedim

\bigskip

At last we have,\\
{\bf The proof of Corollary \ref{altype3.not}.}\\
Let $(\lm_n,\pn)$ be any $\mathcal{Z}$-blow up sequence whose unique blow up point is an interior point, then
by \eqref{psiconc.1} we have,
$$
\frac{|\an|}{2\lm_n}\eps_n^{N-2}=\frac{1+o(1)}{M_{p,0}}.
$$
Consequently
\begin{align}
 \an+2\lm_n E_{\lm_n}
 =& \int\limits_{\om}[\an+\lm_n\pn]_+^{p+1}
    =|\an|\int\limits_{\om}|\an|^{p}[\vn-1]_+^{p}[\vn-1]_+ \nonumber \\
 =& |\an|\frac{1+o(1)}{M_{p,0}}\int\limits_{\om}\mu_n^{N/2}[\vn-1]_+^p[\vn-1]_+ \nonumber \\
 =& |\an|\frac{M_{p+1,0}}{M_{p,0}}(1+o(1)), \nonumber
\end{align}
that is,
$$
E_{\lm_n}=\frac{|\an|}{2\lm_n}(1+\frac{M_{p+1,0}}{M_{p,0}}+o(1)).
$$
Therefore we see that there is no upper bound for the energy in tha plasma region as in Lemma \ref{altype2.intro},
as indeed in view of \eqref{enrg1.0}, for this particular sequence we have,
$$
\int\limits_{\om_+}|\nabla \pn|^2=E_{\lm_n}-\frac{|\an|}{2\lm_n}=\frac{|\an|}{2\lm_n}(\frac{M_{p+1,0}}{M_{p,0}}+o(1)).
$$
On the other side, by Lemma \ref{altype3}, we see that,
$$
\frac{|\an|}{2\lm_n}(1+\frac{M_{p+1,0}}{M_{p,0}}+o(1))\leq
\frac{(p+1)}{4N^2 (\omega_{\sscp N})^{\frac 2 N}}\frac{\mu_{-,n}}{|\om_{+,n}|^{1-\frac{2}{N}}}+ \frac{|\an|}{2\lm_n},
$$
and in particular, since $|\om_{+,n}|=\omega_{\sscp N}(R_0\eps_n)^{N}(1+o(1))$,
\begin{align}
 \frac{|\an|}{2\lm_n}\eps_n^{N-2}(\frac{M_{p+1,0}}{M_{p,0}}+o(1))
 \leq & \frac{(p+1)}{4N^2 (\omega_{\sscp N})^{\frac 2 N}}\frac{\mu_{-,n}}{(\omega_{\sscp N}(R_0)^{N}(1+o(1)))^{1-\frac{2}{N}}} \nonumber \\
 =& \frac{(p+1)}{4N^2}\frac{R_0^2(1+o(1))}{(\omega_{\sscp N}(R_0)^{N})}\mu_{-,n}. \nonumber
\end{align}
Therefore, again by \eqref{psiconc.1},
$$
(\frac{M_{p+1,0}}{M^2_{p,0}}+o(1))\leq
\frac{(p+1)}{4N^2}\frac{R_0^2(1+o(1))}{(\omega_{\sscp N}(R_0)^{N})}\mu_{-,n},
$$
which gives
$$
\mu_{-,n}\geq\frac{4N^2\omega_{\sscp N}(R_0)^{N-2}}{(p+1)}(\frac{M_{p+1,0}}{M^2_{p,0}}+o(1)),
$$
as claimed.
\finedim

\bigskip

\section{Appendix A}\label{appendix}

The uniform bound ~\eqref{eq:assump-2.intro} is rather natural
as can be seen by a closer inspection of the underlying dual formulation arising from physical arguments.
In fact, solutions of $\prl$ can be found (\cite{BeBr}) as solutions of the dual variational problem of minimizing
free energy $\mathcal{F}_{\lm}(\rh)$,

\beq\label{var:1}
\inf \{\mathcal{F}_{\lm}(\rh), \rho \in \mathcal{P}\},\quad \mathcal{P}=
\left\{\rh\in L^{1+\frac1p}(\om)\,|\,\ino \rh=1,\;\; \rh\geq 0\;\mbox{a.e. in}\;\om \right\}
\eeq
$$
\mathcal{F}_{\lm}(\rh)=
{\scriptstyle \frac{p}{p+1}}\ino (\rh)^{1+\frac{1}{p}}-\frac\lm 2 \ino \rho G[\rho],
$$
where $(p+1)\mathcal{S}(\rh)={p}\ino (\rh)^{1+\frac{1}{p}}$ is
the so called nonextensive Entropy, which has been widely used in the description of space plasma physics, see \cite{Leub}, \cite{Liv}
and references quoted therein. Solutions of $\prl$ arising as minimizers of \eqref{var:1} are the so called
variational solutions (\cite{BMar2}). Remark that $\mathcal{P}$ is not the set of admissible functions used in \cite{BeBr}, however by using
a non-convex duality argument (see \cite{BeBr} p.421-422), it can be shown that in fact the variational principle \rife{var:1}
is equivalent to  \rife{var:2} below.
Denoting by $\alpha$ the Lagrange multiplier relative to the mass constraint, solutions of \eqref{var:1}
satisfy the Euler-Lagrange equation,
$$
\rh^{\frac1p}=\left[\alpha+\lm G[\rh] \right]_+, \quad G[\rh]=\ino G(x,y)\rh(y)dy,
$$
which is, putting $\psi=G[\rh]$, nothing but $\prl$.
As mentioned above, in fact \eqref{var:1} is equivalent to
\beq\label{var:2}
\inf \{\mathcal{J}_{\lm}(u), u \in H\},\quad {H}=
\left\{u\in H^{1}(\om)\,|\,\ino [u]_+^p=1,\;\;u=\mbox{constant on }\pa \om\right\}
\eeq
where
$$
\mathcal{J}_{\lm}(u)=\frac{1}{2\lm}\ino |\nabla u|^2-\frac{1}{p+1}\left(\ino [u]_+^{p+1}\right)+ u(\pa \om),
$$
which admits at least one minimizer for any $\lm>0$ (\cite{BeBr}). In particular, if $\rh_0$ is a minimizer of \eqref{var:1} for some $\lm_0>0$,
then there exists a unique $\al_0$ such that $u_0=\al_0+\lm_0 G[\rh_0]$ is a minimizer of \eqref{var:2} and
$\mathcal{F}_{\lm}(\rh_0)=\mathcal{J}_{\lm}(u_0)$. In other words the two problems are equivalent and provide
the same value of the minimum. Now for $\lm>0$ let $\rh_\lm$ be a minimizer of \eqref{var:1},
let $\pl=G[\rh_\lm]$ and define $\all$ such that $u_{\lm}=\all+\lm \pl$ is a minimizer of \eqref{var:2},
then we have that the common value of the minimum reads,
\beq\label{equiv}
\mathcal{J}_{\lm}(u_{\lm})=\mathcal{F}_{\lm}(\rh_\lm)=
\frac{p}{p+1}\ino [\al+\lm \pl]_+^{p+1}-\frac{\lm}{2} \ino|\nabla \pl|^2.
\eeq

We remark however that the second equality in \eqref{equiv} holds whenever $\rl$ is any critical
point of $\mathcal{F}_{\lm}$ and $(\all,\pl)$ is the corresponding solution of $\prl$.
At this point let us observe that, for any solution $(\all,\pl)$, we have (see the proof of Lemma \ref{enrg0} for details)
$$
{\lm} \ino|\nabla \pl|^2=\ino [\al+\lm \pl]_+^{p+1}-\all,
$$
whence putting $\frac1q=\frac{p-1}{p}$, we deduce that,
\beq\label{Fl:sol}
2\mathcal{F}_{\lm}(\rh_\lm)=\frac{p-1}{p+1}\ino [\all+\lm \pl]_+^{p+1}+{\all}\equiv
\frac{1}{q}\mathcal{S}(\rl)+{\all}.
\eeq
Therefore we see that the assumption $|\all|^{-1}\ino [\al+\lm \pl]_+^{p+1}\leq C_1$ just takes the form of
an upper bound for the entropy of the solutions in terms of $\all$,
$$
\mathcal{S}(\rl)\leq \frac{p}{p+1} C_1|\all|,
$$
and that this is just equivalent in turn to impose a control about the linear growth of the free energy in
terms of $\all$,
$$
2\mathcal{F}_{\lm}(\rh_\lm)\leq \frac{p-1}{p+1} C_1|\all|+\all.
$$

\bigskip

\section{Appendix B}
We generalize an argument worked out for $N=2$ in \cite{BSp} and obtain the uniqueness of solutions of $\fbi/\prl$
on balls ($N\geq 3$). Actually the fact that $p<p_{_N}$ plays a crucial role. Without loss of generality we work on balls of unit area $\mathbb{D}_N\subset \R^N$, $N\ge 3$,
whose radius is denoted by $R_N:=\Radius(\mathbb{D}_N)$.
Among other things, in view of the uniqueness, it is readily seen that \eqref{lmaldec} below supports our conjecture about \rife{hyp.intro.pn.1} 
for variational solutions, see Remark \ref{rem9.4}.

\

As in~\eqref{emden}, let~$u\in C^2(B_1(0))\cap C^1_0(\ov{B_1(0)})$ be the unique (\cite{gnn}) solution  of
\beq\label{emden-again}
\graf{-\Delta u =u^p\quad \mbox{in}\;\;B_1(0),\\ \\
u>0 \mbox{ in } B_1(0),\\ \\
u= 0 \mbox{ on } \pa B_1(0),}
\eeq
and recall that that~$u$ is radial and radially decreasing,~$u(x)= u(|x|)=u(r)$ with~$r=|x| \in (0,1)$, whence the equation reduces to
\begin{align}\nonumber
 u''(r)+\frac{N-1}{r}u'(r)= - u^p(r),\qquad r\in (0,1)
\end{align}
with~$u'(0)=0$,~$u(1)=0$. Here~$B_1(0)\subset \R^N$ is the ball of unit radius whose volume is $\omega_N$ and
from Remark~\ref{rem5.4} we have,
$$I_p=\int\limits_{B_{1}(0)}u^p = N\omega_N (-u'(1))>0. $$

For later convenience we introduce the notations,
\begin{align}\nonumber
 I^*(\mathbb{D}_N,p):=\frac{I_p}{R_N^{\frac{N}{p-1}(1-\frac{p}{p_{_N}})}}\equiv \frac{I_p}{R_N^{\frac{2p}{p-1}-N}} =\frac{I_p}{R_N^{\frac{2}{p-1}-N+2}}.
\end{align}
and
\begin{align}\nonumber
 \lm^*(\mathbb{D}_N,p):=I^*(\mathbb{D}_N,p)^{\frac{p-1}{p}}=\frac{I_p^{1-\frac{1}{p}}}{R_N^{\frac{N}{p}(1-\frac{p}{p_{_N}})}}.
\end{align}

\bte\label{thm9.2}
Let $\om=\mathbb{D}_N$ and $p\in (1,p_{_N})$, then for any $I>0$ {\rm (}$\lm>0${\rm )} there exists a unique solution of {\rm $\fbi$}
{\rm (}{\rm $\prl$}{\rm )}
denoted by~${\rm v}$ with boundary value~$\gamma\equiv {\rm v}|_{\pa \mathbb{D}}$.
Furthermore,
\begin{align}\nonumber
 \mbox{the solution is positive}
 &&\Leftrightarrow && \gamma>0
 &&\Leftrightarrow & &  I< I^*(\mathbb{D}_N, p)
 & & \Leftrightarrow && \lambda < \lambda^*(\mathbb{D}_N,p),
\end{align}
while
\begin{align}\nonumber
 \mbox{the free boundary is not empty}
 &&\Leftrightarrow && \gamma<0
 &&\Leftrightarrow & &  I> I^*(\mathbb{D}_N, p)
 & & \Leftrightarrow && \lambda > \lambda^*(\mathbb{D}_N,p),
\end{align}
and in terms of the variables~$(\all,\pl)$, we have
\beq\label{lmaldec}
\all=
c_N\lm \parenthesis{1 -
\parenthesis{\frac{\lambda}{\lambda^*(\mathbb{D}_N,p) } }^{\frac{p}{p_{_N} - p}} }, \qquad \forall \lambda> \lambda^*(\mathbb{D}_N,p),
\eeq
where $c_N=(N(N-2)\omega_N R_N^{N-2})^{-1}$.
\ete
\proof
Since solutions of $\prl$ with $\lm>0$ are in one to one correspondence with those of $\fbi$ with $I>0$ via $I^{p-1}=\lm^p$ and
\eqref{eq:change variable}, then it is enough to prove the uniqueness part of the statement for $\fbi$.\\
By the classical results in \cite{gnn} any solution of $\fbi$ in $\mathbb{D}_N$ is radial
so~${\rm v}_{\gamma}(x)={\rm v}_{\gamma}(r)$
satisfies
\begin{align}\nonumber
 {\rm v}_{\gamma}'' (r)+ \frac{N-1}{r}{\rm v}_{\gamma} = - [{\rm v}_{\gamma}]_+^p, \qquad r\in (0,R_N)
\end{align}
and~${\rm v}_{\gamma}'(0)=0$, ${\rm v}_{\gamma}(R_N)=\gamma$.
Let us define,
\begin{align}\nonumber
 {u}(|x|;R)
 =\begin{cases}
   \frac{1}{R^{\frac{2}{p-1}}} u(\frac{x}{R}) & 0\leq |x| \leq R , \\ \\
   \frac{A}{|x|^{N-2}} + B & R\leq |x| \leq R_N,
  \end{cases}
\end{align}
for some $R>0$ and in case $R<R_N$, observe that
\begin{align}\nonumber
 u(|x|;R)\in C^1
 & &\Leftrightarrow &  &
 \begin{cases}
  \frac{A}{R^{N-2}} + B =0 \\ \\
  -(N-2)\frac{A}{R^{N-1}}= \frac{1}{R^{\frac{2}{p-1}+1}} u'(1)
 \end{cases}
 & & \Leftrightarrow & &
 \begin{cases}
  A=\frac{-u'(1)}{(N-2)}\frac{1}{R^{\frac{2}{p-1}-N+2}}>0 \\ \\
  B=\frac{u'(1)}{(N-2)}\frac{1}{R^{\frac{2}{p-1}}} <0.
 \end{cases}
\end{align}
One readily checks that~$u \in C^2(\ov{B_{R_N}})$ which satisfies the same ODE as~${\rm v}_{_\gamma}$ does. Remark that it may happen that $R\geq R_N$,
in which case we set
${u}(|x|;R)=\frac{1}{R^{\frac{2}{p-1}}} u(\frac{x}{R}), 0\leq |x| \leq R$. Both ${\rm v}_{\gamma}(|x|)$ and $u(|x|;R)$ are radial and of class $C^2$ whence
$$
\frac{d}{dr}{\rm v}_{\gamma}(0)=0=\frac{d}{dr}u(0;R),
$$
and both satisfy $-\Delta {\rm v} =[{\rm v}]_+^p$ in their domain of definition. At this point, for $\gamma \in \R$ fixed, let us choose $R=R_\gamma$ such that,
$$
u(0;R_\gamma)=\frac{1}{R_{\gamma}^{\frac{2}{p-1}}} u(0)={\rm v}_{\gamma}(0).
$$
Obviously for any fixed ${\rm v}_{\gamma}(0)$, there exists one and only one $R_{\gamma}$ which fulfills this condition. Also, by ODE uniqueness theory,
${\rm v}_{\gamma}(|x|)\equiv u(|x|;R_\gamma)$ for $|x|\leq \mbox{\rm min}\{R_\gamma,R_N\}$.
However, by using the fact that ${\rm v}_{\gamma}(|x|)$ and ${u}(|x|;R)$ are strictly decreasing in their
domain of definition, it is easy to check that,
\beq\label{eq9.3}
R_\gamma>R_N \Leftrightarrow \gamma>0,\quad  R_\gamma=R_N \Leftrightarrow \gamma=0,\quad R_\gamma<R_N \Leftrightarrow \gamma<0.
\eeq

For any~$R>0$ let us define,
\begin{align}\nonumber
 \mathcal{I}(R)\coloneqq
 \int\limits_{B_{R_N}} [u(x;R)]_+^p\dx
 =\begin{cases}
   \frac{1}{R^{\frac{N}{p-1}(1-p/p_{_N})}} I_p, & R\leq R_N, \vspace{2mm}\\
   \frac{1}{R^{\frac{N}{p-1}(1-p/p_{_N})}} \int\limits_{B_{R_N/R}} u^p\dx,  & R\geq R_N.
  \end{cases}
\end{align}
Since~$p\in (1,p_{_N})$ then $\mathcal{I}(R)$ is strictly decreasing in~$R$ with range~$(0,+\infty)$.
Therefore for any $I>0$  there exists a unique $R(I)>0$ such that $\mathcal{I}(R(I))=I$.
At this point observe that if any such ${\rm v}_{\gamma}$ has to solve $\fbi$ for some $I>0$, then it must satisfy
$I=\int\limits_{\mathbb{D}_N} [{\rm v}_\gamma]_+^p =\int\limits_{B_{R_N}} [u(|x|;R_\gamma)]_+^p \dx$, that is, in view of \eqref{eq9.3},

\begin{align*}
\gamma\geq 0\Leftrightarrow R_\gamma\geq R_N\;\mbox{\rm and then }
&I=I(\gamma)=\mathcal{I}(R_\gamma)=\frac{1}{R_\gamma^{\frac{N}{p-1}(1-p/p_{_N})}} \int\limits_{B_{R_N/R_\gamma}} u^p\dx,
\end{align*}
\begin{align*}
\gamma< 0\Leftrightarrow R_\gamma< R_N\;\mbox{\rm and then }
&I=I(\gamma)=\mathcal{I}(R_\gamma)=\frac{1}{R_\gamma^{\frac{N}{p-1}(1-p/p_{_N})}}I_p.
\end{align*}

At last we argue by contradiction and assume that there exist $(\gamma_1,{\rm v}_1)\neq(\gamma_2,{\rm v}_2)$ sharing the same value of $I>0$.
Since $\mathcal{I}(R)$ is strictly decreasing,  they should share the same value of $R_{\gamma}$, i.e. $R_{\gamma_1}=R_{\gamma_2}$,
that is ${\rm v}_{1}(0)={\rm v}_{2}(0)$, implying by ODE uniqueness theory that ${\rm v}_1\equiv {\rm v}_2$ and
consequently $\gamma_1= \gamma_2$, which is a contradiction. This fact concludes the proof of the uniqueness.\\
Obviously if $(0,{\rm v}_0)$ solves $\fbi$, i.e. if $\gamma=0$, then necessarily ${\rm v}_0(r)\equiv u(r;R_N)$ in which case
$$
I=I^*(\mathbb{D}_N,p)=\mathcal{I}(R_N)=\frac{I_p}{R_N^{\frac{N}{p-1}(1-p/p_{_N})}}.
$$
In particular we have shown that for any $I$ the unique solution of $\fbi$ takes the form ${\rm v}_{\gamma}(r)=u(r;R_\gamma))$,
where $R_\gamma$ is uniquely defined by $I = \mathcal{I}(R_\gamma)$. Therefore in particular
\begin{align}\nonumber
 \mbox{the free boundary is not empty}
 &&\Leftrightarrow && \gamma<0
 &&\Leftrightarrow & &  I> I^*(\mathbb{D}_N, p)
 & & \Leftrightarrow && \lambda > \lambda^*(\mathbb{D}_N,p),
\end{align}
and similarly for the case $\gamma>0$, as claimed. In particular $\alpha$ is negative
if and only if $I>I^*(\mathbb{D}_N,p)$ which is the same as $R_\gamma<R_N$ and in this case we deduce from
\begin{align}\nonumber
 I=\mathcal{I}(R_\gamma)=\frac{1}{R_\gamma^{\frac{2}{p-1}-N+2}}I_p
\end{align}
that $R_{\gamma}$ satisfies,
\begin{align}\nonumber
 R_\gamma= \parenthesis{\frac{I_p}{I}}^{\frac{1}{\frac{2}{p-1}-N+2}}.
\end{align}
Thus the boundary value of ${\rm v}_{\gamma}$ takes the form,
\begin{align}\nonumber
 \gamma = &\frac{A}{R_N^{N-2}} + B
  = \frac{-u'(1)}{(N-2)} \parenthesis{ \frac{1}{R_N^{N-2}} \frac{1}{R_\gamma^{\frac{2}{p-1}-N+2}} -\frac{1}{R_\gamma^{\frac{2}{p-1}}} } \\
  =& \frac{I_p}{N(N-2)\omega_N}  \parenthesis{\frac{1}{R_N^{N-2}} \frac{I}{I_p} -
  \parenthesis{\frac{I}{I_p}}^{\frac{ \frac{2}{p-1} }{ \frac{2}{p-1} -N+2 }} } \nonumber \\
  =& \frac{1}{N(N-2)\omega_N} \parenthesis{ \frac{I}{R_N^{N-2}} -
  \parenthesis{\frac{I^{\frac{2}{p-1}}}{I_p^{N-2}}}^{\frac{1}{\frac{2}{p-1} -N+2 }} }, \nonumber
\end{align}
and then, by using~$I_p = I^*(\mathbb{D}_N,p) R_N^{\frac{2}{p-1}-N+2}$, we deduce that,
\begin{align}\nonumber
 \gamma = & \frac{1}{N(N-2)\omega_N} \frac{1}{R_N^{N-2}}\parenthesis{I- \parenthesis{\frac{I^{\frac{2}{p-1}}}{I^*(\mathbb{D}_N,p)^{N-2}}}^{\frac{1}{\frac{2}{p-1} -N+2 }} } \\
 = & \frac{1}{N(N-2)\omega_N R_N^{N-2}}\parenthesis{ I - I \parenthesis{\frac{I}{I^*(\mathbb{D}_N,p)}}^{\frac{p-1}{p_{_N} - p}} }, \qquad \forall I > I^*(\mathbb{D}_N, p). \nonumber
\end{align}

In terms of~$\alpha$ and~$\lambda$ this is the same as,
\begin{align} \nonumber
 \alpha=\frac{\gamma}{\lambda^{\frac{1}{p-1}}}
 =& \frac{1}{\lambda^{\frac{1}{p-1}}} \frac{1}{N(N-2)\omega_N R_N^{N-2}}\parenthesis{ \lambda^{\frac{p}{p-1}} - \lambda^{\frac{p}{p-1}} \parenthesis{\frac{\lambda^{\frac{p}{p-1}}}{\lambda^*(\mathbb{D}_N,p)^{\frac{p}{p-1}}}}^{\frac{p-1}{p_{_N} - p}} } \\
 =&\frac{1}{N(N-2)\omega_N R_N^{N-2}} \parenthesis{\lambda - \lambda \parenthesis{\frac{\lambda}{\lambda^*(\mathbb{D}_N,p) } }^{\frac{p}{p_{_N} - p}} },  \qquad \forall \lambda> \lambda^*(\mathbb{D}_N,p). \nonumber
\end{align}
which is~\eqref{lmaldec}.
\finedim

\brm\label{rem9.3} {\it
As far as we are concerned with $\gamma>0$, that is $I< I^*(\mathbb{D}_N, p)$ and $R_\gamma > R_N$,
the integral constraint takes the form,
\begin{align}\nonumber
 I= \mathcal{I}(R_\gamma) = \frac{1}{R_\gamma^{\frac{2}{p-1}-N+2}} \int\limits_{B_{R_N/ R_\gamma}} u^p \dx.
\end{align}
Of course as shown above there is a unique~$R_\gamma>0$ satisfying this relation, 
but since the function~$u$ solving~\eqref{emden-again} is not explicit, we cannot evaluate $R_\gamma$ explicitly.
Moreover, the boundary value of~${\rm v}_{\gamma}$ takes the form
\begin{align}
 \gamma = \frac{1}{R_{\gamma}^{\frac{2}{p-1}}} u\parenthesis{\frac{R_N}{R_{\gamma}}},
\end{align}
and in particular
\begin{itemize}
 \item $I\searrow 0 \Rightarrow R_\gamma\to +\infty$ and consequently $\gamma\to 0$,
 \item $I\nearrow I^*(\mathbb{D}_N,p)\Rightarrow R_\gamma \searrow R_N$ and consequently $\gamma \to 0$ (since~$u(1)=0$).
\end{itemize}
In other words $\gamma=\gamma(I)$ cannot be monotone in $I\in [0,I^*(\mathbb{D}_N,p)]$.
Since Theorem \ref{thm9.2} shows that we have a unique solution on $\mathbb{D}_N$ then we deduce from
Theorem A that in fact $\all$ is monotone decreasing, at least for $\lm<\frac1p\Lambda(\om,2p)$.}
\erm

\brm\label{rem9.4}
{\it
Clearly \eqref{lmaldec} implies that $\alpha\to -\infty$ iff~$\lambda \to +\infty$ and in particular that,
\begin{align}\nonumber
 \lim_{\lambda\to +\infty} \frac{|\alpha|}{\lambda^{\frac{p_{_N}}{p_{_N}- p}}}
 = \frac{1}{N(N-2)\omega_N R_N^{N-2}} \frac{1}{\lambda^*(\mathbb{D}_N,p)^{\frac{p}{p_{_N} -p}} }
 =\frac{1}{N(N-2)\omega_N}\frac{1}{I_p^{\frac{p-1}{p_{_N} - p}}},
\end{align}
whence
\begin{align}\nonumber
 \lim_{\lambda\to+\infty}\parenthesis{\frac{\lambda}{|\alpha|^{1-\frac{p}{p_{_N}}}}}^{\frac{N}{2}}
 =N\omega_N (N-2)^{\frac{N}{2}(1-\frac{p}{p_{_N}})} |u'(1)|^{(p-1)\frac{N-2}{2}}
 = M_{p,0},
\end{align}
which is just the mass quantization identity in this particular case.
}
\erm

\end{document}